\documentclass[a4paper,twoside]{article}

\usepackage{amssymb,amsmath,amsthm,dsfont,amsfonts,color,latexsym}
\allowdisplaybreaks[4]
\numberwithin{equation}{section}
\theoremstyle{plain}
\newtheorem{thm}{Theorem}[section]
\newtheorem{theorem}{Theorem}[section]      
\newtheorem{lemma}[theorem]{Lemma}          
\newtheorem{proposition}[theorem]{Proposition}
\newtheorem{corollary}[theorem]{Corollary}

\theoremstyle{definition}

\theoremstyle{definition}
\newtheorem{remark}[theorem]{Remark}

\usepackage{hyperref}
\usepackage{mathrsfs} 

\definecolor{blue}{rgb}{0.00,0.00,1.00}
\definecolor{red}{rgb}{1.00,0.00,0.00}

\renewcommand{\baselinestretch}{1.2}
\begin{document}


\title{Non-Isothermal Plane Couette Flow for the Boltzmann Equation in Three Dimensions}
\author{
  Zhida Chang,\, Mingying Zhong\\[2mm]
  {\small\itshape School of Mathematics, Guangxi University, China.}\\
  {\small\itshape E-mails: czdzhidachang@163.com, zhongmingying@gxu.edu.cn}\\[5mm]
  }
\date{ }

\pagestyle{myheadings}
\markboth{Non-Isothermal Plane Couette Flow}{Z.-D.~Chang, M.-Y.~Zhong}

\maketitle

\thispagestyle{empty}

\begin{abstract}
\noindent{We consider the Boltzmann equation for non-isothermal plane Couette flow in the three-dimensional channel $\mathbb T_x\times(-1,1)_y\times\mathbb T_z$, with diffuse reflection and cutoff hard and soft potentials, $-3<\gamma\le1$. We develop a new boundary flux transmission framework that reduces the diffuse boundary feedback to closed equations for the outgoing mass fluxes. Strict contraction of the associated boundary-to-boundary transmission operators yields direct representations and estimates for the outgoing mass fluxes in terms of interior sources, boundary sources, and initial data and substantially simplifies the analysis of diffuse reflection. For damped transport equations with diffuse reflection in general bounded $C^1$ domains, we establish this strict contraction property and hence the unique solvability of the boundary flux equations in $L^p$ for $1\le p\le\infty$, even when the damping degenerates at large velocities. For the Couette problem, we construct a unique steady solution near the global Maxwellian when the shear rate and wall temperature perturbations are sufficiently small. We also establish global existence and uniqueness of unsteady solutions for small initial perturbations with polynomial decay in velocity, with exponential decay for $0\le\gamma\le1$ and algebraic decay for $-3<\gamma<0$.}


\end{abstract}

\tableofcontents  
\section{Introduction}
In this paper, we consider a non-isothermal plane Couette flow in the three-dimensional channel
\[
\Omega:=\mathbb T_x\times(-1,1)_y\times\mathbb T_z,
\quad
\partial\Omega:=\mathbb T_x\times\{-1,1\}\times\mathbb T_z,
\]
where $\mathbb T:=\mathbb R/\mathbb Z$. The upper and lower walls $y=\pm1$ move in the $x$-direction with respective velocities $U(\pm1)=\pm\alpha$ and are maintained at temperatures $T(x,\pm1,z)=1+\delta\vartheta(x,\pm1,z)$, where $\vartheta(x,\pm1,z)$ are prescribed functions on $\mathbb T_x\times\mathbb T_z$. Here $\alpha>0$ and $\delta>0$ denote the shear rate and the temperature perturbation parameter, respectively. Throughout the paper, we assume $\sup_{(x,z)\in\mathbb T^2}|\vartheta(x,\pm1,z)|<\infty$.

Let $v=(v_x,v_y,v_z)\in\mathbb R^3$ denote the physical velocity, and let $F_{\alpha,\delta}=F_{\alpha,\delta}(\mathbf x,v)$, with $\mathbf x=(x,y,z)\in\Omega$, be the steady distribution function of the gas particles. The normalized steady problem is
\begin{equation}\label{eq-Boltz-st}
\left\{
\begin{aligned}
&v\cdot\nabla_{\mathbf x}F_{\alpha,\delta}
=Q(F_{\alpha,\delta},F_{\alpha,\delta}),
\\
&F_{\alpha,\delta}(x,\pm1,z,v)\big|_{\pm v_y<0}
=\mathcal M_{T(x,\pm1,z)}\big(v_x-U(\pm1),v_y,v_z\big)
\int_{\pm u_y>0}F_{\alpha,\delta}(x,\pm1,z,u)|u_y|\,du,
\\
&\frac12\int_\Omega\int_{\mathbb R^3}F_{\alpha,\delta}(\mathbf x,v)\,dv d\mathbf x=1.
\end{aligned}
\right.
\end{equation}
Here
\begin{equation}\label{def-wall-Maxw}
\mathcal M_T(v)
:=
\frac{1}{2\pi T^2}e^{-\frac{|v_x|^2+|v_y|^2+|v_z|^2}{2T}}
\end{equation}
is the boundary Maxwellian normalized with respect to the normal mass flux. The bilinear Boltzmann collision operator is defined by
\begin{align}
Q(F,G)(v)
&=
\int_{\mathbb R^3}\int_{\mathbb S^2}
|v-v_*|^\gamma B_0(\cos\theta)
\big[F(v_*')G(v')-F(v_*)G(v)\big]
\,d\omega dv_*
\notag\\
&=:
Q_+(F,G)(v)-Q_-(F,G)(v),
\label{eq-coll-op}
\end{align}
where
\begin{equation}\label{eq-post-vel}
\cos\theta
=
\frac{v-v_*}{|v-v_*|}\cdot\omega,
\quad
v'=v-\big[(v-v_*)\cdot\omega\big]\omega,
\quad
v_*'=v_*+\big[(v-v_*)\cdot\omega\big]\omega.
\end{equation}
The pre- and post-collisional velocities obey the conservation laws
\begin{equation}\label{eq-coll-cons}
v+v_*=v'+v_*',
\quad
|v|^2+|v_*|^2=|v'|^2+|v_*'|^2.
\end{equation}
Throughout the paper, we consider hard and soft potentials with $-3<\gamma\le1$ under Grad's angular cutoff assumption and suppose that
\begin{equation}\label{eq-ang-cutoff}
0\le B_0(\cos\theta)\le C|\cos\theta|.
\end{equation}

The kinetic study of Couette flow is closely connected with the earlier theory of homogeneous shear and homoenergetic flows. For Maxwell molecules, the moment structure associated with shear deformation was already investigated in the classical works of Truesdell and Galkin~\cite{G58,T56}, while Cercignani later established the existence of homoenergetic affine solutions for the Boltzmann equation~\cite{Cer89}. More recent studies addressed non-Maxwellian self-similar profiles under shear and more general deformations~\cite{BNV20,JNV19}. A rigorous nonlinear theory of uniform shear flow was developed by Duan--Liu for Maxwell molecules~\cite{DL21} and subsequently extended to hard potentials~\cite{DL25}. These homogeneous or homoenergetic models capture the kinetic effects of shear without the physical boundaries intrinsic to the Couette problem.

Plane Couette flow, driven by the relative tangential motion of two parallel walls, couples shear deformation with kinetic boundary effects. Early kinetic studies based on moment methods, linearized models, boundary-layer expansions, and numerical computations revealed characteristic phenomena such as velocity slip, temperature variation, and Knudsen-layer corrections~\cite{CC94,GTRTS06,GZ58,LL61,MSG00,Rog16,STO90,SY70,Wil62,YL57} For the full Boltzmann equation in plane geometry, Duan--Liu--Yang~\cite{DLY25} developed a rigorous theory of isothermal plane Couette flow with diffuse reflection for Maxwell molecules. For sufficiently small wall velocity, they constructed a non-equilibrium steady state with polynomially weighted estimates and established its exponential stability under small perturbations independent of the tangential variables. Subsequent studies have addressed plane Couette flow with hard potentials in a finite channel~\cite{MW24}, three-dimensional kinetic Couette flow in the diffusive limit~\cite{DLSY24}, and related multidimensional shear-driven kinetic problems~\cite{CZ26,DLS25,DWW25}.

Non-isothermal boundary forcing constitutes another closely related line of research. Temperature differences between the walls were incorporated into several formal and numerical studies of plane and cylindrical Couette flow~\cite{GTRTS06,LL61,SK99,YL57}. On the rigorous side, Esposito--Lebowitz--Marra studied the stationary Boltzmann equation in a slab with equal wall temperatures and a constant tangential external force, and justified its hydrodynamic limit~\cite{ELM94}. Their subsequent work incorporated different wall temperatures and relative wall motion, together with a tangential external force, and justified the corresponding stationary Navier--Stokes limit~\cite{ELM95}. At the kinetic level, Esposito--Guo--Kim--Marra~\cite{EGKM13} established the existence, uniqueness, and dynamical stability of a non-equilibrium steady state generated by small variations in wall temperature under diffuse reflection in a general bounded domain. Further developments concerning non-isothermal boundaries, long-time behavior, heat transfer, and regularity can be found in~\cite{CK24,DHWZ19,DLYZ22}.

The closest precedents for the present problem are~\cite{DLY25,EGKM13}, which establish the existence and dynamical stability of non-equilibrium steady states for isothermal Couette flow and non-isothermal boundary forcing without Couette shear, respectively. Relative wall motion and different wall temperatures were considered simultaneously in~\cite{ELM95}, where steady solutions were constructed in the small-Knudsen-number regime to justify the stationary Navier--Stokes limit. Motivated by these works, we study steady states and their dynamical stability at a fixed Knudsen number under the combined action of Couette shear and non-isothermal diffuse boundary forcing.

\noindent\textbf{Notation.}
Throughout the paper, $C>0$ denotes a generic constant that may vary from line to line and is independent of the parameters under consideration unless otherwise specified. We write $A\lesssim B$ if $A\le CB$. Constants depending on specific parameters are indicated by subscripts, such as $C_l$. The notation $\varepsilon>0$ denotes a sufficiently small constant.

For $p=2,\infty$, we use $\|\cdot\|_{L^p_{\mathbf x,v}}$,
$\|\cdot\|_{L^p_{\mathbf x}}$, $\|\cdot\|_{L^p_{y,v}}$,
$\|\cdot\|_{L^p_v}$, and $\|\cdot\|_{H^2_{\mathbf x}}$ for the standard norms on
$L^p(\Omega\times\mathbb R_v^3)$, $L^p(\Omega)$,
$L^p((-1,1)_y\times\mathbb R_v^3)$, $L^p(\mathbb R_v^3)$,
and $H^2(\Omega)$, respectively. We also write
\[
\begin{aligned}
&\|f\|_{L^\infty_{x,z}}
:=
\operatorname*{ess\,sup}_{(x,z)\in\mathbb T^2}|f(x,z)|,
\quad
\|f\|_{L^\infty(\partial\Omega)}
:=
\operatorname*{ess\,sup}_{\mathbf x\in\partial\Omega}|f(\mathbf x)|,
\\
&\|f\|_{L^2(\partial\Omega)}^2
:=
\int_{\partial\Omega}|f(\mathbf x)|^2\,dS_{\mathbf x}
=
\sum_\pm\int_{\mathbb T^2}|f(x,\pm1,z)|^2\,dxdz.
\end{aligned}
\]
For the $L^2$ spaces, the corresponding inner products are denoted by
\[
\begin{aligned}
&\langle f,g\rangle_{\mathbf x,v}
:=
\int_\Omega\int_{\mathbb R^3}f(\mathbf x,v)g(\mathbf x,v)\,dv d\mathbf x,
\\
&\langle f,g\rangle_{\mathbf x}
:=
\int_\Omega f(\mathbf x)g(\mathbf x)\,d\mathbf x,
\quad
\langle f,g\rangle_v
:=
\int_{\mathbb R^3}f(v)g(v)\,dv.
\end{aligned}
\]

For $y=\pm1$, let $n(1)=(0,1,0)$ and $n(-1)=(0,-1,0)$. The phase boundary is denoted by $\gamma:=\partial\Omega\times\mathbb R_v^3$ and is decomposed into the outgoing, incoming, and grazing sets:
\[
\left\{
\begin{aligned}
\gamma_+
&:=\big\{(\mathbf x,v)\in\gamma:n(\mathbf x)\cdot v>0\big\}
=\big\{(x,\pm1,z,v):x,z\in\mathbb T,\ \pm v_y>0\big\},
\\
\gamma_-
&:=\big\{(\mathbf x,v)\in\gamma:n(\mathbf x)\cdot v<0\big\}
=\big\{(x,\pm1,z,v):x,z\in\mathbb T,\ \pm v_y<0\big\},
\\
\gamma_0
&:=\big\{(\mathbf x,v)\in\gamma:n(\mathbf x)\cdot v=0\big\}
=\big\{(x,\pm1,z,v):x,z\in\mathbb T,\ v_y=0\big\}.
\end{aligned}
\right.
\]
The kinetic boundary measure on $\gamma$ is $d\gamma:=|n(\mathbf x)\cdot v|\,dv dS_{\mathbf x}=|v_y|\,dvdxdz$. We use the boundary norms
\[
|f|_{L^2(\gamma_\pm)}^2
:=
\int_{\gamma_\pm}|f|^2\,d\gamma,
\quad
|f|_{L^\infty(\gamma_\pm)}
:=
\operatorname*{ess\,sup}_{(\mathbf x,v)\in\gamma_\pm}|f(\mathbf x,v)|.
\]

Let
\begin{equation}\label{def-global-max}
\mu(v)
:=
\frac1{(2\pi)^{3/2}}e^{-\frac{|v|^2}{2}}
\end{equation}
be the global Maxwellian. We define the outgoing mass flux operator $\mathcal J:L^\infty(\gamma_+)\to L^\infty(\partial\Omega)$ by
\begin{equation}\label{def-outgoing-flux}
\mathcal J[f](\mathbf x)
:=
\int_{n(\mathbf x)\cdot u>0}
f(\mathbf x,u)\sqrt\mu\big[n(\mathbf x)\cdot u\big]\,du,
\quad
\mathbf x\in\partial\Omega.
\end{equation}
In particular, at the two walls $y=\pm1$,
\[
\mathcal J[f](x,\pm1,z)
=
\int_{\pm u_y>0}
f(x,\pm1,z,u)\sqrt\mu|u_y|\,du.
\]
The diffuse reflection projection is defined by
\[
P_\gamma f(x,\pm1,z,v)
:=
\sqrt{2\pi}\sqrt\mu\mathcal J[f](x,\pm1,z).
\]
By the normalization of $\mu$, $\mathcal J[P_\gamma f]=\mathcal J[f]$ and $P_\gamma^2=P_\gamma$. Moreover, the restriction of $P_\gamma$ to $\gamma_+$ is an orthogonal projection in $L^2(\gamma_+)$, and hence
\begin{equation}\label{eq-f-L2gamma+}
|f|_{L^2(\gamma_+)}^2
=
|P_\gamma f|_{L^2(\gamma_+)}^2
+
|(I-P_\gamma)f|_{L^2(\gamma_+)}^2.
\end{equation}

We first construct a steady solution $F_{\alpha,\delta}^s=F_{\alpha,\delta}^s(\mathbf x,v)$ to \eqref{eq-Boltz-st} near the global Maxwellian for $(\alpha,\delta)$ close to $(0,0)$, and derive weighted estimates uniform in $\alpha$ and $\delta$. We seek $F_{\alpha,\delta}^s$ in the form
\begin{equation}\label{eq-st-exp}
F_{\alpha,\delta}^s
=
\mu+\sqrt\mu\Big[
\alpha G_1^\alpha+\delta G_1^\delta
+(\alpha^2+\delta^2)G_R^{\alpha,\delta}
\Big],
\end{equation}
where $G_1^\alpha=G_1^\alpha(\mathbf x,v)$ and $G_1^\delta=G_1^\delta(\mathbf x,v)$ are the first-order shear correction and the first-order thermal correction, respectively, and $G_R^{\alpha,\delta}=G_R^{\alpha,\delta}(\mathbf x,v)$ is the higher-order remainder. Accordingly, we impose the zero-mass constraints
\begin{equation}\label{eq-mass-exp}
\int_\Omega\int_{\mathbb R^3}\sqrt\mu G_1^\alpha\,dv d\mathbf x
=
\int_\Omega\int_{\mathbb R^3}\sqrt\mu G_1^\delta\,dv d\mathbf x
=
\int_\Omega\int_{\mathbb R^3}\sqrt\mu G_R^{\alpha,\delta}\,dv d\mathbf x
=
0.
\end{equation}
To derive the boundary conditions associated with the expansion \eqref{eq-st-exp}, we expand the boundary Maxwellian as
\begin{align}
&\mathcal M_{T(x,\pm1,z)}\big(v_x-U(\pm1),v_y,v_z\big)
\notag\\
&=
\mathcal M_1(v)
+
\sqrt{2\pi}\sqrt\mu\Big[
\alpha\eta_1^\alpha(\pm1,v)
+\delta\eta_1^\delta(x,\pm1,z,v)
+(\alpha^2+\delta^2)\eta_R^{\alpha,\delta}(x,\pm1,z,v)
\Big],
\label{eq-MT-exp}
\end{align}
where
\begin{equation}\label{def-eta1}
\eta_1^\alpha(\pm1,v)
:=
\pm v_x\sqrt\mu,
\quad
\eta_1^\delta(x,\pm1,z,v)
:=
\vartheta(x,\pm1,z)\bigg(\frac{|v|^2}{2}-2\bigg)\sqrt\mu
\end{equation}
are the first-order shear boundary correction and the first-order thermal boundary correction, respectively, while $\eta_R^{\alpha,\delta}(x,\pm1,z,v)$ is the higher-order boundary remainder defined through \eqref{eq-MT-exp}.

Substituting the expansion \eqref{eq-st-exp} into \eqref{eq-Boltz-st} and collecting the terms of order $\alpha$, we obtain the following boundary value problem for the first-order shear correction $G_1^\alpha$:
\begin{equation}\label{eq-Ga}
\left\{
\begin{aligned}
&v\cdot\nabla_{\mathbf x}G_1^\alpha+LG_1^\alpha
=
0,
\quad (\mathbf x,v)\in\Omega\times\mathbb R^3,
\\
&G_1^\alpha(x,\pm1,z,v)\big|_{\pm v_y<0}
=
P_\gamma G_1^\alpha(x,\pm1,z,v)
+\eta_1^\alpha(\pm1,v),
\\
&\int_\Omega\int_{\mathbb R^3}
\sqrt\mu G_1^\alpha(\mathbf x,v)\,dv d\mathbf x=0.
\end{aligned}
\right.
\end{equation}
Similarly, collecting the terms of order $\delta$, we obtain the boundary value problem for the first-order thermal correction $G_1^\delta$:
\begin{equation}\label{eq-Gd}
\left\{
\begin{aligned}
&v\cdot\nabla_{\mathbf x}G_1^\delta+LG_1^\delta
=
0,
\quad (\mathbf x,v)\in\Omega\times\mathbb R^3,
\\
&G_1^\delta(x,\pm1,z,v)\big|_{\pm v_y<0}
=
P_\gamma G_1^\delta(x,\pm1,z,v)
+\eta_1^\delta(x,\pm1,z,v),
\\
&\int_\Omega\int_{\mathbb R^3}
\sqrt\mu G_1^\delta(\mathbf x,v)\,dv d\mathbf x=0.
\end{aligned}
\right.
\end{equation}
After subtracting the first-order contributions, the higher-order remainder $G_R^{\alpha,\delta}$ satisfies
\begin{equation}\label{eq-GR}
\left\{
\begin{aligned}
&v\cdot\nabla_{\mathbf x}G_R^{\alpha,\delta}
+LG_R^{\alpha,\delta}
\\
&\quad=
\frac{\alpha^2}{\alpha^2+\delta^2}
\Gamma(G_1^\alpha,G_1^\alpha)
+\frac{\alpha\delta}{\alpha^2+\delta^2}
\Big[
\Gamma(G_1^\alpha,G_1^\delta)
+\Gamma(G_1^\delta,G_1^\alpha)
\Big]
\\
&\qquad+
\frac{\delta^2}{\alpha^2+\delta^2}
\Gamma(G_1^\delta,G_1^\delta)
+\alpha
\Big[
\Gamma(G_1^\alpha,G_R^{\alpha,\delta})
+\Gamma(G_R^{\alpha,\delta},G_1^\alpha)
\Big]
\\
&\qquad+
\delta
\Big[
\Gamma(G_1^\delta,G_R^{\alpha,\delta})
+\Gamma(G_R^{\alpha,\delta},G_1^\delta)
\Big]
+(\alpha^2+\delta^2)
\Gamma(G_R^{\alpha,\delta},G_R^{\alpha,\delta}),
\quad (\mathbf x,v)\in\Omega\times\mathbb R^3,
\\
&G_R^{\alpha,\delta}(x,\pm1,z,v)\big|_{\pm v_y<0}
=
P_\gamma G_R^{\alpha,\delta}(x,\pm1,z,v)
+\eta_R^{\alpha,\delta}(x,\pm1,z,v)
\\
&\qquad+
\sqrt{2\pi}
\Big[
\alpha\eta_1^\alpha(\pm1,v)
+\delta\eta_1^\delta(x,\pm1,z,v)
+(\alpha^2+\delta^2)\eta_R^{\alpha,\delta}(x,\pm1,z,v)
\Big]
\mathcal J[G_R^{\alpha,\delta}](x,\pm1,z)
\\
&\qquad+
\frac{\sqrt{2\pi}}{\alpha^2+\delta^2}
\Big[
\alpha\eta_1^\alpha(\pm1,v)
+\delta\eta_1^\delta(x,\pm1,z,v)
+(\alpha^2+\delta^2)\eta_R^{\alpha,\delta}(x,\pm1,z,v)
\Big]
\\
&\quad\quad\quad
\times\Big[
\alpha\mathcal J[G_1^\alpha](x,\pm1,z)
+\delta\mathcal J[G_1^\delta](x,\pm1,z)
\Big],
\\
&\int_\Omega\int_{\mathbb R^3}
\sqrt\mu G_R^{\alpha,\delta}(\mathbf x,v)\,dv d\mathbf x=0.
\end{aligned}
\right.
\end{equation}
Here, the linearized collision operator $L$ and the bilinear collision operator $\Gamma$ are defined by
\begin{equation}\label{def-LG}
Lf:=-\mu^{-\frac12}
\big[Q(\mu,\sqrt\mu f)+Q(\sqrt\mu f,\mu)\big],
\quad
\Gamma(f,g):=\mu^{-\frac12}Q(\sqrt\mu f,\sqrt\mu g).
\end{equation}

For $l\ge0$, let $w_l(v):=(1+|v|^2)^l$ denote the polynomial velocity weight. We now state the first result of this paper.

\begin{thm}\label{thm-1}
For $-3<\gamma\le1$ and any $l>1$, there exist $\alpha_0=\alpha_0(l)>0$ and $\delta_0=\delta_0(l)>0$ such that, whenever $0<\alpha\le\alpha_0$ and $0<\delta\le\delta_0$, the three-dimensional steady problem \eqref{eq-Boltz-st} admits a unique steady solution
$F_{\alpha,\delta}^s=F_{\alpha,\delta}^s(\mathbf x,v)$ of the form \eqref{eq-st-exp}.
\begin{enumerate}
\item[(a)]
The first-order shear correction $G_1^\alpha=G_1^\alpha(y,v)$ is the unique solution to \eqref{eq-Ga} independent of $x,z$. It satisfies
\begin{equation}\label{eq-st-Ga-sy}
\left\{
\begin{aligned}
&G_1^\alpha(y,-v_x,v_y,v_z)=-G_1^\alpha(y,v_x,v_y,v_z),
\\
&G_1^\alpha(y,v_x,v_y,-v_z)=G_1^\alpha(y,v_x,v_y,v_z),
\end{aligned}
\right.
\end{equation}
and
\begin{equation}\label{est-Linf-st-Ga}
\big\|w_lG_1^\alpha\big\|_{L^\infty_{y,v}}
\le C_l.
\end{equation}

\item[(b)]
The first-order thermal correction $G_1^\delta=G_1^\delta(\mathbf x,v)$ is the unique solution to \eqref{eq-Gd}. It satisfies
\begin{equation}\label{est-Linf-st-Gd}
\big\|w_lG_1^\delta\big\|_{L^\infty_{\mathbf x,v}}
\le C_l.
\end{equation}

\item[(c)]
The higher-order remainder $G_R^{\alpha,\delta}=G_R^{\alpha,\delta}(\mathbf x,v)$ is the unique solution to \eqref{eq-GR} satisfying
\begin{equation}\label{est-Linf-st-GR}
\big\|w_lG_R^{\alpha,\delta}\big\|_{L^\infty_{\mathbf x,v}}
\le C_l.
\end{equation}
\end{enumerate}
Here $C_l>0$ is independent of $\alpha$ and $\delta$.
\end{thm}

\begin{remark}
We make the following remarks on Theorem~\ref{thm-1}.

\noindent{\rm (i)} The estimates \eqref{est-Linf-st-Ga}--\eqref{est-Linf-st-GR} imply
\[
\big|F_{\alpha,\delta}^s(\mathbf x,v)-\mu\big|
\le
C_l(\alpha+\delta)w_l^{-1}\sqrt{\mu},\qquad \forall l>1.
\]
Consequently, the steady solution $F_{\alpha,\delta}^s$ exhibits Gaussian decay in velocity for both hard and soft potentials. Moreover, \eqref{eq-st-Ga-sy} implies that the first-order shear correction contributes neither to the density nor to the energy:
\[
\int_{\mathbb R^3}\sqrt\mu G_1^\alpha\,dv=\int_{\mathbb R^3}|v|^2\sqrt\mu G_1^\alpha\,dv=0.
\]

\noindent{\rm (ii)} The one-dimensional structure of $G_1^\alpha=G_1^\alpha(y,v)$ follows directly from tangential translational invariance and uniqueness. If $\vartheta(x,\pm1,z)=\vartheta(\pm1)$, the full steady-state problem is invariant under tangential translations, and uniqueness implies that $G_1^\delta$ and $G_R^{\alpha,\delta}$ are also independent of $x$ and $z$. Consequently, $F_{\alpha,\delta}^s(\mathbf x,v)=F_{\alpha,\delta}^s(y,v)$. Thus, Theorem~\ref{thm-1} yields a one-dimensional steady solution within the full three-dimensional setting.

\noindent{\rm (iii)} To the best of our knowledge, Theorem~\ref{thm-1} provides the first perturbative existence and uniqueness result for the full three-dimensional steady Boltzmann problem describing a non-isothermal plane Couette flow with diffuse reflection for hard and soft potentials. When $\delta=0$, the problem reduces to the isothermal plane Couette flow, for which the one-dimensional steady theory for Maxwell molecules was established in~\cite{DLY25}. When $\alpha=0$, the problem reduces to the steady Boltzmann equation driven by a non-isothermal diffuse reflection boundary and is closely related to the existing theory for non-isothermal boundary problems~\cite{EGKM13}. When $\alpha=\delta=0$, both shear and thermal driving vanish, and the problem \eqref{eq-Boltz-st} reduces to the classical equilibrium steady Boltzmann problem, with the global Maxwellian $\mu$ as a steady solution.
\end{remark}

Having constructed the steady state in Theorem~\ref{thm-1}, we turn to the following unsteady problem:
\begin{equation}\label{eq-prob-3d}
\left\{
\begin{aligned}
&\partial_tF
+v\cdot\nabla_{\mathbf x}F
=
Q(F,F),
\quad
(t,\mathbf x,v)\in\mathbb R_+\times\Omega\times\mathbb R^3,
\\
&F(t,x,\pm1,z,v)\big|_{\pm v_y<0}
=
\mathcal M_{T(x,\pm1,z)}\big(v_x-U(\pm1),v_y,v_z\big)
\int_{\pm u_y>0}
F(t,x,\pm1,z,u)|u_y|\,du,
\\
&F(0,\mathbf x,v)=F_0(\mathbf x,v),
\quad
\frac12
\int_\Omega\int_{\mathbb R^3}
F(t,\mathbf x,v)\,dv d\mathbf x
=
1.
\end{aligned}
\right.
\end{equation}

\begin{thm}\label{thm-2}
Let $F_{\alpha,\delta}^s(\mathbf x,v)$ be the steady state obtained in Theorem~\ref{thm-1}. There exists $l_0>1$ such that, for any $l\ge l_0$, there exist constants $\alpha_0=\alpha_0(l)>0$, $\delta_0=\delta_0(l)>0$, and $\varepsilon_0=\varepsilon_0(l)>0$ such that, whenever $0<\alpha\le\alpha_0$ and $0<\delta\le\delta_0$, if the initial datum $F_0(\mathbf x,v)\ge0$ satisfies
\[
\big\|w_l\big[F_0(\mathbf x,v)-F_{\alpha,\delta}^s(\mathbf x,v)\big]\big\|_{L^\infty_{\mathbf x,v}}
\le
\varepsilon_0
\]
and
\begin{equation}\label{eq-t0-datum}
\int_\Omega\int_{\mathbb R^3}\big[F_0(\mathbf x,v)-F_{\alpha,\delta}^s(\mathbf x,v)\big]\,dv d\mathbf x=0,
\end{equation}
then the corresponding unsteady boundary value problem \eqref{eq-prob-3d} admits a unique global solution $F(t,\mathbf x,v)\ge0$. Moreover,
\begin{enumerate}
\item[(i)]
For $0\le\gamma\le1$, there exists $\lambda>0$ such that
\begin{equation}\label{est-global-decay}
\big\|w_l\big[F(t,\mathbf x,v)-F_{\alpha,\delta}^s(\mathbf x,v)\big]\big\|_{L^\infty_{\mathbf x,v}}
\le
C_le^{-\lambda t}\big\|w_l\big[F_0(\mathbf x,v)-F_{\alpha,\delta}^s(\mathbf x,v)\big]\big\|_{L^\infty_{\mathbf x,v}}.
\end{equation}

\item[(ii)]
For $-3<\gamma<0$, we have
\begin{equation}\label{est-global-bound-soft}
\big\|w_l\big[F(t,\mathbf x,v)-F_{\alpha,\delta}^s(\mathbf x,v)\big]\big\|_{L^\infty_{\mathbf x,v}}
\le
C_l\big\|w_l\big[F_0(\mathbf x,v)-F_{\alpha,\delta}^s(\mathbf x,v)\big]\big\|_{L^\infty_{\mathbf x,v}}.
\end{equation}
Furthermore, if $l\ge l_0+|\gamma|/2$, then for any $m>0$ satisfying $l-m\ge l_0$,
\begin{equation}\label{est-global-decay-soft}
\big\|w_{l-m}\big[F(t,\mathbf x,v)-F_{\alpha,\delta}^s(\mathbf x,v)\big]\big\|_{L^\infty_{\mathbf x,v}}
\le
C_l(1+t)^{-\frac{2m}{|\gamma|}}
\big\|w_l\big[F_0(\mathbf x,v)-F_{\alpha,\delta}^s(\mathbf x,v)\big]\big\|_{L^\infty_{\mathbf x,v}}.
\end{equation}
\end{enumerate}
Here $\lambda>0$ and $C_l>0$ are independent of $\alpha$ and $\delta$.
\end{thm}

\begin{remark}
We make the following remarks on Theorem~\ref{thm-2}.

\noindent{\rm (i)}
The estimates \eqref{est-global-bound-soft} and \eqref{est-global-decay-soft}
show that, for soft potentials, the perturbation remains uniformly close
to its initial value in the original weighted norm and decays at rate
$(1+t)^{-2m/|\gamma|}$ after any admissible loss $m>0$ of velocity weight.
With the lower velocity weight $w_{l-m}$ fixed, increasing the exponent $l$
of the initial weight $w_l$ enlarges the gap $m$ and hence increases
the algebraic decay exponent $2m/|\gamma|$.
This relation suggests investigating in future work whether exponential
velocity tails of the initial perturbation yield nearly exponential time decay.

\noindent{\rm (ii)}
Theorem~\ref{thm-2} also yields the non-negativity of the steady solution obtained in Theorem~\ref{thm-1}. Indeed, choosing a non-negative initial datum sufficiently close to $F_{\alpha,\delta}^s$, the corresponding global solution remains non-negative and converges to $F_{\alpha,\delta}^s$ as $t\to\infty$, which implies $F_{\alpha,\delta}^s\ge0$
\end{remark}

We now outline the main ideas and strategies of the proof. In the $L^\infty$ analysis, a standard approach to handling diffuse reflection is the Stochastic Cycles developed in~\cite{Guo10}, which iterates the boundary condition along successive backward characteristics.
In the present work, we develop a new boundary flux transmission framework
for both steady and time-dependent problems. Its key ingredient is the strict contraction of the associated boundary-to-boundary transmission operators. Under diffuse reflection, the incoming trace can be written schematically as
\[
h|_{\gamma_-}
=
\underbrace{\Phi_{\mathrm{bd}}}_{\text{diffuse boundary profile}}
\times
\underbrace{\mathcal J[h]}_{\text{unknown outgoing mass flux}}
+\underbrace{r}_{\text{prescribed boundary source}}.
\]
Integrating the characteristic representation against the outgoing mass-flux measure yields a closed boundary relation for $\mathcal J[h]$ involving boundary-to-boundary transmission operators. The strict contraction of these operators yields direct estimates for $\mathcal J[h]$, and hence for $P_\gamma h$, without first estimating the full outgoing boundary trace. This mechanism underlies the wall flux estimates in Lemmas~\ref{lem-wall-flux} and~\ref{lem-wall-flux-time}.

For the steady problem \eqref{eq-wall-transport}, the boundary flux transmission mechanism takes a particularly explicit form in the slab geometry. Since the boundary consists of the two walls $y=\pm1$, the relevant boundary unknowns reduce to the two outgoing mass fluxes $\mathcal J[h](x,1,z)$ and $\mathcal J[h](x,-1,z)$. Along characteristics, the outgoing mass flux associated with one wall is transmitted to the opposite wall through the boundary-to-boundary flux transmission operators $\mathscr B_\lambda(\pm1):L^\infty(\partial\Omega)\to L^\infty(\partial\Omega)$ defined in \eqref{def-wall-B}. More precisely,
\[
\big(\mathscr B_\lambda(\pm1)f\big)(x,\mp1,z)
=
\int_{\pm v_y<0}
e^{-\lambda(v)t_b}
\Phi(x-v_xt_b,\pm1,z-v_zt_b,v)
\sqrt\mu|v_y|
f(x-v_xt_b,\pm1,z-v_zt_b)\,dv,
\]
where $t_b=2/|v_y|$ is the travel time between the two walls, and $e^{-\lambda(v)t_b}$ represents the damping accumulated along the flight. The two wall fluxes satisfy the closed operator system \eqref{eq-wall-system},
\[
\begin{pmatrix}
I & -\mathscr B_\lambda(-1)\\
-\mathscr B_\lambda(1) & I
\end{pmatrix}
\begin{pmatrix}
\mathcal J[h](x,1,z)\\
\mathcal J[h](x,-1,z)
\end{pmatrix}
=
\begin{pmatrix}
\mathcal R_{\lambda,+}(x,z)\\
\mathcal R_{\lambda,-}(x,z)
\end{pmatrix},
\]
where $\mathcal R_{\lambda,\pm}$, defined in \eqref{def-wall-R}, collect the contributions of the interior and incoming boundary sources.

The key idea is the strict contraction condition
\[
\|\mathscr B_\lambda(\pm1)\|_{L^\infty(\partial\Omega)\to L^\infty(\partial\Omega)}<1,
\]
which is verified for the steady linear problem in \eqref{est-Gd-wall-B}. Consequently, the round-trip operators $\mathscr B_\lambda(-1)\mathscr B_\lambda(1)$ and $\mathscr B_\lambda(1)\mathscr B_\lambda(-1)$ are also strictly contractive. The above $2\times2$ system can therefore be solved directly, showing that the two outgoing mass fluxes are uniquely determined by the interior and prescribed incoming boundary sources, as given by \eqref{eq-wall-flux-form} in Lemma~\ref{lem-wall-flux}. These representations are used to derive the boundary flux estimate \eqref{est-Gd-J} and \eqref{est-Gd-J-L2} for the steady linear problem. The time-dependent estimate in Lemma~\ref{lem-wall-flux-time} is based on the same boundary-to-boundary transmission mechanism.

The slab geometry provides an explicit realization of a boundary flux transmission mechanism that extends to damped transport equations with diffuse reflection in general bounded $C^1$ domains. For the steady problem \eqref{eq-general-steady}, Proposition~\ref{prop-general-steady-flux} shows that the outgoing mass flux satisfies the boundary flux equation
\[
\mathcal J[h]
=
\mathscr B_\lambda\mathcal J[h]
+
\mathcal R_\lambda,
\quad
\mathbf x\in\partial\Omega^{'},
\]
where $\mathscr B_\lambda$ is the boundary-to-boundary flux transmission operator defined in \eqref{def-B-lambda}, and $\mathcal R_\lambda$ is defined in \eqref{def-R-lambda}. Unlike the slab geometry, where the closure is represented by an explicit $2\times2$ operator system associated with the two walls, here $\mathcal J[h]$ is a function on the entire physical boundary $\partial\Omega^{'}$, and
\[
\mathscr B_\lambda:
L^\infty(\partial\Omega^{'})
\longrightarrow
L^\infty(\partial\Omega^{'})
\]
describes one boundary-to-boundary transmission of this flux along a characteristic flight.

The key point in the general-domain setting is again the strict contraction of $\mathscr B_\lambda$. The proof of this property in Proposition~\ref{prop-general-steady-flux} is based on the decomposition of $\Gamma_+$ into the non-grazing, bounded-velocity region $\Gamma_+^\varepsilon$ defined in \eqref{def-Gamma-eps} and its complement. The crucial geometric input is Lemma~\ref{lem-general-ft}, which shows that on $\Gamma_+^\varepsilon$ the backward exit time admits a uniform positive lower bound $\tau_{\varepsilon,\Omega^{'}}>0$. Together with $\lambda_\varepsilon:=\inf_{\mathbf x\in\Omega^{'},\,|v|\le\varepsilon^{-1}}\lambda(\mathbf x,v)>0$, this yields uniform attenuation on $\Gamma_+^\varepsilon$. On the complementary region, corresponding to near-grazing or large velocities, the Maxwellian flux measure is small, as quantified by $\delta_\varepsilon$ in \eqref{def-B-delta}. Combining these two mechanisms yields the strict contraction estimate
\[
\|\mathscr B_\lambda f\|_{L^\infty(\partial\Omega^{'})}
\le
\kappa_0\|f\|_{L^\infty(\partial\Omega^{'})},
\quad
0<\kappa_0<1,
\]
where $\kappa_0$ is defined in \eqref{def-B-kappa}. It follows that $I-\mathscr B_\lambda$ is invertible on $L^\infty(\partial\Omega^{'})$. Applying this to the boundary flux equation \eqref{eq-general-flux}, we obtain
\[
\mathcal J[h]
=
(I-\mathscr B_\lambda)^{-1}\mathcal R_\lambda
=
\sum_{j=0}^\infty
\mathscr B_\lambda^j\mathcal R_\lambda
=
\underbrace{\mathcal R_\lambda}_{0\text{ transmissions}}
+
\underbrace{\mathscr B_\lambda\mathcal R_\lambda}_{1\text{ transmission}}
+
\underbrace{\mathscr B_\lambda^2\mathcal R_\lambda}_{2\text{ transmissions}}
+\cdots.
\]
Thus, the interior source and the prescribed boundary source uniquely determine the outgoing mass flux. The time-dependent formulation in Proposition~\ref{prop-general-time-flux} is governed by the same boundary flux equation \eqref{eq-general-flux-time} and strict contraction mechanism.

We finally record two further features of the framework. First, the argument does not rely on the specific form of the Boltzmann collision frequency. It suffices that the damping has a uniform positive lower bound on the bounded-velocity regions considered above, while it may degenerate at large velocities. Second, the strict contraction property extends beyond $L^\infty$. Corollary~\ref{cor-B-Lp} gives
\[
\|\mathscr B_\lambda f\|_{L^p(\partial\Omega^{'})}
\le
\kappa_0\|f\|_{L^p(\partial\Omega^{'})},
\quad
1\le p\le\infty,
\]
and the corresponding time-dependent estimate is given in \eqref{est-B-time-Lp}. This $L^p$ extension may be useful in future applications of the boundary flux transmission framework.

We now turn to the construction of the steady solution in Theorem~\ref{thm-1}.  We first establish a unified linear theory for the three steady corrections by studying the boundary value problem \eqref{eq-Gd-linear}--\eqref{eq-Gd-linear-mass}. To this end, we introduce parameters $0<\epsilon\le1$ and $0\le\sigma\le1$ and consider the approximate boundary value problem \eqref{eq-Gd-eps-sigma}. Combining Guo's $L^2$--$L^\infty$ framework with the boundary flux estimate in Lemma~\ref{lem-wall-flux}, we derive the uniform a priori estimate \eqref{est-Linf-Gd-apriori} in Lemma~\ref{lem-Gd-apriori}, uniformly in $\epsilon$ and $\sigma$. Next, starting from the solvability at $\sigma=0$, we continue the solution in $\sigma$ up to $\sigma=1$, as established in Lemma~\ref{lem-Gd-sigma-cont}. We then pass to the limit $\epsilon\to0$ and recover the zero-mass condition \eqref{eq-Gd-linear-mass}, obtaining the linear solvability result in Lemma~\ref{lem-Gd-linear}. Finally, we apply this linear theory to obtain the first-order shear and thermal corrections in Propositions~\ref{prop-Ga} and~\ref{prop-Gd} and construct the higher-order remainder by Picard iteration, as detailed in Proposition~\ref{prop-GR}.

We now explain several key points in this construction. The approximate problem
\eqref{eq-Gd-eps-sigma} reads
\[
\left\{
\begin{aligned}
&\epsilon\nu(v)G+v\cdot\nabla_{\mathbf x}G+\nu(v)G
=\sigma KG+\mathfrak F,
\\
&G(x,\pm1,z,v)\big|_{\pm v_y<0}
=
\Phi(x,\pm1,z,v)\mathcal J[G](x,\pm1,z)
+r(x,\pm1,z,v).
\end{aligned}
\right.
\]
The penalty $\epsilon\nu(v)G$ is chosen to match the source norm
$\|\nu^{-1}w_l\mathfrak F\|_{L^\infty_{\mathbf x,v}}$
in the uniform a priori estimate. A standard choice is the constant penalty $\epsilon G$, for which
the difference source $(\epsilon'-\epsilon)G_{\epsilon'}$ satisfies $\big\|\nu^{-1}w_l(\epsilon'-\epsilon)G_{\epsilon'}\big\|_{L^\infty_{\mathbf x,v}}=|\epsilon'-\epsilon|\|\nu^{-1}w_lG_{\epsilon'}\|_{L^\infty_{\mathbf x,v}}$.
For soft potentials, the right-hand side is not controlled by the
available uniform bound on $\|w_lG_{\epsilon'}\|_{L^\infty_{\mathbf x,v}}$. One can instead obtain convergence in a weaker velocity weight, as in~\cite{DHWZ19}. We instead choose $\epsilon\nu(v)G$,
which gives the Cauchy estimate
\[
\|w_l(G_\epsilon-G_{\epsilon'})\|_{L^\infty_{\mathbf x,v}}
\le
C_l\big\|\nu^{-1}w_l
(\epsilon'-\epsilon)\nu(v)G_{\epsilon'}\big\|_{L^\infty_{\mathbf x,v}}
=
C_l|\epsilon'-\epsilon|
\|w_lG_{\epsilon'}\|_{L^\infty_{\mathbf x,v}},
\]
as shown in \eqref{est-Linf-Gd-eps-Cauchy}. The same choice preserves the weighted zero-mass condition \eqref{eq-Gd-mass} for all approximation parameters $(\epsilon,\sigma)$. This constraint also controls the mean of the macroscopic density $a$.
Indeed,
\[
\begin{aligned}
0
=
\int_\Omega\int_{\mathbb R^3}
\nu(v)\sqrt\mu G(\mathbf x,v)\,dv d\mathbf x
={}&
\bigg(\int_{\mathbb R^3}\nu(v)\mu\,dv\bigg)
\int_\Omega a\,d\mathbf x
+
\bigg(\underbrace{\int_{\mathbb R^3}\nu(v)\mu v\,dv}_{0}\bigg)
\cdot\int_\Omega\mathbf b\,d\mathbf x
\\
&+
\bigg(\int_{\mathbb R^3}\nu(v)\mu\frac{|v|^2-3}{\sqrt6}\,dv\bigg)
\int_\Omega c\,d\mathbf x
+
\int_\Omega\int_{\mathbb R^3}\nu(v)\sqrt\mu\mathbf P_1G\,dv d\mathbf x.
\end{aligned}
\]
Thus the mean of $a$ is controlled by $\|c\|_{L^2_{\mathbf x}}$
and $\|\nu^{\frac12}\mathbf P_1G\|_{L^2_{\mathbf x,v}}$,
with no contribution from $\mathbf b$. It readily follows that
\[
\|a\|_{L^2_{\mathbf x}}^2
\lesssim
\|\Delta_{\mathbf x}\phi_a\|_{L^2_{\mathbf x}}^2
+
\|c\|_{L^2_{\mathbf x}}^2
+
\|\nu^{\frac12}\mathbf P_1G\|_{L^2_{\mathbf x,v}}^2,
\]
as shown in \eqref{est-L2-Gd-a-mean}. Combining this with \eqref{est-L2-Gd-a-phi} yields the density estimate \eqref{est-L2-Gd-a}. Together with the remaining estimates, this give the a priori estimate in Lemma~\ref{lem-Gd-apriori}, uniformly in $(\epsilon,\sigma)$, before carrying out the continuation.

To recover the usual zero-mass condition \eqref{eq-Gd-linear-mass},
we use linearity to adjust the solution $G_\nu$ of
\eqref{eq-Gd-linear-nu}. If $G_h$ solves the homogeneous problem
\eqref{eq-Gd-hom}, then $G:=G_\nu-c_0G_h$ satisfies the equation
and boundary condition in \eqref{eq-Gd-linear-nu} for any constant $c_0$.
Since $G_h$ has nonzero mass by \eqref{est-Gd-hom-mass},
choosing $c_0$ as in \eqref{def-Gd-c} ensures that $G$ has zero mass.
For uniqueness, we similarly adjust the difference of two solutions
to have weighted zero mass and apply the uniqueness for
\eqref{eq-Gd-linear-nu}. This yields the unique solution to
\eqref{eq-Gd-linear}--\eqref{eq-Gd-linear-mass}
with the estimate \eqref{est-Linf-Gd-linear}.

Having completed the steady analysis, we turn to the unsteady problem
\eqref{eq-prob-3d}. The key to treating initial perturbations with polynomial velocity decay is Caflisch's decomposition
$F(t,\mathbf x,v)-F_{\alpha,\delta}^s(\mathbf x,v)=g_1+\sqrt\mu\,g_2$, as in~\cite{Cf80,DLY25}. This transforms \eqref{eq-g-unsteady}
into \eqref{eq-g1-unsteady}--\eqref{eq-g2-unsteady}.
Here $g_1(0)=F_0-F_{\alpha,\delta}^s$ and $g_2(0)=0$, so the
initial weighted norm does not involve the factor $\mu^{-\frac12}$. Applying Lemma~\ref{lem-wall-flux-time}, we obtain the boundary flux
representation \eqref{eq-g-wall-flux-time}, which expresses the outgoing
mass fluxes in terms of the initial data and interior sources.
This enables a direct contraction mapping proof of local existence
in Lemma~\ref{lem-local-existence}.
Under the a priori assumption \eqref{apriori}, we combine Guo's $L^2$--$L^\infty$
framework with the boundary flux estimates \eqref{est-global-flux-R}
to derive \eqref{est-global-close}. The local theory and a standard
continuation argument then yield global existence, exponential decay
for hard potentials, and the uniform bound for soft potentials. 

For soft potentials, we seek a time weight $W(t)$ in
\eqref{eq-soft-Ag}. Velocity interpolation controls the terms
$W'(t)g_1$ and $W'(t)g_2$ and determines $\beta$ through
\eqref{def-soft-decay-beta}. Solving for $W(t)$ under
\eqref{def-soft-W}, \eqref{est-soft-W-source}, and
\eqref{est-soft-W-L2} gives $W(t)=(1+\varepsilon t)^{2(l-l_0)/|\gamma|}$
for sufficiently small $\varepsilon>0$. The interpolation inequality \eqref{est-soft-interp} then converts the loss $m$ of velocity weight into the algebraic decay rate $(1+t)^{-2m/|\gamma|}$ in \eqref{est-global-decay-soft},
completing the proof of Theorem~\ref{thm-2}.

The rest of the paper is organized as follows. Section~\ref{sec-2} introduces the basic operators and collects the related estimates. Section~\ref{sec-3} establishes the steady and time-dependent wall flux estimates in the Couette channel and develops the boundary flux transmission framework for damped transport equations with diffuse reflection in general bounded $C^1$ domains. Section~\ref{sec-4} establishes the linear theory for the steady corrections, while Section~\ref{sec-5} constructs the first-order shear and thermal corrections and the higher-order remainder. Section~\ref{sec-6} establishes local existence for the unsteady problem, and Section~\ref{sec-7} proves global existence and stability, with exponential decay for hard potentials and algebraic decay for soft potentials.

\section{Preliminaries}\label{sec-2}

We first record several basic identities for the collision operators that will be used throughout the paper. Under the convention adopted in \eqref{eq-coll-op}, the operators introduced above admit the following representations:
\begin{equation}\label{def-coll-ops}
\left\{
\begin{aligned}
&\Gamma(f,g)(v)
=
\int_{\mathbb R^3}\int_{\mathbb S^2}|v-v_*|^\gamma B_0(\cos\theta)\sqrt{\mu}(v_*)
\big[f(v_*')g(v')-f(v_*)g(v)\big]\,d\omega dv_*,
\\
&Lf(v)
=
\nu(v)f(v)-Kf(v),
\\
&\nu(v)
=
\int_{\mathbb R^3}\int_{\mathbb S^2}|v-v_*|^\gamma B_0(\cos\theta)\mu(v_*)\,d\omega dv_*,
\\
&Kf(v)
=
\mu^{-\frac12}(v)\big[Q(\sqrt\mu f,\mu)(v)
+
Q_+(\mu,\sqrt\mu f)(v)\big]
=
\int_{\mathbb R^3}\mathbf{k}(v,v_*)f(v_*)\,dv_*.
\end{aligned}
\right.
\end{equation}
Here the collision frequency $\nu(v)$ satisfies
\[
\nu(v)\sim(1+|v|)^\gamma,\qquad -3<\gamma\le1.
\]

For later use, we also introduce the operators
\begin{equation}\label{def-calL-calK}
\left\{
\begin{aligned}
\mathcal Lf(v)
&:=
-\big[Q(f,\mu)(v)+Q(\mu,f)(v)\big]
=
\nu(v) f(v)-\mathcal Kf(v),
\\
\mathcal Kf(v)
&:=
Q(f,\mu)(v)+Q_+(\mu,f)(v)
=
\sqrt\mu K\big(\mu^{-\frac12}f\big)(v).
\end{aligned}
\right.
\end{equation}
We next introduce the orthogonal decomposition associated with $L$. The null space of $L$ is the five-dimensional space
\begin{equation}\label{def-ker-L}
\ker L
=
\operatorname{span}
\big\{
\chi_0,\chi_1,\chi_2,\chi_3,\chi_4
\big\},
\end{equation}
where
\begin{equation}\label{def-chi}
\chi_0
=
\sqrt\mu,
\quad
\chi_1
=
v_x\sqrt\mu,
\quad
\chi_2
=
v_y\sqrt\mu,
\quad
\chi_3
=
v_z\sqrt\mu,
\quad
\chi_4
=
\frac{|v|^2-3}{\sqrt6}\sqrt\mu.
\end{equation}
The functions $\{\chi_j\}_{j=0}^4$ form an orthonormal basis of $\ker L$ in $L^2(\mathbb R_v^3)$. We define the corresponding macroscopic and microscopic projections by
\begin{equation}\label{def-P01}
\mathbf P_0f
:=
\sum_{j=0}^4\langle f,\chi_j\rangle_v\chi_j,
\quad
\mathbf P_1f
:=
f-\mathbf P_0f.
\end{equation}

We recall the pointwise kernel estimate for the integral operator $K$
from \cite[Lemma~2.2]{DHWY17}. The weighted kernel estimate follows
readily from \cite[Lemmas~2.4 and~2.5]{DHWY17}.

\begin{lemma}\label{lem-K-kernel}
Let $K$ and $\mathbf{k}$ be defined by \eqref{def-coll-ops}.
For $-3<\gamma\le1$, it holds that
\begin{align}
|\mathbf{k}(v,v_*)|
&\le
C|v-v_*|^\gamma
\exp\bigg\{-\frac14(|v|^2+|v_*|^2)\bigg\}
\notag\\
&\quad+
\frac{C}{|v-v_*|^{\frac{3-\gamma}{2}}}
\exp\bigg\{
-\frac18|v-v_*|^2
-\frac18\frac{\big(|v|^2-|v_*|^2\big)^2}{|v-v_*|^2}
\bigg\}.
\label{est-kernel}
\end{align}
For $l\ge0$, define
\begin{equation}\label{def-k-weight}
\mathbf{k}_w(v,v_*)
:=
w_l(v)\mathbf{k}(v,v_*)w_l^{-1}(v_*).
\end{equation}
Then
\begin{equation}\label{est-kernel-weight}
\int_{\mathbb R^3}|\mathbf{k}_w(v,v_*)|\,dv_*
\le
\frac{C_l\nu(v)}{(1+|v|)^2}.
\end{equation}
\end{lemma}

We next collect the velocity-weighted estimates for the nonlinear
collision operators. The $L^2_v$ estimate for $\Gamma$ follows from
\cite[Lemma A.2]{DL25} for $0\le\gamma\le1$ and
\cite[Lemma 2.5]{DHWZ19} for $-3<\gamma<0$.
The $L^\infty_v$ estimate for $\Gamma$ follows from
\eqref{def-coll-ops}, while the
corresponding estimate for $Q$ follows from
\cite[Lemmas 2.9 and 2.11]{Cao22}.

\begin{lemma}\label{lem-Gamma}
For $-3<\gamma\le1$, let $l>1$. Then
\begin{equation}\label{est-L2-Gamma}
\|\nu^{-\frac12}\Gamma(f,g)\|_{L^2_v}
\le
C_l\|w_lf\|_{L^\infty_v}\|\nu^{\frac12}g\|_{L^2_v},
\end{equation}
and
\begin{equation}\label{est-Linf-Gamma}
\|\nu^{-1}w_l\Gamma(f,g)\|_{L^\infty_v}
\le
C_l\|w_lf\|_{L^\infty_v}\|w_lg\|_{L^\infty_v}.
\end{equation}
Moreover, for $l> \dfrac12\max\{3,3+\gamma\}$,
\begin{equation}\label{est-Linf-Q}
\|\nu^{-1}w_lQ(F_1,F_2)\|_{L^\infty_v}
\le
C_l\|w_lF_1\|_{L^\infty_v}\|w_lF_2\|_{L^\infty_v}.
\end{equation}
\end{lemma}

We also recall the coercivity estimate for the linearized collision
operator established in \cite{Guo10} for hard potentials and
in \cite[(2.7)]{DHWZ19} for soft potentials.

\begin{lemma}\label{lem-L-coer}
For $-3<\gamma\le1$, there exists $\lambda_0>0$ such that
\begin{equation}\label{est-L2-L-coer}
\langle Lf,f\rangle_v
=
\langle L\mathbf P_1f,\mathbf P_1f\rangle_v
\ge
\lambda_0\|\nu^{\frac12}\mathbf P_1f\|_{L^2_v}^2.
\end{equation}
\end{lemma}

Finally, we recall the large-velocity smallness estimate for the
integral operator $\mathcal K$. The estimate with the polynomial
velocity weight $w_l$ follows from \cite[Lemma 6.3]{Cao22}, while
the corresponding estimate with the weight $(1+|v|)^{2l}$ is given
in \cite[Lemma 2.2]{DLL26}.

\begin{lemma}\label{lem-calK-high}
Let $\mathcal K$ be defined by \eqref{def-calL-calK}.
For $-3<\gamma\le1$, $l>\dfrac12\max\{3,3+\gamma\}$ and $M>0$,
it holds that
\begin{equation}\label{est-Linf-calK-high}
\sup_{|v|\ge M}
\nu^{-1}(v)w_l(v)|\mathcal Kf(v)|
\le
\bigg(
\frac{C}{l^{\frac{\gamma+3}{4}}}
+
\frac{C_l}{M^2}
\bigg)
\|w_lf\|_{L^\infty_v},
\end{equation}
where $C$ is independent of $l$ and $M$, and $C_l$ is
independent of $M$.
\end{lemma}

\section{Boundary flux transmission framework}\label{sec-3}

In this section, we develop the boundary flux transmission framework. We first establish steady and time-dependent estimates for the outgoing mass
fluxes at the two walls of the Couette channel. We then derive the
boundary flux equation for damped transport with diffuse reflection
in general bounded $C^1$ domains and prove strict contraction of the
associated boundary operator.

\subsection{Steady wall flux estimate}

Let $\lambda=\lambda(v)>0$. Let $\Phi(x,\pm1,z,v)$ be nonnegative wall data satisfying
\begin{equation}\label{eq-wall-Phi-norm}
\int_{\pm v_y<0}
\Phi(x,\pm1,z,v)\sqrt\mu|v_y|\,dv
=
1.
\end{equation}
To capture the boundary flux structure needed in the subsequent steady analysis, we consider the following boundary value problem for $h=h(\mathbf x,v)$:
\begin{equation}\label{eq-wall-transport}
\left\{
\begin{aligned}
&v\cdot\nabla_{\mathbf x}h
+\lambda(v)h
=
\mathcal S(\mathbf x,v),
\quad
(\mathbf x,v)\in\Omega\times\mathbb R^3,
\\
&h(x,\pm1,z,v)\big|_{\pm v_y<0}
=
\Phi(x,\pm1,z,v)\mathcal J[h](x,\pm1,z)
+r(x,\pm1,z,v).
\end{aligned}
\right.
\end{equation}
Here $\mathcal S$ and $r$ are the prescribed
interior and incoming boundary sources, respectively.

Let $s\ge0$ denote the backward travel time from $(\mathbf x,v)$. Then the characteristics associated with \eqref{eq-wall-transport} are given by
\begin{equation}\label{def-wall-char}
\mathbf X(s;\mathbf x,v)=\mathbf x-sv,
\quad
V(s;\mathbf x,v)=v.
\end{equation}

For $(\mathbf x,v)\in(\Omega\times\mathbb R^3)\cup\gamma_+$, define the backward exit time by
\begin{equation}\label{def-wall-tb}
t_b(\mathbf x,v)
:=
\inf\big\{\tau>0:\mathbf X(\tau;\mathbf x,v)\in\partial\Omega\big\}.
\end{equation}

For the boundary-to-boundary characteristics below, $t_b=2/|v_y|$. We define the boundary-to-boundary flux transmission operators $\mathscr B_\lambda(\pm1):L^\infty(\partial\Omega)\to L^\infty(\partial\Omega)$ by
\begin{align}
\big(\mathscr B_\lambda(\pm1)f\big)(x,\mp1,z)
&=
\int_{\pm v_y<0}
e^{-\lambda(v)t_b}
\Phi(x-v_xt_b,\pm1,z-v_zt_b,v)
\sqrt\mu|v_y|
f(x-v_xt_b,\pm1,z-v_zt_b)\,dv,
\label{def-wall-B}
\end{align}
where each $\mathscr B_\lambda(\pm1)$ transfers data from the wall $y=\pm1$ to the opposite wall $y=\mp1$. For the boundary and interior source contributions, set
\begin{align}
\mathcal R_{\lambda,\pm}(x,z)
&=
\int_{\pm v_y>0}\int_0^{t_b}
e^{-\lambda(v)s}
\mathcal S\big(\mathbf X(s;(x,\pm1,z),v),v\big)
\sqrt\mu|v_y|\,dsdv
\notag\\
&\quad+
\int_{\pm v_y>0}
e^{-\lambda(v)t_b}
r(x-v_xt_b,\mp1,z-v_zt_b,v)
\sqrt\mu|v_y|\,dv.
\label{def-wall-R}
\end{align}

The following lemma gives explicit representations and estimates of the outgoing mass fluxes at the two walls for \eqref{eq-wall-transport}.

\begin{lemma}\label{lem-wall-flux}
Let $h$ be a mild solution of \eqref{eq-wall-transport}. Assume that
\begin{equation}\label{est-wall-B}
\max_{\pm}\|\mathscr B_\lambda(\pm1)\|_{L^\infty(\partial\Omega)\to L^\infty(\partial\Omega)}<1.
\end{equation}
Then the two outgoing mass fluxes satisfy
\begin{equation}\label{eq-wall-flux-form}
\left\{
\begin{aligned}
\mathcal J[h](x,1,z)
={}&
\big[I-\mathscr B_\lambda(-1)\mathscr B_\lambda(1)\big]^{-1}
\big[\mathcal R_{\lambda,+}
+\mathscr B_\lambda(-1)\mathcal R_{\lambda,-}\big](x,1,z),
\\
\mathcal J[h](x,-1,z)
={}&
\big[I-\mathscr B_\lambda(1)\mathscr B_\lambda(-1)\big]^{-1}
\big[\mathcal R_{\lambda,-}
+\mathscr B_\lambda(1)\mathcal R_{\lambda,+}\big](x,-1,z).
\end{aligned}
\right.
\end{equation}
Moreover,
\begin{equation}\label{est-wall-flux-Linf}
\|\mathcal J[h]\|_{L^\infty(\partial\Omega)}
\le
C
\Big[
\|\lambda^{-1}\mathcal S\|_{L^\infty_{\mathbf x,v}}
+
|r|_{L^\infty(\gamma_-)}
\Big],
\end{equation}
and
\begin{equation}\label{est-wall-flux-L2}
\|\mathcal J[h]\|_{L^2(\partial\Omega)}^2
\le
C\Big[
\|\lambda^{-\frac12}\mathcal S\|_{L^2_{\mathbf x,v}}^2
+
|r|_{L^2(\gamma_-)}^2
\Big].
\end{equation}
Here $C>0$ depends on $\lambda$ and $\Phi$.
\end{lemma}

\begin{proof}
Along the backward characteristics \eqref{def-wall-char}, we have
\begin{equation}\label{eq-duh-h}
\frac{d}{ds}
h\big(\mathbf X(s;\mathbf x,v),v\big)
=
\lambda(v)h\big(\mathbf X(s;\mathbf x,v),v\big)
-
\mathcal S\big(\mathbf X(s;\mathbf x,v),v\big).
\end{equation}

For $v_y>0$, consider the outgoing trace at the upper wall $y=1$. The backward characteristic issued from $((x,1,z),v)$ reaches the lower wall $y=-1$ after the backward travel time $t_b=2/v_y$ at
\[
\mathbf X(t_b)
=
(x-v_xt_b,-1,z-v_zt_b).
\]
Applying Duhamel's formula for \eqref{eq-duh-h} and using the incoming boundary condition in \eqref{eq-wall-transport}, we obtain
\begin{align}
h(x,1,z,v)
&=
e^{-\lambda(v)t_b}\Big[
\Phi\big(X(t_b),-1,Z(t_b),v\big)
\mathcal J[h]\big(X(t_b),-1,Z(t_b)\big)
+r\big(X(t_b),-1,Z(t_b),v\big)\Big]
\notag\\
&\quad+
\int_0^{t_b}e^{-\lambda(v)s}
\mathcal S\big(\mathbf X(s;(x,1,z),v),v\big)\,ds.
\label{eq-wall-Duhamel--1}
\end{align}
Integrating \eqref{eq-wall-Duhamel--1} over $v_y>0$ against $\sqrt\mu v_y\,dv$, we obtain
\begin{equation}\label{eq-wall-flux-plus}
\mathcal J[h](x,1,z)
-
\mathscr B_\lambda(-1)\mathcal J[h](x,-1,z)
=
\mathcal R_{\lambda,+}(x,z).
\end{equation}

For $v_y<0$, consider the outgoing trace at the lower wall $y=-1$. Similarly, the backward characteristic from $((x,-1,z),v)$ reaches the upper wall at
\[
\mathbf X(t_b)
=
(x-v_xt_b,1,z-v_zt_b).
\]
Proceeding as above, we obtain
\begin{equation}\label{eq-wall-flux-minus}
\mathcal J[h](x,-1,z)
-
\mathscr B_\lambda(1)\mathcal J[h](x,1,z)
=
\mathcal R_{\lambda,-}(x,z).
\end{equation}
Equivalently, \eqref{eq-wall-flux-plus}--\eqref{eq-wall-flux-minus} can be written as
\begin{equation}\label{eq-wall-system}
\begin{pmatrix}
I & -\mathscr B_\lambda(-1)\\
-\mathscr B_\lambda(1) & I
\end{pmatrix}
\begin{pmatrix}
\mathcal J[h](x,1,z)\\
\mathcal J[h](x,-1,z)
\end{pmatrix}
=
\begin{pmatrix}
\mathcal R_{\lambda,+}(x,z)\\
\mathcal R_{\lambda,-}(x,z)
\end{pmatrix}.
\end{equation}
By the contraction condition \eqref{est-wall-B},
\[
\begin{aligned}
&\max\Big\{
\|\mathscr B_\lambda(-1)\mathscr B_\lambda(1)\|_{L^\infty(\partial\Omega)\to L^\infty(\partial\Omega)},
\|\mathscr B_\lambda(1)\mathscr B_\lambda(-1)\|_{L^\infty(\partial\Omega)\to L^\infty(\partial\Omega)}
\Big\}
\\
&\le
\Big[
\max_{\pm}\|\mathscr B_\lambda(\pm1)\|_{L^\infty(\partial\Omega)\to L^\infty(\partial\Omega)}
\Big]^2<1.
\end{aligned}
\]
Thus both $I-\mathscr B_\lambda(-1)\mathscr B_\lambda(1)$ and $I-\mathscr B_\lambda(1)\mathscr B_\lambda(-1)$ are invertible on $L^\infty(\partial\Omega)$. Eliminating either wall flux from \eqref{eq-wall-flux-plus}--\eqref{eq-wall-flux-minus} yields \eqref{eq-wall-flux-form}.

Moreover, by \eqref{def-wall-R} and
\[
\int_0^{t_b}e^{-\lambda(v)s}\lambda(v)\,ds
\le1,
\]
we have
\[
\|\mathcal R_{\lambda,\pm}\|_{L^\infty_{x,z}}
\lesssim
\|\lambda^{-1}\mathcal S\|_{L^\infty_{\mathbf x,v}}
+
|r|_{L^\infty(\gamma_-)}.
\]
Then \eqref{est-wall-B}--\eqref{eq-wall-flux-form} yield \eqref{est-wall-flux-Linf}. 

For the $L^2$ estimate, by the Cauchy--Schwarz inequality,
\[
\begin{aligned}
\big|\mathscr B_\lambda(\pm1)f\big|^2
&\le
\big[\mathscr B_\lambda(\pm1)1\big]
\big[\mathscr B_\lambda(\pm1)|f|^2\big]
\\
&\le
\|\mathscr B_\lambda(\pm1)\|_{L^\infty(\partial\Omega)\to L^\infty(\partial\Omega)}
\mathscr B_\lambda(\pm1)|f|^2.
\end{aligned}
\]
Integrating over $\partial\Omega$ and using the change of variables
$(x,z)\mapsto(x-v_xt_b,z-v_zt_b)$, and the normalization \eqref{eq-wall-Phi-norm}, we obtain
\begin{equation}\label{est-wall-B-L2}
\|\mathscr B_\lambda(\pm1)f\|_{L^2(\partial\Omega)}^2
\le
\|\mathscr B_\lambda(\pm1)\|_{L^\infty(\partial\Omega)\to L^\infty(\partial\Omega)}
\|f\|_{L^2(\partial\Omega)}^2.
\end{equation}
Furthermore,
\begin{align}
\int_{\mathbb T^2}|\mathcal R_{\lambda,\pm}(x,z)|^2\,dxdz
&\lesssim
\int_{\mathbb T^2}\int_{\pm v_y>0}\int_0^{t_b}
e^{-\lambda(v)s}\lambda^{-1}(v)
\big|\mathcal S\big(\mathbf X(s;(x,\pm1,z),v),v\big)\big|^2
|v_y|\,dsdv dxdz
\notag\\
&\quad+
\int_{\mathbb T^2}\int_{\pm v_y>0}
\big|r(x-v_xt_b,\mp1,z-v_zt_b,v)\big|^2
|v_y|\,dv dxdz
\notag\\
&\lesssim
\|\lambda^{-\frac12}\mathcal S\|_{L^2_{\mathbf x,v}}^2
+
|r|_{L^2(\gamma_-)}^2,
\label{est-wall-R-L2}
\end{align}
where we have used
\[
\mathbf X(s;(x,\pm1,z),v)=(x-v_xs,\pm1-v_ys,z-v_zs),
\qquad
d\mathbf X=|v_y|\,dsdxdz.
\]

Combining \eqref{eq-wall-flux-form}, \eqref{est-wall-B-L2}, and
\eqref{est-wall-R-L2} yields \eqref{est-wall-flux-L2}.
This completes the proof of Lemma~\ref{lem-wall-flux}.
\end{proof}

\subsection{Time-dependent wall flux estimate}

We next consider the corresponding time-dependent problem for $h=h(t,\mathbf x,v)$:
\begin{equation}\label{eq-wall-transport-time}
\left\{
\begin{aligned}
&\partial_t h
+v\cdot\nabla_{\mathbf x}h
+\lambda(v)h
=
\mathcal S(t,\mathbf x,v),
\quad
(t,\mathbf x,v)\in\mathbb R_+\times\Omega\times\mathbb R^3,
\\
&h(t,x,\pm1,z,v)\big|_{\pm v_y<0}
=
\Phi(x,\pm1,z,v)\mathcal J[h](t,x,\pm1,z)
+r(t,x,\pm1,z,v),
\\
&h(0,\mathbf x,v)
=
h_0(\mathbf x,v).
\end{aligned}
\right.
\end{equation}

For $0\le s\le t$, let $s$ denote the physical time along the backward characteristic terminating at $(\mathbf x,v)$ at time $t$. Then the characteristics associated with \eqref{eq-wall-transport-time} are given by
\begin{equation}\label{def-wall-char-time}
\mathbf X(s;t,\mathbf x,v)=\mathbf x-(t-s)v,
\quad
V(s;t,\mathbf x,v)=v.
\end{equation}

For the boundary-to-boundary characteristics below, $t_b=2/|v_y|$. We define the time-dependent boundary-to-boundary flux transmission operators by
\begin{align}
\big(\mathscr B_\lambda^t(\pm1)f\big)(t,x,\mp1,z)
&=
\int_{\pm v_y<0}
\mathbf 1_{\{t_b<t\}}e^{-\lambda(v)t_b}
\Phi(x-v_xt_b,\pm1,z-v_zt_b,v)
\sqrt\mu|v_y|
\notag\\
&\quad\times
f(t-t_b,x-v_xt_b,\pm1,z-v_zt_b)\,dv.
\label{def-wall-B-time}
\end{align}
For the boundary, initial, and interior source contributions, set
\begin{align}
\mathcal R_{\lambda,\pm}(t,x,z)
&=
\int_{\pm v_y>0}\int_{\max\{0,t-t_b\}}^t
e^{-\lambda(v)(t-\tau)}
\mathcal S\big(\tau,\mathbf X(\tau;t,(x,\pm1,z),v),v\big)
\sqrt\mu|v_y|\,d\tau dv
\notag\\
&\quad+
\int_{\pm v_y>0}
\mathbf 1_{\{t_b\ge t\}}e^{-\lambda(v)t}
h_0\big(\mathbf X(0;t,(x,\pm1,z),v),v\big)
\sqrt\mu|v_y|\,dv
\notag\\
&\quad+
\int_{\pm v_y>0}
\mathbf 1_{\{t_b<t\}}e^{-\lambda(v)t_b}
r(t-t_b,x-v_xt_b,\mp1,z-v_zt_b,v)
\sqrt\mu|v_y|\,dv.
\label{def-wall-R-time}
\end{align}

The following lemma gives explicit representations and estimates of the time-dependent outgoing mass fluxes at the two walls for \eqref{eq-wall-transport-time}.

\begin{lemma}\label{lem-wall-flux-time}
Let $T>0$ and let $h$ be a mild solution of \eqref{eq-wall-transport-time}. Assume that
\begin{equation}\label{est-wall-B-time}
\sup_{t>0}\max_{\pm}
\|\mathscr B_\lambda^t(\pm1)\|_{L^\infty(0,t;L^\infty(\partial\Omega))\to L^\infty(0,t;L^\infty(\partial\Omega))}
<1.
\end{equation}
Then the two outgoing mass fluxes satisfy
\begin{equation}\label{eq-wall-flux-form-time}
\left\{
\begin{aligned}
\mathcal J[h](t,x,1,z)
={}&
\big[I-\mathscr B_\lambda^t(-1)\mathscr B_\lambda^t(1)\big]^{-1}
\big[\mathcal R_{\lambda,+}
+\mathscr B_\lambda^t(-1)\mathcal R_{\lambda,-}\big](t,x,1,z),
\\
\mathcal J[h](t,x,-1,z)
={}&
\big[I-\mathscr B_\lambda^t(1)\mathscr B_\lambda^t(-1)\big]^{-1}
\big[\mathcal R_{\lambda,-}
+\mathscr B_\lambda^t(1)\mathcal R_{\lambda,+}\big](t,x,-1,z).
\end{aligned}
\right.
\end{equation}
Moreover,
\begin{align}
&\|\mathcal J[h]\|_{L^\infty(0,T;L^\infty(\partial\Omega))}
\le
C\Big[
\|h_0\|_{L^\infty_{\mathbf x,v}}
+
\|\lambda^{-1}\mathcal S\|_{L^\infty(0,T;L^\infty_{\mathbf x,v})}
+
|r|_{L^\infty((0,T)\times\gamma_-)}
\Big],
\label{est-wall-flux-time-Linf}
\end{align}
and
\begin{align}
\int_0^T\|\mathcal J[h](t)\|_{L^2(\partial\Omega)}^2\,dt
\le
C\bigg[
&\|h_0\|_{L^2_{\mathbf x,v}}^2
+
\int_0^T\|\lambda^{-\frac12}\mathcal S(t)\|_{L^2_{\mathbf x,v}}^2\,dt
+
\int_0^T|r(t)|_{L^2(\gamma_-)}^2\,dt
\bigg].
\label{est-wall-flux-time-L2}
\end{align}
Here $C>0$ depends on $\lambda$ and $\Phi$ but is independent of $T$.
\end{lemma}

\begin{proof}
The argument is analogous to that of Lemma~\ref{lem-wall-flux}, with the additional distinction between whether the backward characteristic reaches the opposite wall or the initial surface first.

For $v_y>0$, we have
\begin{align}
h(t,x,1,z,v)
&=
\mathbf 1_{\{t_b<t\}}e^{-\lambda(v)t_b}
\Big[
\Phi\big(X(t-t_b),-1,Z(t-t_b),v\big)
\notag\\
&\qquad\times
\mathcal J[h]\big(t-t_b,X(t-t_b),-1,Z(t-t_b)\big)
+r\big(t-t_b,X(t-t_b),-1,Z(t-t_b),v\big)
\Big]
\notag\\
&\quad+
\mathbf 1_{\{t_b\ge t\}}e^{-\lambda(v)t}
h_0\big(\mathbf X(0;t,(x,1,z),v),v\big)
\notag\\
&\quad+
\int_{\max\{0,t-t_b\}}^t
e^{-\lambda(v)(t-\tau)}
\mathcal S\big(\tau,\mathbf X(\tau;t,(x,1,z),v),v\big)\,d\tau.
\label{eq-wall-Duhamel-time--1}
\end{align}
Integrating \eqref{eq-wall-Duhamel-time--1} over $v_y>0$ against $\sqrt\mu v_y\,dv$, we obtain
\begin{equation}\label{eq-wall-flux-plus-time}
\mathcal J[h](t,x,1,z)
-\mathscr B_\lambda^t(-1)\mathcal J[h](t,x,-1,z)
=
\mathcal R_{\lambda,+}(t,x,z).
\end{equation}
Similarly, for $v_y<0$, we obtain
\begin{equation}\label{eq-wall-flux-minus-time}
\mathcal J[h](t,x,-1,z)
-\mathscr B_\lambda^t(1)\mathcal J[h](t,x,1,z)
=
\mathcal R_{\lambda,-}(t,x,z).
\end{equation}
Equivalently,
\begin{equation}\label{eq-wall-system-time}
\begin{pmatrix}
I & -\mathscr B_\lambda^t(-1)\\
-\mathscr B_\lambda^t(1) & I
\end{pmatrix}
\begin{pmatrix}
\mathcal J[h](t,x,1,z)\\
\mathcal J[h](t,x,-1,z)
\end{pmatrix}
=
\begin{pmatrix}
\mathcal R_{\lambda,+}(t,x,z)\\
\mathcal R_{\lambda,-}(t,x,z)
\end{pmatrix}.
\end{equation}
By the contraction condition \eqref{est-wall-B-time}, and \eqref{eq-wall-flux-plus-time}--\eqref{eq-wall-flux-minus-time} yields \eqref{eq-wall-flux-form-time}.

Moreover, by \eqref{def-wall-R-time}, we have
\[
\|\mathcal R_{\lambda,\pm}\|_{L^\infty(0,T;L^\infty_{x,z})}
\lesssim
\|h_0\|_{L^\infty_{\mathbf x,v}}
+
\|\lambda^{-1}\mathcal S\|_{L^\infty(0,T;L^\infty_{\mathbf x,v})}
+
|r|_{L^\infty((0,T)\times\gamma_-)}.
\]
Then \eqref{eq-wall-flux-form-time} and \eqref{est-wall-B-time} yield \eqref{est-wall-flux-time-Linf}.

Arguing as in the proof of \eqref{est-wall-B-L2}, we use the translation $(t,x,z)\mapsto(t-t_b,x-v_xt_b,z-v_zt_b)$.
Thus, by \eqref{eq-wall-Phi-norm} and \eqref{est-wall-B-time}, we obtain
\begin{align}
\|\mathscr B_\lambda^t(\pm1)f\|_{L^2((0,T)\times\partial\Omega)}^2
&\le
\|\mathscr B_\lambda^t(\pm1)\|_{L^\infty(0,T;L^\infty(\partial\Omega))\to L^\infty(0,T;L^\infty(\partial\Omega))}
\|f\|_{L^2((0,T)\times\partial\Omega)}^2
\notag\\
&<\|f\|_{L^2((0,T)\times\partial\Omega)}^2,
\label{est-wall-B-time-L2}
\end{align}
Moreover, by \eqref{def-wall-R-time},
\begin{align}
&\int_0^T\int_{\mathbb T^2}|\mathcal R_{\lambda,\pm}(t,x,z)|^2\,dxdzdt
\notag\\
&\lesssim
\int_0^T\int_{\mathbb T^2}\int_{\pm v_y>0}\int_{\max\{0,t-t_b\}}^t
e^{-\lambda(v)(t-\tau)}\lambda^{-1}(v)
\big|\mathcal S\big(\tau,\mathbf X(\tau;t,(x,\pm1,z),v),v\big)\big|^2
|v_y|\,d\tau dv dxdzdt
\notag\\
&\quad+
\int_0^T\int_{\mathbb T^2}\int_{\pm v_y>0}
\mathbf 1_{\{t_b\ge t\}}
\big|h_0\big(\mathbf X(0;t,(x,\pm1,z),v),v\big)\big|^2
|v_y|\,dv dxdzdt
\notag\\
&\quad+
\int_0^T\int_{\mathbb T^2}\int_{\pm v_y>0}
\mathbf 1_{\{t_b<t\}}
\big|r(t-t_b,x-v_xt_b,\mp1,z-v_zt_b,v)\big|^2
|v_y|\,dv dxdzdt
\notag\\
&\lesssim
\int_0^T\|\lambda^{-\frac12}\mathcal S(t)\|_{L^2_{\mathbf x,v}}^2\,dt
+
\|h_0\|_{L^2_{\mathbf x,v}}^2
+
\int_0^T|r(t)|_{L^2(\gamma_-)}^2\,dt,
\label{est-wall-R-time-L2}
\end{align}
where we have use
\[
\mathbf X(\tau;t,(x,\pm1,z),v)
=(x-v_x(t-\tau),\pm1-v_y(t-\tau),z-v_z(t-\tau)),\qquad
d\mathbf X=|v_y|\,dtdxdz
\]
for fixed $\tau$ and $v$.

Combining \eqref{eq-wall-flux-form-time}, \eqref{est-wall-B-time-L2}, and \eqref{est-wall-R-time-L2} yields \eqref{est-wall-flux-time-L2}. This completes the proof of Lemma~\ref{lem-wall-flux-time}.
\end{proof}

\subsection{Boundary flux transmission in general bounded domains}\label{subsec3.3}

The preceding wall flux estimates exhibit, in the slab geometry, the boundary-to-boundary transmission of the outgoing mass flux generated by diffuse reflection. We now extend this mechanism to damped transport equations with diffuse reflection in general bounded $C^1$ domains. For the steady problem, we derive the boundary flux equation and prove the strict contraction of the associated boundary-to-boundary flux transmission operator. The time-dependent problem is treated subsequently by the same mechanism.

To distinguish the notation in the general-domain setting from that used for the Couette channel, let $\Omega^{'}\subset\mathbb R^3$ be a bounded $C^1$ domain, and let $n(\mathbf x)$ denote the outward unit normal at $\mathbf x\in\partial\Omega^{'}$. We set $\Gamma:=\partial\Omega^{'}\times\mathbb R_v^3$ and decompose it into the outgoing, incoming, and grazing sets
\[
\left\{
\begin{aligned}
\Gamma_+
&:=\{(\mathbf x,v)\in\partial\Omega^{'}\times\mathbb R^3:n(\mathbf x)\cdot v>0\},\\
\Gamma_-
&:=\{(\mathbf x,v)\in\partial\Omega^{'}\times\mathbb R^3:n(\mathbf x)\cdot v<0\},\\
\Gamma_0
&:=\{(\mathbf x,v)\in\partial\Omega^{'}\times\mathbb R^3:n(\mathbf x)\cdot v=0\}.
\end{aligned}
\right.
\]
The kinetic boundary measure on $\Gamma$ is $d\Gamma:=|n(\mathbf x)\cdot v|\,dv dS_{\mathbf x}$, where $dS_{\mathbf x}$ denotes the surface measure on $\partial\Omega^{'}$. Throughout this subsection, all other norms are understood in the standard sense, and their definitions are omitted for brevity.

We begin with the steady problem for $h=h(\mathbf x,v)$:
\begin{equation}\label{eq-general-steady}
\left\{
\begin{aligned}
&v\cdot\nabla_{\mathbf x}h+\lambda(\mathbf x,v)h
=\mathcal S(\mathbf x,v),
\quad
(\mathbf x,v)\in\Omega^{'}\times\mathbb R^3,\\
&h(\mathbf x,v)\big|_{\Gamma_-}
=\sqrt{2\pi}\sqrt\mu\mathcal J[h](\mathbf x)+r(\mathbf x,v).
\end{aligned}
\right.
\end{equation}
Here $\mathcal J[h]$ denotes the outgoing mass flux defined in \eqref{def-outgoing-flux}, with $\Omega$ replaced by $\Omega^{'}$. The associated backward characteristics are
\begin{equation}\label{def-general-char}
\mathbf X(s;\mathbf x,v)=\mathbf x-sv,
\quad
V(s;\mathbf x,v)=v,
\quad
s\ge0.
\end{equation}
For $(\mathbf x,v)\in(\Omega^{'}\times\mathbb R^3)\cup\Gamma_+$, let $t_b(\mathbf x,v)$ denote the backward exit time associated with $\Omega^{'}$, defined as in \eqref{def-wall-tb} with $\Omega$ replaced by $\Omega^{'}$.

The following lemma provides the key geometric ingredient for establishing the strict contraction of the boundary-to-boundary flux transmission operators: a uniform positive lower bound for the backward exit time away from the near-grazing and large-velocity regions.

\begin{lemma}\label{lem-general-ft}
Let $\Omega^{'}\subset\mathbb R^3$ be a bounded $C^1$ domain. For any
$0<\varepsilon<1$, define
\begin{equation}\label{def-Gamma-eps}
\Gamma_+^\varepsilon
:=
\left\{
(\mathbf x,v)\in\Gamma_+:
n(\mathbf x)\cdot v\ge\varepsilon,\quad
|v|\le\frac1\varepsilon
\right\}.
\end{equation}
Then there exists a constant $\tau_{\varepsilon,\Omega^{'}}>0$, depending only on $\varepsilon$ and the geometry of $\Omega^{'}$, such that
\begin{equation}\label{est-general-flight-time}
t_b(\mathbf x,v)\ge\tau_{\varepsilon,\Omega^{'}},
\quad
(\mathbf x,v)\in\Gamma_+^\varepsilon.
\end{equation}
\end{lemma}

\begin{proof}
By \cite[Lemma~4.3(1)]{BG16}, the backward exit time $t_b(\mathbf x,v)$ is lower semicontinuous with respect to $(\mathbf x,v)$.

Suppose, to the contrary, that
\[
\inf_{(\mathbf x,v)\in\Gamma_+^\varepsilon}t_b(\mathbf x,v)=0.
\]
Then there exists a sequence $(\mathbf x_n,v_n)\in\Gamma_+^\varepsilon$ such that
\[
t_b(\mathbf x_n,v_n)\to0
\quad\text{as }n\to\infty.
\]
Since $\Gamma_+^\varepsilon$ is compact, up to a subsequence,
\[
(\mathbf x_n,v_n)\to(\mathbf x_*,v_*)\in\Gamma_+^\varepsilon.
\]
By the lower semicontinuity of $t_b$,
\[
t_b(\mathbf x_*,v_*)
\le
\liminf_{n\to\infty}t_b(\mathbf x_n,v_n)
=0.
\]
Thus $t_b(\mathbf x_*,v_*)=0$. However, since $n(\mathbf x_*)\cdot v_*\ge\varepsilon>0$, the backward characteristic $\mathbf x_*-sv_*$ lies in $\Omega^{'}$ for all sufficiently small $s>0$. Therefore, $t_b(\mathbf x_*,v_*)>0$, a contradiction. This completes the proof of Lemma~\ref{lem-general-ft}.
\end{proof}

We define the boundary-to-boundary flux transmission operator by
\begin{equation}\label{def-B-lambda}
\mathscr B_\lambda f(\mathbf x)
:=
\sqrt{2\pi}
\int_{n(\mathbf x)\cdot v>0}
e^{-\int_0^{t_b(\mathbf x,v)}\lambda(\mathbf X(\tau;\mathbf x,v),v)\,d\tau}
\mu\big[n(\mathbf x)\cdot v\big]
f\big(\mathbf X(t_b(\mathbf x,v);\mathbf x,v)\big)\,dv.
\end{equation}
We also define
\begin{align}
\mathcal R_\lambda(\mathbf x)
={}&
\int_{n(\mathbf x)\cdot v>0}\int_0^{t_b(\mathbf x,v)}
e^{-\int_0^s\lambda(\mathbf X(\tau;\mathbf x,v),v)\,d\tau}
\mathcal S(\mathbf X(s;\mathbf x,v),v)
\sqrt\mu\big[n(\mathbf x)\cdot v\big]\,dsdv \notag\\
&+
\int_{n(\mathbf x)\cdot v>0}
e^{-\int_0^{t_b(\mathbf x,v)}\lambda(\mathbf X(\tau;\mathbf x,v),v)\,d\tau}
r\big(\mathbf X(t_b(\mathbf x,v);\mathbf x,v),v\big)
\sqrt\mu\big[n(\mathbf x)\cdot v\big]\,dv.
\label{def-R-lambda}
\end{align}

The following proposition gives a closed boundary flux equation for the outgoing mass flux in a general bounded $C^1$ domain.

\begin{proposition}\label{prop-general-steady-flux}
Let $\lambda=\lambda(\mathbf x,v)>0$ satisfy, for any $0<\varepsilon<1$ as in Lemma~\ref{lem-general-ft},
\[
\lambda_\varepsilon:=\inf_{\mathbf x\in\Omega^{'},\,|v|\le\varepsilon^{-1}}\lambda(\mathbf x,v)>0.
\]
Let $h$ be a mild solution of \eqref{eq-general-steady}. Then the outgoing mass flux $\mathcal J[h]$ satisfies
\begin{equation}\label{eq-general-flux}
\mathcal J[h]
=
\mathscr B_\lambda\mathcal J[h]+\mathcal R_\lambda,
\quad
\mathbf x\in\partial\Omega^{'}.
\end{equation}
Moreover, there exists a constant $\kappa_0\in(0,1)$ such that
\begin{equation}\label{est-B-contraction}
\|\mathscr B_\lambda\|_{L^\infty(\partial\Omega^{'})\to L^\infty(\partial\Omega^{'})}
\le
\kappa_0.
\end{equation}
Consequently, $I-\mathscr B_\lambda$ is invertible on $L^\infty(\partial\Omega^{'})$, and
\begin{equation}\label{eq-flux-iteration}
\mathcal J[h]
=
\sum_{j=0}^\infty\mathscr B_\lambda^j\mathcal R_\lambda.
\end{equation}
In particular,
\begin{equation}\label{est-general-steady-flux}
\|\mathcal J[h]\|_{L^\infty(\partial\Omega^{'})}
\le
C_{\Omega^{'},\kappa_0}
\Big[
\|\mathcal S\|_{L^\infty(\Omega^{'}\times\mathbb R^3)}
+
|r|_{L^\infty(\Gamma_-)}
\Big].
\end{equation}
\end{proposition}

\begin{proof}
For $(\mathbf x,v)\in\Gamma_+$, write $t_b=t_b(\mathbf x,v)$. The Duhamel formula along the backward characteristic, together with the boundary condition in \eqref{eq-general-steady}, gives
\[
\begin{aligned}
h(\mathbf x,v)
={}&
\sqrt{2\pi}\sqrt\mu
e^{-\int_0^{t_b}\lambda(\mathbf X(\tau;\mathbf x,v),v)\,d\tau}
\mathcal J[h]\big(\mathbf X(t_b;\mathbf x,v)\big)
+
e^{-\int_0^{t_b}\lambda(\mathbf X(\tau;\mathbf x,v),v)\,d\tau}
r\big(\mathbf X(t_b;\mathbf x,v),v\big)
\\
&+
\int_0^{t_b}
e^{-\int_0^s\lambda(\mathbf X(\tau;\mathbf x,v),v)\,d\tau}
\mathcal S(\mathbf X(s;\mathbf x,v),v)\,ds.
\end{aligned}
\]
Multiplying by $\sqrt\mu\big[n(\mathbf x)\cdot v\big]$ and integrating over $n(\mathbf x)\cdot v>0$, we obtain \eqref{eq-general-flux} from \eqref{def-B-lambda} and \eqref{def-R-lambda}.

We next prove the strict contraction property \eqref{est-B-contraction}. For $0<\varepsilon<1$, set
\begin{equation}\label{def-B-delta}
\delta_\varepsilon
:=
\sup_{\mathbf x\in\partial\Omega^{'}}
\sqrt{2\pi}
\int_{\{v:(\mathbf x,v)\in\Gamma_+\setminus\Gamma_+^\varepsilon\}}
\mu\big[n(\mathbf x)\cdot v\big]\,dv.
\end{equation}
A direct estimate yields
\[
\delta_\varepsilon
\le
C\varepsilon+Ce^{-\frac{c}{\varepsilon^2}},
\]
and hence $\delta_\varepsilon\to0$ as $\varepsilon\to0$. Fix $0<\varepsilon_0<1$ sufficiently small such that $\delta_{\varepsilon_0}\le1/2$.
By \eqref{est-general-flight-time} in Lemma~\ref{lem-general-ft},
\[
t_b(\mathbf x,v)\ge\tau_{\varepsilon_0,\Omega^{'}}>0,
\quad
(\mathbf x,v)\in\Gamma_+^{\varepsilon_0}.
\]
Since $|v|\le\varepsilon_0^{-1}$ on $\Gamma_+^{\varepsilon_0}$, the definition of $\lambda_{\varepsilon_0}$ and \eqref{est-general-flight-time} imply that the attenuation factor on $\Gamma_+^{\varepsilon_0}$ is bounded by $e^{-\lambda_{\varepsilon_0}\tau_{\varepsilon_0,\Omega^{'}}}$. Consequently, for any $f\in L^\infty(\partial\Omega^{'})$,
\begin{align}
\big|(\mathscr B_\lambda f)(\mathbf x)\big|
={}&
\sqrt{2\pi}\bigg|
\int_{n(\mathbf x)\cdot v>0}
e^{-\int_0^{t_b(\mathbf x,v)}\lambda(\mathbf X(\tau;\mathbf x,v),v)\,d\tau}
\mu[n(\mathbf x)\cdot v]
f\big(\mathbf X(t_b(\mathbf x,v);\mathbf x,v)\big)\,dv
\bigg|
\notag\\
\le{}&
\sqrt{2\pi}\|f\|_{L^\infty(\partial\Omega^{'})}
\bigg[
e^{-\lambda_{\varepsilon_0}\tau_{\varepsilon_0,\Omega^{'}}}
\int_{\{v:(\mathbf x,v)\in\Gamma_+^{\varepsilon_0}\}}
\mu[n(\mathbf x)\cdot v]\,dv
+
\int_{\{v:(\mathbf x,v)\in\Gamma_+\setminus\Gamma_+^{\varepsilon_0}\}}
\mu[n(\mathbf x)\cdot v]\,dv
\bigg]
\notag\\
\le{}&
\Big[
e^{-\lambda_{\varepsilon_0}\tau_{\varepsilon_0,\Omega^{'}}}
(1-\delta_{\varepsilon_0})
+\delta_{\varepsilon_0}
\Big]
\|f\|_{L^\infty(\partial\Omega^{'})}.
\end{align}
Set
\begin{equation}\label{def-B-kappa}
\kappa_0
:=
e^{-\lambda_{\varepsilon_0}\tau_{\varepsilon_0,\Omega^{'}}}
(1-\delta_{\varepsilon_0})
+
\delta_{\varepsilon_0}.
\end{equation}
Since $\lambda_{\varepsilon_0}>0$, $\tau_{\varepsilon_0,\Omega^{'}}>0$, and $\delta_{\varepsilon_0}<1$, we have $0<\kappa_0<1$, which proves \eqref{est-B-contraction}.

Moreover, since $\Omega^{'}$ is bounded,
\[
t_b(\mathbf x,v)\big[n(\mathbf x)\cdot v\big]
\le
t_b(\mathbf x,v)|v|
=
\big|\mathbf x-\mathbf X(t_b(\mathbf x,v);\mathbf x,v)\big|
\le
C_{\Omega^{'}}.
\]
It follows from \eqref{def-R-lambda} that
\begin{align}
\|\mathcal R_\lambda\|_{L^\infty(\partial\Omega^{'})}
\le{}&
\|\mathcal S\|_{L^\infty(\Omega^{'}\times\mathbb R^3)}
\sup_{\mathbf x\in\partial\Omega^{'}}
\int_{n(\mathbf x)\cdot v>0}\int_0^{t_b(\mathbf x,v)}
\sqrt\mu\big[n(\mathbf x)\cdot v\big]\,dsdv
\notag\\
&+
|r|_{L^\infty(\Gamma_-)}
\sup_{\mathbf x\in\partial\Omega^{'}}
\int_{n(\mathbf x)\cdot v>0}
\sqrt\mu\big[n(\mathbf x)\cdot v\big]\,dv
\notag\\
\le{}&
C_{\Omega^{'}}
\Big[
\|\mathcal S\|_{L^\infty(\Omega^{'}\times\mathbb R^3)}
+
|r|_{L^\infty(\Gamma_-)}
\Big].
\end{align}
Moreover, \eqref{est-B-contraction} implies that $I-\mathscr B_\lambda$ is invertible on $L^\infty(\partial\Omega^{'})$, with
\[
(I-\mathscr B_\lambda)^{-1}
=
\sum_{j=0}^\infty\mathscr B_\lambda^j.
\]
Hence \eqref{eq-general-flux} yields \eqref{eq-flux-iteration} and
\[
\|\mathcal J[h]\|_{L^\infty(\partial\Omega^{'})}
\le
\frac1{1-\kappa_0}
\|\mathcal R_\lambda\|_{L^\infty(\partial\Omega^{'})},
\]
which gives \eqref{est-general-steady-flux}. This completes the proof of Proposition~\ref{prop-general-steady-flux}.
\end{proof}

We next consider the corresponding time-dependent problem for $h=h(t,\mathbf x,v)$:
\begin{equation}\label{eq-general-time}
\left\{
\begin{aligned}
&\partial_t h+v\cdot\nabla_{\mathbf x}h+\lambda(\mathbf x,v)h
=\mathcal S(t,\mathbf x,v),
\quad (t,\mathbf x,v)\in\mathbb R_+\times\Omega^{'}\times\mathbb R^3,\\
&h(t,\mathbf x,v)\big|_{\Gamma_-}
=\sqrt{2\pi}\sqrt\mu\mathcal J[h](t,\mathbf x)+r(t,\mathbf x,v),\\
&h(0,\mathbf x,v)
=h_0(\mathbf x,v).
\end{aligned}
\right.
\end{equation}
For $0\le s\le t$, the backward characteristic terminating at $(\mathbf x,v)$ at time $t$ is given by
\begin{equation}\label{def-general-char-time}
\mathbf X(s;t,\mathbf x,v)
=
\mathbf x-(t-s)v,
\quad
V(s;t,\mathbf x,v)
=
v.
\end{equation}
For $(\mathbf x,v)\in\Gamma_+$, write $t_b=t_b(\mathbf x,v)$. For $f\in L^\infty(\mathbb R_+\times\partial\Omega^{'})$, define the time-dependent boundary-to-boundary flux transmission operator by
\begin{align}
\mathscr B_\lambda^t f(t,\mathbf x)
:={}&
\sqrt{2\pi}
\int_{n(\mathbf x)\cdot v>0}
\mathbf 1_{\{t_b<t\}}
e^{-\int_{t-t_b}^t\lambda(\mathbf X(s;t,\mathbf x,v),v)\,ds}
\mu\big[n(\mathbf x)\cdot v\big]
\notag\\
&\quad\times
f\big(
t-t_b,
\mathbf X(t-t_b;t,\mathbf x,v)
\big)\,dv.
\label{def-B-lambda-time}
\end{align}
We also define
\begin{align}
\mathcal R_\lambda(t,\mathbf x)
={}&
\int_{n(\mathbf x)\cdot v>0}
\int_{\max\{0,t-t_b\}}^t
e^{-\int_\tau^t\lambda(\mathbf X(s;t,\mathbf x,v),v)\,ds}
\mathcal S\big(\tau,\mathbf X(\tau;t,\mathbf x,v),v\big)
\sqrt\mu\big[n(\mathbf x)\cdot v\big]\,d\tau dv
\notag\\
&+
\int_{n(\mathbf x)\cdot v>0}
\mathbf 1_{\{t_b\ge t\}}
e^{-\int_0^t\lambda(\mathbf X(s;t,\mathbf x,v),v)\,ds}
h_0\big(\mathbf X(0;t,\mathbf x,v),v\big)
\sqrt\mu\big[n(\mathbf x)\cdot v\big]\,dv
\notag\\
&+
\int_{n(\mathbf x)\cdot v>0}
\mathbf 1_{\{t_b<t\}}
e^{-\int_{t-t_b}^t\lambda(\mathbf X(s;t,\mathbf x,v),v)\,ds}
r\big(t-t_b,\mathbf X(t-t_b;t,\mathbf x,v),v\big)
\sqrt\mu\big[n(\mathbf x)\cdot v\big]\,dv.
\label{def-R-lambda-time}
\end{align}

The following proposition gives a closed time-dependent boundary flux equation for the outgoing mass flux in a general bounded $C^1$ domain.

\begin{proposition}\label{prop-general-time-flux}
Let $\lambda=\lambda(\mathbf x,v)>0$ satisfy, for any $0<\varepsilon<1$ as in Lemma~\ref{lem-general-ft},
\[
\lambda_\varepsilon:=\inf_{\mathbf x\in\Omega^{'},\,|v|\le\varepsilon^{-1}}\lambda(\mathbf x,v)>0.
\]
Let $T>0$ and let $h$ be a mild solution of \eqref{eq-general-time}. Then the outgoing mass flux $\mathcal J[h]$ satisfies
\begin{equation}\label{eq-general-flux-time}
\mathcal J[h]
=
\mathscr B_\lambda^t\mathcal J[h]+\mathcal R_\lambda,
\quad
(t,\mathbf x)\in[0,T]\times\partial\Omega^{'}.
\end{equation}
Moreover, with the same constant $\kappa_0\in(0,1)$ as in Proposition~\ref{prop-general-steady-flux},
\begin{equation}\label{est-B-contraction-time}
\|\mathscr B_\lambda^t\|_{L^\infty([0,T]\times\partial\Omega^{'})\to L^\infty([0,T]\times\partial\Omega^{'})}
\le
\kappa_0.
\end{equation}
Consequently, $I-\mathscr B_\lambda^t$ is invertible on $L^\infty([0,T]\times\partial\Omega^{'})$, and
\begin{equation}\label{eq-flux-iteration-time}
\mathcal J[h]
=
\sum_{j=0}^\infty(\mathscr B_\lambda^t)^j\mathcal R_\lambda.
\end{equation}
In particular,
\begin{align}
\|\mathcal J[h]\|_{L^\infty([0,T]\times\partial\Omega^{'})}
\le
C_{\Omega^{'},\kappa_0}
\Big[
&\|\mathcal S\|_{L^\infty([0,T]\times\Omega^{'}\times\mathbb R^3)}
+
|r|_{L^\infty([0,T]\times\Gamma_-)}
+
\|h_0\|_{L^\infty(\Omega^{'}\times\mathbb R^3)}
\Big].
\label{est-general-time-flux}
\end{align}
\end{proposition}

\begin{proof}
The proof follows the same argument as that of Proposition~\ref{prop-general-steady-flux}. For $(t,\mathbf x,v)\in[0,T]\times\Gamma_+$,
\[
\begin{aligned}
h(t,\mathbf x,v)
={}&
\mathbf 1_{\{t_b<t\}}
e^{-\int_{t-t_b}^t\lambda(\mathbf X(s;t,\mathbf x,v),v)\,ds}
\sqrt{2\pi}\sqrt\mu\mathcal J[h]
\big(t-t_b,\mathbf X(t-t_b;t,\mathbf x,v)\big)
\\
&
+\mathbf 1_{\{t_b<t\}}
e^{-\int_{t-t_b}^t\lambda(\mathbf X(s;t,\mathbf x,v),v)\,ds}
r\big(t-t_b,\mathbf X(t-t_b;t,\mathbf x,v),v\big)
\\
&+
\mathbf 1_{\{t_b\ge t\}}
e^{-\int_0^t\lambda(\mathbf X(s;t,\mathbf x,v),v)\,ds}
h_0\big(\mathbf X(0;t,\mathbf x,v),v\big)
\\
&+
\int_{\max\{0,t-t_b\}}^t
e^{-\int_\tau^t\lambda(\mathbf X(s;t,\mathbf x,v),v)\,ds}
\mathcal S\big(\tau,\mathbf X(\tau;t,\mathbf x,v),v\big)\,d\tau.
\end{aligned}
\]
Multiplying by $\sqrt\mu\big[n(\mathbf x)\cdot v\big]$ and integrating over $n(\mathbf x)\cdot v>0$ gives \eqref{eq-general-flux-time}.

For the contraction estimate, we decompose $\Gamma_+$ into $\Gamma_+^\varepsilon$ defined in \eqref{def-Gamma-eps} and $\Gamma_+\setminus\Gamma_+^\varepsilon$, and use $\delta_\varepsilon$ defined in \eqref{def-B-delta}. Since the factor $\mathbf 1_{\{t_b<t\}}$ only restricts the boundary-to-boundary characteristics contributing to $\mathscr B_\lambda^t$, the same argument as in the proof of Proposition~\ref{prop-general-steady-flux} yields
\[
\|\mathscr B_\lambda^t f\|_{L^\infty([0,T]\times\partial\Omega^{'})}
\le
\kappa_0
\|f\|_{L^\infty([0,T]\times\partial\Omega^{'})},
\quad
\forall f\in L^\infty([0,T]\times\partial\Omega^{'}).
\]
where $\kappa_0\in(0,1)$ is defined in \eqref{def-B-kappa}.

Furthermore, since $\lambda>0$, the attenuation factors in \eqref{def-R-lambda-time} are bounded by one. Since $\Omega^{'}$ is bounded,
\[
\min\{t,t_b(\mathbf x,v)\}\big[n(\mathbf x)\cdot v\big]
\le
t_b(\mathbf x,v)\big[n(\mathbf x)\cdot v\big]
\le
C_{\Omega^{'}}.
\]
It follows from \eqref{def-R-lambda-time} that
\[
\|\mathcal R_\lambda\|_{L^\infty([0,T]\times\partial\Omega^{'})}
\le
C_{\Omega^{'}}
\Big[
\|\mathcal S\|_{L^\infty([0,T]\times\Omega^{'}\times\mathbb R^3)}
+
|r|_{L^\infty([0,T]\times\Gamma_-)}
+
\|h_0\|_{L^\infty(\Omega^{'}\times\mathbb R^3)}
\Big].
\]
The remainder of the proof is identical to that of Proposition~\ref{prop-general-steady-flux}. This completes the proof of Proposition~\ref{prop-general-time-flux}.
\end{proof}

The strict contraction property established above for the boundary-to-boundary flux transmission operators defined in \eqref{def-B-lambda} and \eqref{def-B-lambda-time} is not restricted to the $L^\infty$ setting. In fact, both operators remain strictly contractive throughout the full range $1\le p\le\infty$, as stated in the following corollary.

\begin{corollary}\label{cor-B-Lp}
Under the assumptions of Propositions~\ref{prop-general-steady-flux} and \ref{prop-general-time-flux}, for every $1\le p\le\infty$, the boundary-to-boundary flux transmission operator satisfies
\begin{equation}\label{est-B-Lp}
\|\mathscr B_\lambda\|_{L^p(\partial\Omega^{'})\to L^p(\partial\Omega^{'})}
\le
\kappa_0.
\end{equation}
Moreover, for every $T>0$,
\begin{equation}\label{est-B-time-Lp}
\|\mathscr B_\lambda^t\|_{L^p((0,T)\times\partial\Omega^{'})\to L^p((0,T)\times\partial\Omega^{'})}
\le
\kappa_0.
\end{equation}
Consequently, $I-\mathscr B_\lambda$ is invertible on $L^p(\partial\Omega^{'})$, and $I-\mathscr B_\lambda^t$ is invertible on $L^p((0,T)\times\partial\Omega^{'})$.
\end{corollary}

\begin{proof}
The case $p=\infty$ follows from Propositions~\ref{prop-general-steady-flux} and \ref{prop-general-time-flux}. We first prove the steady $L^1$ estimate.

For $(\mathbf x,v)\in\Gamma_+$, set
\[
\mathbf x'
:=
\mathbf X(t_b(\mathbf x,v);\mathbf x,v)
=
\mathbf x-t_b(\mathbf x,v)v.
\]
Then $(\mathbf x',v)\in\Gamma_-$. Define the corresponding forward exit time by $t_f(\mathbf x',v):=\inf\{s>0:\mathbf x'+sv\in\partial\Omega^{'}\}$. The boundary change of variables $\mathbf x\mapsto\mathbf x'$ satisfies
\begin{equation}\label{eq-boundary-change}
|n(\mathbf x)\cdot v|\,dS_{\mathbf x}
=
|n(\mathbf x')\cdot v|\,dS_{\mathbf x'}.
\end{equation}

For $f\in L^1(\partial\Omega^{'})$, it follows from \eqref{def-B-lambda} and \eqref{eq-boundary-change} that
\begin{align}
\|\mathscr B_\lambda f\|_{L^1(\partial\Omega^{'})}
\le{}&
\sqrt{2\pi}
\int_{\partial\Omega^{'}}
\int_{n(\mathbf x)\cdot v>0}
e^{-\int_0^{t_b(\mathbf x,v)}
\lambda(\mathbf x-\tau v,v)\,d\tau}
\mu\big[n(\mathbf x)\cdot v\big]
\big|f\big(\mathbf x-t_b(\mathbf x,v)v\big)\big|\,dv dS_{\mathbf x}
\notag\\
={}&
\int_{\partial\Omega^{'}}
|f(\mathbf x')|
\underbrace{\bigg[
\sqrt{2\pi}
\int_{n(\mathbf x')\cdot v<0}
e^{-\int_0^{t_f(\mathbf x',v)}
\lambda(\mathbf x'+\tau v,v)\,d\tau}
\mu|n(\mathbf x')\cdot v|\,dv
\bigg]}_{I(\mathbf x')}
\,dS_{\mathbf x'}.
\label{est-B-L1-reduction}
\end{align}
Setting $u=-v$ and noting that
$t_f(\mathbf x',-u)=t_b(\mathbf x',u)$, we obtain
\begin{equation}\label{eq-I-outgoing}
I(\mathbf x')
=
\sqrt{2\pi}
\int_{n(\mathbf x')\cdot u>0}
e^{-\int_0^{t_b(\mathbf x',u)}
\lambda(\mathbf x'-\tau u,-u)\,d\tau}
\mu\big[n(\mathbf x')\cdot u\big]\,du.
\end{equation}
We split the velocity integral in \eqref{eq-I-outgoing} over $\Gamma_+^{\varepsilon_0}$ and $\Gamma_+\setminus\Gamma_+^{\varepsilon_0}$. On $\Gamma_+^{\varepsilon_0}$, Lemma~\ref{lem-general-ft} gives
\[
t_b(\mathbf x',u)
\ge
\tau_{\varepsilon_0,\Omega^{'}}>0,
\]
and the definition of $\lambda_{\varepsilon_0}$ yields
\[
e^{-\int_0^{t_b(\mathbf x',u)}
\lambda(\mathbf x'-\tau u,-u)\,d\tau}
\le
e^{-\lambda_{\varepsilon_0}\tau_{\varepsilon_0,\Omega^{'}}}.
\]
Hence, using \eqref{def-B-delta}, and \eqref{def-B-kappa}, we obtain from \eqref{est-B-L1-reduction} that
\begin{equation}\label{est-B-L1}
\|\mathscr B_\lambda f\|_{L^1(\partial\Omega^{'})}
\le
\kappa_0\|f\|_{L^1(\partial\Omega^{'})},
\quad
\forall f\in L^1(\partial\Omega^{'}).
\end{equation}
where $\kappa_0\in(0,1)$ is defined by \eqref{def-B-kappa}. Combining \eqref{est-B-L1} with \eqref{est-B-contraction}, the Riesz--Thorin interpolation theorem yields \eqref{est-B-Lp}.

The time-dependent $L^1$ estimate follows in the same way, and we omit the details. Finally, since $\kappa_0<1$, the stated invertibility of $I-\mathscr B_\lambda$ and $I-\mathscr B_\lambda^t$ follows from the Neumann series. This completes the proof of Corollary~\ref{cor-B-Lp}.
\end{proof}

\section{Linear problem for the steady corrections}\label{sec-4}

In this section, we consider the following linear steady boundary value
problem associated with the three steady correction equations
\eqref{eq-Ga}, \eqref{eq-Gd}, and \eqref{eq-GR}:
\begin{equation}\label{eq-Gd-linear}
\left\{
\begin{aligned}
&v\cdot\nabla_{\mathbf x}G+\nu(v) G
=
KG+\mathfrak F,
\\
&G(x,\pm1,z,v)\big|_{\pm v_y<0}
=
\Phi(x,\pm1,z,v)\mathcal J[G](x,\pm1,z)
+r(x,\pm1,z,v),
\end{aligned}
\right.
\end{equation}
with the zero-mass condition
\begin{equation}\label{eq-Gd-linear-mass}
\int_\Omega\int_{\mathbb R^3}
\sqrt\mu G(\mathbf x,v)\,dv d\mathbf x=0.
\end{equation}
Here we take $\Phi=\sqrt{2\pi}\sqrt\mu$ for \eqref{eq-Ga} and
\eqref{eq-Gd}, and
$\Phi=\mu^{-\frac12}\mathcal M_{T(x,\pm1,z)}(v_x-U(\pm1),v_y,v_z)$
for \eqref{eq-GR}. Thus,
\begin{equation}\label{eq-Gd-Phi-norm}
\int_{\pm v_y<0}
\Phi(x,\pm1,z,v)\sqrt\mu|v_y|\,dv=1,
\qquad (x,z)\in\mathbb T_x\times\mathbb T_z.
\end{equation}
Moreover, by \eqref{eq-MT-exp},
\begin{equation}\label{eq-Gd-Phi-small}
\big|w_l(\Phi-\sqrt{2\pi}\sqrt\mu)\big|_{L^\infty(\gamma_-)}
\le C_l(\alpha+\delta).
\end{equation}
Let $\mathfrak F=\mathfrak F(\mathbf x,v)$ and
$r=r(x,\pm1,z,v)$ be prescribed interior and incoming boundary
sources, respectively, satisfying the mass compatibility
conditions
\begin{equation}\label{eq-Gd-source-mass}
\left\{
\begin{aligned}
&\int_\Omega\int_{\mathbb R^3}\sqrt\mu\mathfrak F\,dv d\mathbf x=0,
\\
&\int_{\pm v_y<0}r(x,\pm1,z,v)\sqrt\mu|v_y|\,dv=0.
\end{aligned}
\right.
\end{equation}

For $l>1$, set
\[
\mathcal X_l
:=
\Bigg\{
\mathsf g=\mathsf g(\mathbf x,v)\ \Bigg|\
\|w_l\mathsf g\|_{L^\infty_{\mathbf x,v}}<\infty,\quad
\int_\Omega\int_{\mathbb R^3}
\sqrt\mu\mathsf g(\mathbf x,v)\,dv d\mathbf x=0
\Bigg\}.
\]

\begin{lemma}\label{lem-Gd-linear}
Let $l>1$, and let $\alpha,\delta>0$ be sufficiently small.
Assume that $\mathfrak F$ and $r$ satisfy the mass compatibility
conditions \eqref{eq-Gd-source-mass}. Then the boundary value problem \eqref{eq-Gd-linear} admits a unique solution $G=G(\mathbf x,v)\in\mathcal X_l$. Moreover,
\begin{equation}\label{est-Linf-Gd-linear}
\|w_lG\|_{L^\infty_{\mathbf x,v}}
\le
C_l\Big[
\|\nu^{-1}w_l\mathfrak F\|_{L^\infty_{\mathbf x,v}}
+
|w_lr|_{L^\infty(\gamma_-)}
\Big],
\end{equation}
where $C_l>0$ is independent of $\alpha$ and $\delta$.
\end{lemma}

\subsection{Uniform a priori estimates}

To solve \eqref{eq-Gd-linear}, let $0<\epsilon\le1$ and
$0\le\sigma\le1$. We consider the following approximate boundary value problem with \eqref{eq-Gd-source-mass} for $G=G(\mathbf x,v)$:
\begin{equation}\label{eq-Gd-eps-sigma}
\left\{
\begin{aligned}
&\epsilon\nu(v)G
+v\cdot\nabla_{\mathbf x}G
+\nu(v)G
=
\sigma KG+\mathfrak F,
\\
&G(x,\pm1,z,v)\big|_{\pm v_y<0}
=
\Phi(x,\pm1,z,v)\mathcal J[G](x,\pm1,z)
+r(x,\pm1,z,v).
\end{aligned}
\right.
\end{equation}
Moreover, by \eqref{eq-Gd-Phi-norm}, \eqref{eq-Gd-source-mass} and the boundary condition in \eqref{eq-Gd-eps-sigma}, we have
\begin{align}
\int_{\mathbb R^3}v_yG(x,1,z,v)\sqrt\mu\,dv
={}&
\int_{v_y>0}G(x,1,z,v)\sqrt\mu|v_y|\,dv
-\int_{v_y<0}G(x,1,z,v)\sqrt\mu|v_y|\,dv
\notag\\
={}&
\mathcal J[G](x,1,z)
-\int_{v_y<0}
\Big[\Phi(x,1,z,v)\mathcal J[G](x,1,z)+r(x,1,z,v)\Big]
\sqrt\mu|v_y|\,dv
\notag\\
={}&
\mathcal J[G](x,1,z)
\bigg(1-\int_{v_y<0}\Phi(x,1,z,v)\sqrt\mu|v_y|\,dv\bigg)
-
\int_{v_y<0}r(x,1,z,v)\sqrt\mu|v_y|\,dv
\notag\\
={}&0.
\label{eq-Gd-by-bc}
\end{align}
The same argument applies at $y=-1$,
$\int_{\mathbb R^3}v_yG(x,-1,z,v)\sqrt\mu\,dv=0$.
Taking the $L^2(\Omega\times\mathbb R_v^3)$ inner product of the equation
in \eqref{eq-Gd-eps-sigma} with $\sqrt\mu$, and using
\eqref{eq-Gd-source-mass}--\eqref{eq-Gd-by-bc}, we obtain
\[
\big[\epsilon+(1-\sigma)\big]
\int_\Omega\int_{\mathbb R^3}\nu(v)\sqrt\mu G(\mathbf x,v)\,dv d\mathbf x=0.
\]
Since $\epsilon+(1-\sigma)>0$, every solution satisfies
\begin{equation}\label{eq-Gd-mass}
\int_\Omega\int_{\mathbb R^3}\nu(v)\sqrt\mu G(\mathbf x,v)\,dv d\mathbf x=0.
\end{equation}

We begin by deriving an a priori estimate uniform in $\epsilon$ and $\sigma$, assuming only sufficient regularity of the solution to justify the computations below.

\begin{lemma}[Uniform a priori estimate]\label{lem-Gd-apriori}
Let $G=G(\mathbf x,v)$ solve \eqref{eq-Gd-eps-sigma} under the mass
compatibility condition \eqref{eq-Gd-source-mass}.
Then, for any $l>1$, if $\alpha$ and $\delta$ sufficiently small,
\begin{equation}\label{est-Linf-Gd-apriori}
\|w_lG\|_{L^\infty_{\mathbf x,v}}
\le
C_l\Big[
\|\nu^{-1}w_l\mathfrak F\|_{L^\infty_{\mathbf x,v}}
+
|w_lr|_{L^\infty(\gamma_-)}
\Big],
\end{equation}
where $C_l>0$ is independent of $\epsilon$ and $\sigma$.
\end{lemma}

\begin{proof}
\noindent{\textbf{Step 1. Weighted $L^\infty$ Estimates.}}
Fix $l>1$. We define
\[
\mathcal G(\mathbf x,v):=w_l(v)G(\mathbf x,v).
\]
Multiplying \eqref{eq-Gd-eps-sigma} by $w_l$, we obtain
\begin{equation}\label{eq-Gd-weighted-m}
\left\{
\begin{aligned}
&\epsilon\nu(v)\mathcal G
+v\cdot\nabla_{\mathbf x}\mathcal G
+\nu(v)\mathcal G
=
\sigma w_lK\big(w_l^{-1}\mathcal G\big)
+w_l\mathfrak F,
\\
&\mathcal G(x,\pm1,z,v)\big|_{\pm v_y<0}
=
w_l\Phi(x,\pm1,z,v)\mathcal J[G](x,\pm1,z)
+w_lr(x,\pm1,z,v).
\end{aligned}
\right.
\end{equation}

Applying Duhamel's formula to \eqref{eq-Gd-weighted-m} yields
\begin{align}
\mathcal G(\mathbf x,v)
={}&
e^{-(1+\epsilon)\nu(v)t_b}
\Big[
w_l\Phi\mathcal J[G]+w_lr
\Big]\big(\mathbf X(t_b),v\big)
\notag\\
&+
\sigma\int_0^{t_b}e^{-(1+\epsilon)\nu(v)s}
\Big[w_lK\big(w_l^{-1}\mathcal G\big)\Big]
\big(\mathbf X(s),v\big)\,ds
\notag\\
&+
\int_0^{t_b}e^{-(1+\epsilon)\nu(v)s}
(w_l\mathfrak F)\big(\mathbf X(s),v\big)\,ds
\notag\\
=:{}&
\sum_{i=1}^3I_i(\mathbf x,v).
\label{eq-Gd-Duhamel}
\end{align}
Here $\mathbf X(s):=\mathbf x-sv$, and the backward exit time defined in \eqref{def-wall-tb} is given by
\begin{equation}\label{def-tby}
t_b=t_b(\mathbf x,v)
=
\left\{
\begin{aligned}
&\frac{y+1}{v_y}, &&v_y>0,\\
&\frac{1-y}{|v_y|}, &&v_y<0.
\end{aligned}
\right.
\end{equation}

We first estimate $I_3$. Since
$\int_0^{t_b}e^{-(1+\epsilon)\nu(v)s}\,ds\le\nu^{-1}(v)$, we have
\begin{equation}\label{est-Gd-I4}
|I_3(\mathbf x,v)|
\le
\big\|\nu^{-1}w_l\mathfrak F\big\|_{L^\infty_{\mathbf x,v}}.
\end{equation}

We next estimate the boundary term $I_1$. To apply
Lemma~\ref{lem-wall-flux} to $G$, we take
\[
\lambda(v)=(1+\epsilon)\nu(v),
\qquad
\mathcal S=\sigma KG+\mathfrak F.
\]
Consequently, by \eqref{eq-Gd-Phi-small}, for $\alpha$ and $\delta$ sufficiently small, we have
\begin{align}
\max_{\pm}
\|\mathscr B_\lambda(\pm1)\|_{L^\infty(\partial\Omega)\to L^\infty(\partial\Omega)}
&\le
\max_{\pm}\operatorname*{ess\,sup}_{x,z}
\int_{\pm u_y<0}
e^{-\frac{2(1+\epsilon)\nu(u)}{|u_y|}}
\Phi\bigg(x-\frac{2u_x}{|u_y|},\pm1,z-\frac{2u_z}{|u_y|},u\bigg)
\sqrt\mu|u_y|\,du
\notag\\
&\le
\max_{\pm}\int_{\pm u_y<0}
e^{-\frac{2(1+\epsilon)\nu(u)}{|u_y|}}
\Big[\sqrt{2\pi}\mu+C_l(\alpha+\delta)w_l^{-1}\sqrt\mu\Big]
|u_y|\,du
\notag\\
&\le
\sqrt{2\pi}\int_{u_y>0}
e^{-\frac{2\nu(u)}{u_y}}\mu u_y\,du
+C_l(\alpha+\delta)\int_{u_y>0}w_l^{-1}\sqrt\mu u_y\,du
\notag\\
&\le
\sqrt{2\pi}\int_{u_y>0}
e^{-\frac{2\nu(u)}{u_y}}\mu u_y\,du
+C_l(\alpha+\delta)
<1.
\label{est-Gd-wall-B}
\end{align}
Here we used the fact that,
\[
\sqrt{2\pi}\int_{u_y>0} e^{-\frac{2\nu(u)}{u_y}}\mu u_y\,du
<
\sqrt{2\pi}\int_{u_y>0}\mu u_y\,du
=1.
\]
Thus \eqref{est-wall-B} holds uniformly in $\epsilon$ and $\sigma$. Using \eqref{eq-wall-flux-form} in Lemma~\ref{lem-wall-flux}, we obtain
\begin{equation}\label{eq-Gd-wall-flux-formulas}
\left\{
\begin{aligned}
\mathcal J[G](x,1,z)
={}&
\big[I-\mathscr B_\lambda(-1)\mathscr B_\lambda(1)\big]^{-1}
\big[\mathcal R_{\epsilon,+}(x,z)
+\mathscr B_\lambda(-1)\mathcal R_{\epsilon,-}(x,z)\big],
\\
\mathcal J[G](x,-1,z)
={}&
\big[I-\mathscr B_\lambda(1)\mathscr B_\lambda(-1)\big]^{-1}
\big[\mathcal R_{\epsilon,-}(x,z)
+\mathscr B_\lambda(1)\mathcal R_{\epsilon,+}(x,z)\big],
\end{aligned}
\right.
\end{equation}
where, with $t_b=2/|u_y|$,
\begin{align}
\mathcal R_{\epsilon,\pm}(x,z)
={}&
\int_{\pm u_y>0}e^{-(1+\epsilon)\nu(u)t_b}
r(x-u_xt_b,\mp1,z-u_zt_b,u)
\sqrt\mu|u_y|\,du
\notag\\
&+
\sigma\int_{\pm u_y>0}\int_0^{t_b}e^{-(1+\epsilon)\nu(u)s}
KG(x-u_xs,\pm1-u_ys,z-u_zs,u)
\sqrt\mu|u_y|\,dsdu
\notag\\
&+
\int_{\pm u_y>0}\int_0^{t_b}e^{-(1+\epsilon)\nu(u)s}
\mathfrak F(x-u_xs,\pm1-u_ys,z-u_zs,u)
\sqrt\mu|u_y|\,dsdu.
\label{def-Gd-R-eps-pm}
\end{align}
By \eqref{est-Gd-wall-B} and the Neumann series,
\begin{equation}\label{est-Gd-J-1}
\big|\mathcal J[G](x,\pm1,z)\big|
\le
C\big(
\|\mathcal R_{\epsilon,+}\|_{L^\infty_{x,z}}
+
\|\mathcal R_{\epsilon,-}\|_{L^\infty_{x,z}}
\big).
\end{equation}

To estimate \eqref{def-Gd-R-eps-pm}, decompose
\[
\mathcal R_{\epsilon,\pm}
=
\mathcal R_{\epsilon,\pm}^{r}
+
\mathcal R_{\epsilon,\pm}^{K}
+
\mathcal R_{\epsilon,\pm}^{\mathfrak F}.
\]
By the Cauchy--Schwarz inequality, we obtain
\begin{equation}\label{est-Gd-R-eta}
\|\mathcal R_{\epsilon,\pm}^{r}\|_{L^\infty_{x,z}}
\le
\bigg(
\int_{\pm u_y>0}
\big\|r(x,\mp1,z,u)\big\|_{L^\infty_{x,z}}^2
|u_y|\,du
\bigg)^{\frac12}
\bigg(
\int_{\pm u_y>0}\mu|u_y|\,du
\bigg)^{\frac12}
\lesssim
|w_lr|_{L^\infty(\gamma_-)}.
\end{equation}
Using the changes of variables
\begin{equation}\label{eq-Gd-wall-y-change}
\left\{
\begin{aligned}
&y'=1-u_ys,\quad u_y>0,
\\
&y'=-1+|u_y|s,\quad u_y<0,
\end{aligned}
\right.
\end{equation}
we obtain
\begin{equation}\label{est-Gd-R-F}
\|\mathcal R_{\epsilon,\pm}^{\mathfrak F}\|_{L^\infty_{x,z}}
\le
\|\nu^{-1}w_l\mathfrak F\|_{L^\infty_{\mathbf x,v}}
\int_{-1}^1\int_{\pm u_y>0}
e^{-\frac{(1+\epsilon)\nu(u)(1\mp y')}{|u_y|}}
\nu(u)w_l^{-1}\sqrt\mu\,dudy'
\lesssim
\|\nu^{-1}w_l\mathfrak F\|_{L^\infty_{\mathbf x,v}}.
\end{equation}

For $\mathcal R_{\epsilon,\pm}^{K}$, let
$\varepsilon_1>0$ be sufficiently small and $M_1>0$ sufficiently large.
We decompose
\begin{equation}\label{eq-Gd-RK-decomp}
\mathcal R_{\epsilon,\pm}^{K}
=
\sum_{j=0}^3\mathcal R_{\epsilon,\pm}^{K,j},
\end{equation}
where
\[
\left\{
\begin{aligned}
\mathcal R_{\epsilon,\pm}^{K,0}(x,z)
:={}&
\sigma\int_{\pm u_y>0}\int_0^{\min\{t_b,\varepsilon_1\}}
e^{-(1+\epsilon)\nu(u)s}\sqrt\mu
KG\big(
\underbrace{x-u_xs,\pm1-u_ys,z-u_zs}_{\mathbf X(s)},u
\big)|u_y|\,dsdu,
\\
\mathcal R_{\epsilon,\pm}^{K,1}(x,z)
:={}&
\sigma\int_{\pm u_y>0}
\int_{\min\{t_b,\varepsilon_1\}}^{t_b}e^{-(1+\epsilon)\nu(u)s}
\mathbf 1_{\{|u|\ge M_1\}}\sqrt\mu(u)
\int_{\mathbb R^3}\mathbf k(u,v_*)
G\big(\mathbf X(s),v_*\big)\,dv_*|u_y|\,dsdu,
\\
\mathcal R_{\epsilon,\pm}^{K,2}(x,z)
:={}&
\sigma\int_{\pm u_y>0}
\int_{\min\{t_b,\varepsilon_1\}}^{t_b}e^{-(1+\epsilon)\nu(u)s}
\mathbf 1_{\{|u|<M_1\}}\sqrt\mu(u)
\int_{|u-v_*|<M_1^{-1}}\mathbf k(u,v_*)
G\big(\mathbf X(s),v_*\big)\,dv_*|u_y|\,dsdu,
\\
\mathcal R_{\epsilon,\pm}^{K,3}(x,z)
:={}&
\sigma\int_{\pm u_y>0}
\int_{\min\{t_b,\varepsilon_1\}}^{t_b}e^{-(1+\epsilon)\nu(u)s}
\mathbf 1_{\{|u|<M_1\}}\sqrt\mu(u)
\int_{|u-v_*|\ge M_1^{-1}}\mathbf k(u,v_*)
G\big(\mathbf X(s),v_*\big)\,dv_*|u_y|\,dsdu.
\end{aligned}
\right.
\]
It follows directly from Lemma~\ref{lem-K-kernel} that
\begin{align}
\|\mathcal R_{\epsilon,\pm}^{K,0}\|_{L^\infty_{x,z}}
\le{}&
\varepsilon_1\|w_lG\|_{L^\infty_{\mathbf x,v}}
\int_{\pm u_y>0}\sqrt\mu(u)w_l^{-1}(u)|u_y|
\int_{\mathbb R^3}|\mathbf k_w(u,v_*)|\,dv_*du
\notag\\
\le{}&
C_l\varepsilon_1
\|w_lG\|_{L^\infty_{\mathbf x,v}}.
\label{est-Gd-RK-short}
\end{align}
Using \eqref{eq-Gd-wall-y-change} and
\eqref{est-kernel}--\eqref{est-kernel-weight} in
Lemma~\ref{lem-K-kernel}, we obtain
\begin{align}
&
\|\mathcal R_{\epsilon,\pm}^{K,1}\|_{L^\infty_{x,z}}
+
\|\mathcal R_{\epsilon,\pm}^{K,2}\|_{L^\infty_{x,z}}
\notag\\
&\le
C_l\|w_lG\|_{L^\infty_{\mathbf x,v}}
\int_{-1}^1\int_{\pm u_y>0}
\mathbf 1_{\{|u|\ge M_1\}}\sqrt\mu(u)w_l^{-1}(u)
\int_{\mathbb R^3}|\mathbf k_w(u,v_*)|\,dv_*dudy'
\notag\\
&\quad+
C_l\|w_lG\|_{L^\infty_{\mathbf x,v}}
\int_{-1}^1\int_{\pm u_y>0}
\mathbf 1_{\{|u|<M_1\}}\sqrt\mu(u)
\int_{|u-v_*|<M_1^{-1}}
|\mathbf k(u,v_*)|w_l^{-1}(v_*)\,dv_*dudy'
\notag\\
&\le
C_l e^{-\frac{M_1^2}{16}}
\|w_lG\|_{L^\infty_{\mathbf x,v}}
+
C_l\|w_lG\|_{L^\infty_{\mathbf x,v}}
\sup_{|u|<M_1}
\int_{|u-v_*|<M_1^{-1}}
|\mathbf k_w(u,v_*)|\,dv_*
\notag\\
&\le
C_l\bigg(e^{-\frac{M_1^2}{16}}
+\frac{1}{M_1^{\gamma+3}}
+\frac{1}{M_1^{\frac{\gamma+3}{2}}}\bigg)
\|w_lG\|_{L^\infty_{\mathbf x,v}}
\notag\\
&\le
\frac{C_l}{M_1^{\frac{\gamma+3}{2}}}
\|w_lG\|_{L^\infty_{\mathbf x,v}},
\qquad M_1\gg1.
\label{est-Gd-RK-small}
\end{align}
For $\mathcal R_{\epsilon,+}^{K,3}$, the support of
the integrand satisfies $|u|<M_1$ and $|u-v_*|\ge M_1^{-1}$.
Hence, \eqref{est-kernel} gives
\[
\sqrt\mu(u)|\mathbf k(u,v_*)|
\le
C_{M_1}e^{-\frac18|v_*|^2}.
\]
Using the above estimate and applying the Cauchy--Schwarz inequality
first in $v_*$ and then in $(s,u)$, we obtain
\begin{align}
\big|\mathcal R_{\epsilon,+}^{K,3}(x,z)\big|^2
\le{}&
C_{M_1}\bigg[
\int_{u_y>0}\mathbf 1_{\{t_b>\varepsilon_1\}}
\int_{\varepsilon_1}^{t_b}
e^{-(1+\epsilon)\nu(u)s}\mathbf 1_{\{|u|<M_1\}}
\Big(
\int_{\mathbb R^3}
\nu(v_*)\big|G\big(\mathbf X(s),v_*\big)\big|^2\,dv_*
\Big)^{\frac12}
u_y\,dsdu
\bigg]^2
\notag\\
\le{}&
C_{M_1}\bigg[
\int_{u_y>0}\mathbf 1_{\{t_b>\varepsilon_1\}}
\int_{\varepsilon_1}^{t_b}
e^{-(1+\epsilon)\nu(u)s}\mathbf 1_{\{|u|<M_1\}}
u_y\,dsdu
\bigg]
\notag\\
&\quad\times
\bigg[
\int_{u_y>0}\mathbf 1_{\{t_b>\varepsilon_1\}}
\int_{\varepsilon_1}^{t_b}
e^{-(1+\epsilon)\nu(u)s}\mathbf 1_{\{|u|<M_1\}}
u_y\int_{\mathbb R^3}
\nu(v_*)\big|G\big(\mathbf X(s),v_*\big)\big|^2
\,dv_*dsdu
\bigg]
\notag\\
\le{}&
C_{M_1}
\int_{u_y>0}\mathbf 1_{\{t_b>\varepsilon_1\}}
\int_{\varepsilon_1}^{t_b}
e^{-c_{M_1}s}\mathbf 1_{\{|u|<M_1\}}
u_y\int_{\mathbb R^3}
\nu(v_*)\big|G\big(\mathbf X(s),v_*\big)\big|^2
\,dv_*dsdu
\notag\\
={}&
C_{M_1}\int_{\varepsilon_1}^{\infty}
e^{-c_{M_1}s}
\int_{\mathbf E_s}\frac{1-y'}{s^4}
\int_{\mathbb R^3}
\nu(v_*)|G(\mathbf x',v_*)|^2
\,dv_*d\mathbf x'ds
\notag\\
\le{}&
C_{M_1}\int_{\varepsilon_1}^{\infty}
e^{-c_{M_1}s}\frac{(1+s)^2}{s^4}\,ds
\|\nu^{\frac12}G\|_{L^2_{\mathbf x,v}}^2
\notag\\
\le{}&
C_{\varepsilon_1,M_1}
\|\nu^{\frac12}G\|_{L^2_{\mathbf x,v}}^2.
\label{est-Gd-RK3-point}
\end{align}
Here we denoted $c_{M_1}:=\inf_{|u|\le M_1}\nu(u)>0$. For each $s\ge\varepsilon_1$, we used the change of variables
\[
(x',y',z')=\mathbf X(s)=(x-u_xs,1-u_ys,z-u_zs),
\quad
\bigg|
\det\frac{\partial(x',y',z')}{\partial(u_x,u_y,u_z)}
\bigg|=s^3,
\]
and denoted its image by
\[
\begin{aligned}
\mathbf E_s
&:=
\Big\{
(x-u_xs,1-u_ys,z-u_zs):
u_y>0,\ s<t_b,\ |u|<M_1
\Big\}
\\
&\subset
[x-M_1s,x+M_1s]\times(-1,1)
\times[z-M_1s,z+M_1s].
\end{aligned}
\]
The same estimate holds for $\mathcal R_{\epsilon,-}^{K,3}$.
Combining \eqref{est-Gd-RK-short}, \eqref{est-Gd-RK-small}, and
\eqref{est-Gd-RK3-point}, we obtain
\begin{equation}\label{est-Gd-RK}
\|\mathcal R_{\epsilon,\pm}^{K}\|_{L^\infty_{x,z}}
\le
C_l\Big(\varepsilon_1+\frac1{M_1^{\frac{\gamma+3}{2}}}\Big)
\|w_lG\|_{L^\infty_{\mathbf x,v}}
+
C_{\varepsilon_1,M_1}
\|\nu^{\frac12}G\|_{L^2_{\mathbf x,v}}.
\end{equation}

For any fixed $l>1$ and any $\varepsilon>0$, combining
\eqref{est-Gd-J-1}, \eqref{est-Gd-R-eta}, \eqref{est-Gd-R-F}, and
\eqref{est-Gd-RK}, and choosing $\varepsilon_1$ sufficiently small
and $M_1$ sufficiently large, we obtain
\begin{align}
\|\mathcal J[G]\|_{L^\infty(\partial\Omega)}
\le
\varepsilon\|w_lG\|_{L^\infty_{\mathbf x,v}}
+
C_{\varepsilon,l}\Big[
\|\nu^{\frac12}G\|_{L^2_{\mathbf x,v}}
+
\|\nu^{-1}w_l\mathfrak F\|_{L^\infty_{\mathbf x,v}}
+
|w_lr|_{L^\infty(\gamma_-)}
\Big].
\label{est-Gd-J}
\end{align}
By \eqref{eq-Gd-Phi-small} and \eqref{est-Gd-J}, we obtain,
for any $\varepsilon>0$,
\begin{align}
|I_1(\mathbf x,v)|
\le{}&
\varepsilon\|\mathcal G\|_{L^\infty_{\mathbf x,v}}
+
C_{\varepsilon,l}\Big[
\|\nu^{\frac12}G\|_{L^2_{\mathbf x,v}}
+
\|\nu^{-1}w_l\mathfrak F\|_{L^\infty_{\mathbf x,v}}
+
|w_lr|_{L^\infty(\gamma_-)}
\Big].
\label{est-Gd-I1}
\end{align}

Recalling the weighted kernel representation in Lemma~\ref{lem-K-kernel},
\[
w_lK\big(w_l^{-1}\mathcal G\big)(\mathbf x,v)
=
\int_{\mathbb R^3}\mathbf k_w(v,v_*)\mathcal G(\mathbf x,v_*)\,dv_*,
\quad
\mathbf k_w(v,v_*)
=
w_l(v)\mathbf k(v,v_*)w_l^{-1}(v_*).
\]
We next estimate $I_2$. For sufficiently large $M_2>0$ and $M_3\gg M_2$, we decompose
\begin{equation}\label{eq-Gd-I2-decomp}
I_2=\sum_{j=1}^4I_{2,j},
\end{equation}
where
\[
\left\{
\begin{aligned}
I_{2,1}
:={}&
\sigma\int_0^{t_b}e^{-(1+\epsilon)\nu(v)s}
\mathbf 1_{\{|v|\ge M_2\}}
\int_{\mathbb R^3}
\mathbf k_w(v,v_*)
\mathcal G\big(\mathbf X(s),v_*\big)\,dv_*ds,
\\
I_{2,2}
:={}&
\sigma\int_0^{t_b}e^{-(1+\epsilon)\nu(v)s}
\mathbf 1_{\{|v|<M_2\}}
\int_{|v_*|\ge2M_2}
\mathbf k_w(v,v_*)
\mathcal G\big(\mathbf X(s),v_*\big)\,dv_*ds,
\\
I_{2,3}
:={}&
\sigma\int_0^{t_b}e^{-(1+\epsilon)\nu(v)s}
\mathbf 1_{\{|v|<M_2\}}
\int_{\substack{|v_*|<2M_2\\|v-v_*|<M_3^{-1}}}
\mathbf k_w(v,v_*)
\mathcal G\big(\mathbf X(s),v_*\big)\,dv_*ds,
\\
I_{2,4}
:={}&
\sigma\int_0^{t_b}e^{-(1+\epsilon)\nu(v)s}
\mathbf 1_{\{|v|<M_2\}}
\int_{\substack{|v_*|<2M_2\\|v-v_*|\ge M_3^{-1}}}
\mathbf k_w(v,v_*)
\mathcal G\big(\mathbf X(s),v_*\big)\,dv_*ds.
\end{aligned}
\right.
\]
Similarly to \eqref{est-Gd-RK-small}, using
\eqref{est-kernel}--\eqref{est-kernel-weight}, we obtain
\begin{align}
\sum_{j=1}^3|I_{2,j}|
\lesssim{}&
\|\mathcal G\|_{L^\infty_{\mathbf x,v}}
\sup_{|v|\ge M_2}
\nu^{-1}(v)
\int_{\mathbb R^3}|\mathbf k_w(v,v_*)|\,dv_*
+
\|\mathcal G\|_{L^\infty_{\mathbf x,v}}
\sup_{|v|<M_2}
\nu^{-1}(v)
\int_{|v_*|\ge2M_2}|\mathbf k_w(v,v_*)|\,dv_*
\notag\\
&+
\|\mathcal G\|_{L^\infty_{\mathbf x,v}}
\sup_{|v|<M_2}
\nu^{-1}(v)
\int_{\substack{|v_*|<2M_2\\|v-v_*|<M_3^{-1}}}
|\mathbf k_w(v,v_*)|\,dv_*
\notag\\
\le{}&
C_l\|\mathcal G\|_{L^\infty_{\mathbf x,v}}
\bigg[
\frac1{M_2^2}
+
(1+M_2^2)^{\frac32}e^{-\frac1{16}M_2^2}
+
\frac{(1+M_2^2)^{\frac32}}{M_3^{\frac{\gamma+3}{2}}}
\bigg]
\notag\\
\le{}&
\frac{C_l}{M_2}
\|\mathcal G\|_{L^\infty_{\mathbf x,v}},
\qquad M_2\gg1,
\label{est-Gd-I2-small}
\end{align}
where we used
\[
w_l(v)w_l^{-1}(v_*)\le C_l
\quad\text{on}\quad
\big\{(v,v_*):|v-v_*|<M_3^{-1}\big\},
\]
and chose $M_3$ sufficiently large depending on $M_2$ such that
\[
\frac{(1+M_2^2)^{\frac32}}{M_3^{\frac{\gamma+3}{2}}}
\le\frac1{M_2}.
\]
Consequently, combining \eqref{eq-Gd-Duhamel}, \eqref{est-Gd-I4},
\eqref{est-Gd-I1}, and \eqref{est-Gd-I2-small}, we obtain
\begin{equation}\label{est-Gd-pointwise}
|\mathcal G(\mathbf x,v)|
\le
\mathcal R
+
|I_{2,4}|,
\end{equation}
where
\begin{align}
\mathcal R
:={}&
\Big(\varepsilon+\frac{C_l}{M_2}\Big)
\|\mathcal G\|_{L^\infty_{\mathbf x,v}}
+
C_{\varepsilon,l}\Big[
\|\nu^{\frac12}G\|_{L^2_{\mathbf x,v}}
+
\|\nu^{-1}w_l\mathfrak F\|_{L^\infty_{\mathbf x,v}}
+
|w_lr|_{L^\infty(\gamma_-)}
\Big].
\label{def-Gd-Rm}
\end{align}
For fixed $(s,v_*)$, let
\[
\widetilde{\mathbf X}(\widetilde s)
=
\mathbf X(s)-\widetilde s v_*,
\quad
\widetilde t_b
=
t_b\big(\mathbf X(s),v_*\big)
\]
denote the second backward characteristic starting from $(\mathbf X(s),v_*)$ and its backward exit time.
Applying the pointwise estimate \eqref{est-Gd-pointwise} at $(\mathbf X(s),v_*)$
and splitting the second characteristic parameter at $\widetilde s=\varepsilon_2$,
we obtain, for $0<\varepsilon_2\le\varepsilon$ sufficiently small and fixed,
\begin{align}
|\mathcal G(\mathbf x,v)|
\le{}&
C_l\mathcal R
+
C\int_0^{t_b}e^{-(1+\epsilon)\nu(v)s}
\mathbf 1_{\{|v|<M_2\}}
\int_{\substack{|v_*|<2M_2\\|v-v_*|\ge M_3^{-1}}}
|\mathbf k_w(v,v_*)|
\notag\\
&\quad\times
\int_0^{\min\{\varepsilon_2,\widetilde t_b\}}
e^{-(1+\epsilon)\nu(v_*)\widetilde s}
\int_{\mathbb R^3}
|\mathbf k_w(v_*,v_{**})|
\big|\mathcal G\big(\widetilde{\mathbf X}(\widetilde s),v_{**}\big)\big|
\,dv_{**}d\widetilde s dv_*ds
\notag\\
&+
C\int_0^{t_b}e^{-(1+\epsilon)\nu(v)s}
\mathbf 1_{\{|v|<M_2\}}
\int_{\substack{|v_*|<2M_2\\|v-v_*|\ge M_3^{-1}}}
|\mathbf k_w(v,v_*)|
\mathbf 1_{\{\widetilde t_b>\varepsilon_2\}}
\notag\\
&\quad\times
\int_{\varepsilon_2}^{\widetilde t_b}
e^{-(1+\epsilon)\nu(v_*)\widetilde s}
\mathbf 1_{\{|v_*|<M_2\}}
\int_{\substack{|v_{**}|<2M_2\\|v_*-v_{**}|\ge M_3^{-1}}}
|\mathbf k_w(v_*,v_{**})|
\notag\\
&\quad\times
\big|\mathcal G\big(\widetilde{\mathbf X}(\widetilde s),v_{**}\big)\big|
\,dv_{**}d\widetilde s dv_*ds
\notag\\
\le{}&
C_l\mathcal R
+
C_{M_2,M_3,l}\mathbf 1_{\{|v|<M_2\}}
\int_0^{t_b}e^{-(1+\epsilon)\nu(v)s}
\int_{|v_*|<2M_2}
\mathbf 1_{\{\widetilde t_b>\varepsilon_2\}}
\int_{\varepsilon_2}^{\widetilde t_b}
e^{-(1+\epsilon)\nu(v_*)\widetilde s}
\notag\\
&\quad\quad\times
\int_{|v_{**}|<2M_2}
\nu^{\frac12}(v_{**})
\big|G\big(\widetilde{\mathbf X}(\widetilde s),v_{**}\big)\big|
\,dv_{**}d\widetilde s dv_*ds.
\label{est-Gd-double-Duhamel}
\end{align}
Here we have used \eqref{est-kernel-weight} and, on the corresponding integration regions,
\[
|\mathbf k_w(v,v_*)|
+
|\mathbf k_w(v_*,v_{**})|
\le C_{M_2,M_3,l},
\quad
w_l(v_{**})\nu^{-\frac12}(v_{**})\le C_{M_2,l}.
\]

For the last term in \eqref{est-Gd-double-Duhamel}, we make the change of variables
\begin{equation}\label{eq-cha-X}
\widetilde{\mathbf x}
=
\widetilde{\mathbf X}(\widetilde s)
=
\mathbf X(s)-\widetilde s v_*,
\quad
\bigg|
\det\frac{\partial\widetilde{\mathbf x}}{\partial v_*}
\bigg|
=
\widetilde s^3.
\end{equation}
Since $\widetilde s<\widetilde t_b$, we have $\widetilde y\in(-1,1)$. Denoting the image by
\[
\begin{aligned}
\mathbf E_{s,\widetilde s}
&:=
\Big\{
\mathbf X(s)-\widetilde s v_*:
|v_*|<2M_2,\ \widetilde s<\widetilde t_b
\Big\}
\\
&\subset
[x-v_xs-2M_2\widetilde s,x-v_xs+2M_2\widetilde s]
\times(-1,1)
\times
[z-v_zs-2M_2\widetilde s,z-v_zs+2M_2\widetilde s].
\end{aligned}
\]
Let $c_{2M_2}:=\inf_{|u|\le2M_2}\nu(u)>0$.
Therefore, the Cauchy--Schwarz inequality and the change of variables
\eqref{eq-cha-X} yield
\begin{align}
|\mathcal G(\mathbf x,v)|
\le{}&
C_l\mathcal R
+
C_{M_2,M_3,l}\int_0^{t_b}e^{-c_{2M_2}s}
\int_{\varepsilon_2}^\infty e^{-c_{2M_2}\widetilde s}
\notag\\
&\quad\times
\bigg[
\int_{|v_*|<2M_2}
\mathbf 1_{\{\widetilde s<\widetilde t_b\}}
\int_{|v_{**}|<2M_2}
\nu(v_{**})
\big|G\big(\widetilde{\mathbf X}(\widetilde s),v_{**}\big)\big|^2
\,dv_{**}dv_*
\bigg]^{\frac12}
\,d\widetilde s ds
\notag\\
={}&
C_l\mathcal R
+
C_{M_2,M_3,l}\int_0^{t_b}e^{-c_{2M_2}s}
\int_{\varepsilon_2}^\infty
\frac{e^{-c_{2M_2}\widetilde s}}{\widetilde s^{3/2}}
\notag\\
&\quad\times
\bigg[
\int_{\mathbf E_{s,\widetilde s}}
\int_{|v_{**}|<2M_2}
\nu(v_{**})
|G(\widetilde{\mathbf x},v_{**})|^2
\,dv_{**}d\widetilde{\mathbf x}
\bigg]^{\frac12}
\,d\widetilde s ds
\notag\\
\le{}&
C_l\mathcal R
+
C_{M_2,M_3,l}
\int_0^\infty e^{-c_{2M_2}s}\,ds
\int_{\varepsilon_2}^\infty
e^{-c_{2M_2}\widetilde s}
\frac{1+\widetilde s}{\widetilde s^{3/2}}\,d\widetilde s
\|\nu^{\frac12}G\|_{L^2_{\mathbf x,v}}
\notag\\
\le{}&
C_l\mathcal R
+
C_{\varepsilon_2,M_2,M_3,l}
\|\nu^{\frac12}G\|_{L^2_{\mathbf x,v}}.
\label{est-Gd-double-Duh-L2}
\end{align}

For any fixed $l>1$, choosing $1\ll M_2\ll M_3$ and
$\varepsilon>0$ sufficiently small, and substituting
\eqref{def-Gd-Rm} into \eqref{est-Gd-double-Duh-L2}, we obtain
\begin{align}
\|w_lG\|_{L^\infty_{\mathbf x,v}}
\le{}&
C_l\Big[
\|\nu^{\frac12}G\|_{L^2_{\mathbf x,v}}
+
\|\nu^{-1}w_l\mathfrak F\|_{L^\infty_{\mathbf x,v}}
+
|w_lr|_{L^\infty(\gamma_-)}
\Big],
\label{est-Linf-Gd-m}
\end{align}
where $C_l>0$ is independent of $\epsilon$ and $\sigma$.

\noindent{\textbf{Step 2. $L^2$ Estimates for $G$.}}
Taking the $L^2(\Omega\times\mathbb R^3_v)$ inner product of the first
equation in \eqref{eq-Gd-eps-sigma} with $G$, applying
Lemma~\ref{lem-L-coer}, we obtain
\begin{align}
&
(\epsilon+1-\sigma)\|\nu^{\frac12}G\|_{L^2_{\mathbf x,v}}^2
+
\sigma\lambda_0\|\nu^{\frac12}\mathbf P_1G\|_{L^2_{\mathbf x,v}}^2
+
\frac12|(I-P_\gamma)G|_{L^2(\gamma_+)}^2
\notag\\
&\le
\varepsilon\|\nu^{\frac12}G\|_{L^2_{\mathbf x,v}}^2
+
C_\varepsilon\|\nu^{-\frac12}\mathfrak F\|_{L^2_{\mathbf x,v}}^2
+
C_l(\alpha+\delta)^2\|\mathcal J[G]\|_{L^2(\partial\Omega)}^2
+
C|r|_{L^2(\gamma_-)}^2,
\label{est-L2-Gd-energy}
\end{align}
where $\lambda_0>0$ is the constant in Lemma~\ref{lem-L-coer}. We have used the identity
\[
\begin{aligned}
&
\int_{\mathbb T^2}\int_{\mathbb R^3}v_y
\Big[
|G(x,1,z,v)|^2
-
|G(x,-1,z,v)|^2
\Big]\,dvdxdz
\\
&=
\int_{\mathbb T^2}\bigg\{
\int_{v_y>0}v_y|G(x,1,z,v)|^2\,dv
-
\int_{v_y<0}
\big|\Phi(x,1,z,v)\mathcal J[G](x,1,z)+r(x,1,z,v)\big|^2|v_y|\,dv
\\
&\qquad+
\int_{v_y<0}|v_y||G(x,-1,z,v)|^2\,dv
-
\int_{v_y>0}\big|\Phi(x,-1,z,v)\mathcal J[G](x,-1,z)+r(x,-1,z,v)\big|^2v_y\,dv
\bigg\}\,dxdz
\\
&=
|G|_{L^2(\gamma_+)}^2-|P_\gamma G|_{L^2(\gamma_-)}^2-\big|(\Phi-\sqrt{2\pi}\sqrt\mu)\mathcal J[G]+r\big|_{L^2(\gamma_-)}^2
\\
&=
|G|_{L^2(\gamma_+)}^2-|P_\gamma G|_{L^2(\gamma_+)}^2
-
\big|(\Phi-\sqrt{2\pi}\sqrt\mu)\mathcal J[G]+r\big|_{L^2(\gamma_-)}^2
\\
&=
|(I-P_\gamma)G|_{L^2(\gamma_+)}^2
-
\big|(\Phi-\sqrt{2\pi}\sqrt\mu)\mathcal J[G]+r\big|_{L^2(\gamma_-)}^2.
\end{aligned}
\]
Here the cross term on $\gamma_-$ vanishes by
\eqref{eq-Gd-Phi-norm} and \eqref{eq-Gd-source-mass}.
where the second equality follows from \eqref{eq-Gd-Phi-norm} and \eqref{eq-Gd-source-mass}. Moreover, by \eqref{eq-Gd-Phi-small},
\[
\big|(\Phi-\sqrt{2\pi}\sqrt\mu)\mathcal J[G]+r\big|_{L^2(\gamma_-)}^2
\le
C_l(\alpha+\delta)^2\|\mathcal J[G]\|_{L^2(\partial\Omega)}^2
+
C|r|_{L^2(\gamma_-)}^2.
\]

We next introduce the macroscopic quantities associated with $G$.
Recalling \eqref{def-chi} and \eqref{def-P01}, we write
\[
\left\{
\begin{aligned}
\mathbf P_0G
&=\bigg[a+\mathbf b\cdot v
+c\frac{|v|^2-3}{\sqrt6}\bigg]\sqrt\mu,\quad
\mathbf b=(b_x,b_y,b_z),
\\
a&=\langle G,\chi_0\rangle_v,\quad
b_x=\langle G,\chi_1\rangle_v,
\quad
b_y=\langle G,\chi_2\rangle_v,
\\
b_z&=\langle G,\chi_3\rangle_v,
\quad
c=\langle G,\chi_4\rangle_v.
\end{aligned}
\right.
\]

\textbf{Estimate of $a$.}
Let $\phi_a=\phi_a(\mathbf x)$ be the unique zero-mean
solution of
\[
\left\{
\begin{aligned}
-\Delta_{\mathbf x}\phi_a(\mathbf x)
&=
a-\frac1{|\Omega|}\int_\Omega a\,d\mathbf x,
\quad \mathbf x\in\Omega,
\\
\partial_y\phi_a(\mathbf x)\big|_{\partial\Omega}
&=0,\quad
\int_\Omega\phi_a(\mathbf x)\,d\mathbf x=0.
\end{aligned}
\right.
\]
By the standard elliptic estimate and the trace theorem,
\begin{equation}\label{est-L2-Gd-phia}
\|\phi_a\|_{H^2_{\mathbf x}}
+
\|\partial_x\phi_a\|_{L^2(\partial\Omega)}
+
\|\partial_z\phi_a\|_{L^2(\partial\Omega)}
\lesssim
\|\Delta_{\mathbf x}\phi_a\|_{L^2_{\mathbf x}}.
\end{equation}
Inspired by \cite{EGKM13}, we choose the test function
$\psi_a$ as
\[
\psi_a(\mathbf x,v)
=
(10-|v|^2)\sqrt\mu
v\cdot\nabla_{\mathbf x}\phi_a(\mathbf x).
\]
Taking the $L^2(\Omega\times\mathbb R_v^3)$ inner product of the first
equation in \eqref{eq-Gd-eps-sigma} with $\psi_a$, we obtain
\begin{equation}\label{eq-Gd-test-a}
\langle v\cdot\nabla_{\mathbf x}G,\psi_a\rangle_{\mathbf x,v}
+
(\epsilon+1-\sigma)
\langle\nu(v) G,\psi_a\rangle_{\mathbf x,v}
+
\sigma\langle LG,\psi_a\rangle_{\mathbf x,v}
=
\langle\mathfrak F,\psi_a\rangle_{\mathbf x,v}.
\end{equation}
We estimate the terms in \eqref{eq-Gd-test-a} separately. For the transport term,
\begin{align}
\langle v\cdot\nabla_{\mathbf x}G,\psi_a\rangle_{\mathbf x,v}
={}&
5\|\Delta_{\mathbf x}\phi_a\|_{L^2_{\mathbf x}}^2
-
\langle\mathbf P_1G,v\cdot\nabla_{\mathbf x}\psi_a\rangle_{\mathbf x,v}
\notag\\
&+
\underbrace{
\int_{\gamma_+}G(x,\pm1,z,v)\psi_a(x,\pm1,z,v)\,d\gamma
-
\int_{\gamma_-}G(x,\pm1,z,v)\psi_a(x,\pm1,z,v)\,d\gamma
}_{\mathcal B_a}.
\label{eq-Gd-a-ts}
\end{align}
By the Neumann boundary condition for $\phi_a$,
\[
\psi_a(x,\pm1,z,v)
=
(10-|v|^2)\sqrt\mu
\Big[
v_x\partial_x\phi_a(x,\pm1,z)
+
v_z\partial_z\phi_a(x,\pm1,z)
\Big].
\]
On $\gamma_-$, using the boundary condition in
\eqref{eq-Gd-eps-sigma}, we obtain
\begin{align}
\int_{\gamma_-}G(x,\pm1,z,v)\psi_a(x,\pm1,z,v)\,d\gamma
={}&
\sum_{\pm}\int_{\mathbb T^2}
\partial_x\phi_a(x,\pm1,z)
\int_{\pm v_y<0}
\Big[
\Phi(x,\pm1,z,v)\mathcal J[G](x,\pm1,z)
\notag\\
&\quad+
r(x,\pm1,z,v)
\Big]
(10-|v|^2)\sqrt\mu v_x|v_y|\,dvdxdz
\notag\\
&+
\sum_{\pm}\int_{\mathbb T^2}
\partial_z\phi_a(x,\pm1,z)
\int_{\pm v_y<0}
\Big[
\Phi(x,\pm1,z,v)\mathcal J[G](x,\pm1,z)
\notag\\
&\quad+
r(x,\pm1,z,v)
\Big]
(10-|v|^2)\sqrt\mu v_z|v_y|\,dvdxdz
\notag\\
={}&
\int_{\gamma_-}
\Big[
\big(\Phi(x,\pm1,z,v)-\sqrt{2\pi}\sqrt\mu\big)
\mathcal J[G](x,\pm1,z)
+
r(x,\pm1,z,v)
\Big]
\notag\\
&\quad\times
\psi_a(x,\pm1,z,v)\,d\gamma,
\label{est-L2-G1d-incom}
\end{align}
where we have used
\[
\int_{\gamma_-}P_\gamma G(x,\pm1,z,v)
\psi_a(x,\pm1,z,v)\,d\gamma
=0
\]
by oddness in $v_x$ and $v_z$.
On $\gamma_+$, the same argument gives
\[
\int_{\gamma_+}P_\gamma G(x,\pm1,z,v)
\psi_a(x,\pm1,z,v)\,d\gamma
=
0.
\]
Consequently,
\[
\begin{aligned}
|\mathcal B_a|
&\le
|(I-P_\gamma)G|_{L^2(\gamma_+)}
\bigg[
\int_{\gamma_+}|\psi_a(x,\pm1,z,v)|^2\,d\gamma
\bigg]^{\frac12}
\\
&\quad+
|(\Phi-\sqrt{2\pi}\sqrt\mu)\mathcal J[G]+r|_{L^2(\gamma_-)}
\bigg[
\int_{\gamma_-}|\psi_a(x,\pm1,z,v)|^2\,d\gamma
\bigg]^{\frac12}
\\
&\lesssim
\Big[
|(I-P_\gamma)G|_{L^2(\gamma_+)}
+
|(\Phi-\sqrt{2\pi}\sqrt\mu)\mathcal J[G]+r|_{L^2(\gamma_-)}
\Big]
\\
&\quad\times
\bigg\{
\sum_{\pm}\int_{\mathbb T^2}
\Big[
|\partial_x\phi_a(x,\pm1,z)|^2
+
|\partial_z\phi_a(x,\pm1,z)|^2
\Big]\,dxdz
\bigg\}^{\frac12}.
\end{aligned}
\]
Moreover, by \eqref{eq-Gd-Phi-small} and \eqref{est-L2-Gd-phia}, we obtain
\begin{align}
&
|\mathcal B_a|
+
\big|\langle\mathbf P_1G,v\cdot\nabla_{\mathbf x}\psi_a\rangle_{\mathbf x,v}\big|
\notag\\
&\lesssim
\Big[
\|\nu^{\frac12}\mathbf P_1G\|_{L^2_{\mathbf x,v}}
+
|(I-P_\gamma)G|_{L^2(\gamma_+)}
+
(\alpha+\delta)\|\mathcal J[G]\|_{L^2(\partial\Omega)}
+
|r|_{L^2(\gamma_-)}
\Big]
\|\Delta_{\mathbf x}\phi_a\|_{L^2_{\mathbf x}}
\notag\\
&\le
\varepsilon\|\Delta_{\mathbf x}\phi_a\|_{L^2_{\mathbf x}}^2
+
C_\varepsilon\Big[
\|\nu^{\frac12}\mathbf P_1G\|_{L^2_{\mathbf x,v}}^2
+
|(I-P_\gamma)G|_{L^2(\gamma_+)}^2
\notag\\
&\qquad+
(\alpha+\delta)^2\|\mathcal J[G]\|_{L^2(\partial\Omega)}^2
+
|r|_{L^2(\gamma_-)}^2
\Big].
\label{est-L2-Gd-a-ts}
\end{align}

For the second term on the left-hand side of \eqref{eq-Gd-test-a}, noticed that $(\epsilon+1-\sigma)^2\le2(\epsilon+1-\sigma)$, we have
\begin{align}
(\epsilon+1-\sigma)
\big|\langle\nu(v) G,\psi_a\rangle_{\mathbf x,v}\big|
&\lesssim
(\epsilon+1-\sigma)
\|\nu^{\frac12}G\|_{L^2_{\mathbf x,v}}
\|\Delta_{\mathbf x}\phi_a\|_{L^2_{\mathbf x}}
\notag\\
&\le
\varepsilon\|\Delta_{\mathbf x}\phi_a\|_{L^2_{\mathbf x}}^2
+
C_\varepsilon(\epsilon+1-\sigma)
\|\nu^{\frac12}G\|_{L^2_{\mathbf x,v}}^2.
\label{est-L2-Gd-a-zs}
\end{align}

By \eqref{est-kernel-weight} with $l=0$, we have
\[
\begin{aligned}
\|\nu^{-\frac12}K\mathbf P_1G\|_{L^2_{\mathbf x,v}}^2
&\le
\int_\Omega\int_{\mathbb R^3}
\frac{1}{\nu(v)}
\bigg(\int_{\mathbb R^3}|\mathbf k(v,v_*)|\,dv_*\bigg)
\bigg(\int_{\mathbb R^3}
|\mathbf k(v,v_*)|\,|\mathbf P_1G(\mathbf x,v_*)|^2\,dv_*
\bigg)\,dv\,d\mathbf x
\\
&\lesssim
\int_\Omega\int_{\mathbb R^3}
|\mathbf P_1G(\mathbf x,v_*)|^2
\bigg(\int_{\mathbb R^3}|\mathbf k(v,v_*)|\,dv\bigg)\,dv_*d\mathbf x
\\
&\lesssim
\|\nu^{\frac12}\mathbf P_1G\|_{L^2_{\mathbf x,v}}^2.
\end{aligned}
\]
Then the remaining terms in \eqref{eq-Gd-test-a} satisfy
\begin{align}
&
\sigma\big|\langle LG,\psi_a\rangle_{\mathbf x,v}\big|
+
\big|\langle\mathfrak F,\psi_a\rangle_{\mathbf x,v}\big|
\notag\\
&\le
\Big(
\|\nu^{\frac12}\mathbf P_1G\|_{L^2_{\mathbf x,v}}
+
\|\nu^{-\frac12}K\mathbf P_1G\|_{L^2_{\mathbf x,v}}
+
\|\nu^{-\frac12}\mathfrak F\|_{L^2_{\mathbf x,v}}
\Big)
\|\nu^{\frac12}\psi_a\|_{L^2_{\mathbf x,v}}
\notag\\
&\lesssim
\Big(
\|\nu^{\frac12}\mathbf P_1G\|_{L^2_{\mathbf x,v}}
+
\|\nu^{-\frac12}\mathfrak F\|_{L^2_{\mathbf x,v}}
\Big)
\|\Delta_{\mathbf x}\phi_a\|_{L^2_{\mathbf x}}
\notag\\
&\le
\varepsilon\|\Delta_{\mathbf x}\phi_a\|_{L^2_{\mathbf x}}^2
+
C_\varepsilon\Big(
\|\nu^{\frac12}\mathbf P_1G\|_{L^2_{\mathbf x,v}}^2
+
\|\nu^{-\frac12}\mathfrak F\|_{L^2_{\mathbf x,v}}^2
\Big).
\label{est-L2-Gd-a-rem}
\end{align}

Combining \eqref{eq-Gd-test-a}--\eqref{est-L2-Gd-a-ts}
and \eqref{est-L2-Gd-a-zs}--\eqref{est-L2-Gd-a-rem},
and choosing $\varepsilon>0$ sufficiently small, we obtain
\begin{align}
\|\Delta_{\mathbf x}\phi_a\|_{L^2_{\mathbf x}}^2
\lesssim{}&
\|\nu^{\frac12}\mathbf P_1G\|_{L^2_{\mathbf x,v}}^2
+
(\epsilon+1-\sigma)\|\nu^{\frac12}G\|_{L^2_{\mathbf x,v}}^2
+
|(I-P_\gamma)G|_{L^2(\gamma_+)}^2
\notag\\
&+
\|\nu^{-\frac12}\mathfrak F\|_{L^2_{\mathbf x,v}}^2
+
(\alpha+\delta)^2\|\mathcal J[G]\|_{L^2(\partial\Omega)}^2
+
|r|_{L^2(\gamma_-)}^2.
\label{est-L2-Gd-a-phi}
\end{align}
On the other hand, the mass constraint \eqref{eq-Gd-mass} gives
\begin{align}
0
={}&
\int_{\Omega\times\mathbb R^3}
\nu(v)\sqrt\mu G\,dv d\mathbf x
\notag\\
={}&
\bigg(\int_{\mathbb R^3}\nu(v)\mu\,dv\bigg)
\int_\Omega a\,d\mathbf x
+
\bigg(\int_{\mathbb R^3}\nu(v)\mu\frac{|v|^2-3}{\sqrt6}\,dv\bigg)
\int_\Omega c\,d\mathbf x
+
\int_{\Omega\times\mathbb R^3}
\nu(v)\sqrt\mu\mathbf P_1G\,dv d\mathbf x.
\label{eq-Gd-mass-decomp}
\end{align}
Using \eqref{eq-Gd-mass-decomp}, together with the periodicity in
$x$ and $z$ and $\partial_y\phi_a(\mathbf x)\big|_{\partial\Omega}=0$, we obtain
\begin{align}
\|a\|_{L^2_{\mathbf x}}^2
={}&
\|\Delta_{\mathbf x}\phi_a\|_{L^2_{\mathbf x}}^2
+
\frac1{|\Omega|}\bigg|\int_\Omega a\,d\mathbf x\bigg|^2
\notag\\
\lesssim{}&
\|\Delta_{\mathbf x}\phi_a\|_{L^2_{\mathbf x}}^2
+
\|c\|_{L^2_{\mathbf x}}^2
+
\|\nu^{\frac12}\mathbf P_1G\|_{L^2_{\mathbf x,v}}^2.
\label{est-L2-Gd-a-mean}
\end{align}
Combining \eqref{est-L2-Gd-a-phi} and \eqref{est-L2-Gd-a-mean},
we obtain
\begin{align}
\|a\|_{L^2_{\mathbf x}}^2
\lesssim{}&
\|\nu^{\frac12}\mathbf P_1G\|_{L^2_{\mathbf x,v}}^2
+
(\epsilon+1-\sigma)\|\nu^{\frac12}G\|_{L^2_{\mathbf x,v}}^2
+
|(I-P_\gamma)G|_{L^2(\gamma_+)}^2
\notag\\
&+
\|\nu^{-\frac12}\mathfrak F\|_{L^2_{\mathbf x,v}}^2
+
(\alpha+\delta)^2\|\mathcal J[G]\|_{L^2(\partial\Omega)}^2
+
|r|_{L^2(\gamma_-)}^2
+
\|c\|_{L^2_{\mathbf x}}^2.
\label{est-L2-Gd-a}
\end{align}

\textbf{Estimate of $\mathbf b$.}
For $j=x,y,z$, let $\phi_j=\phi_j(\mathbf x)$ solve
\[
\left\{
\begin{aligned}
-\Delta_{\mathbf x}\phi_j(\mathbf x)
&=b_j,\quad \mathbf x\in\Omega,
\\
\phi_j(\mathbf x)\big|_{\partial\Omega}
&=0.
\end{aligned}
\right.
\]
Then,
\[
\|\phi_j\|_{H^2_{\mathbf x}}
+
\|\partial_y\phi_j\|_{L^2(\partial\Omega)}
\lesssim
\|b_j\|_{L^2_{\mathbf x}}.
\]
We choose the test function $\psi_b$ as
\[
\psi_b(\mathbf x,v)
=
\sum_{j=x,y,z}
\Big[
2\partial_j\phi_j
-
v_jv\cdot\nabla_{\mathbf x}\phi_j
\Big]\sqrt\mu.
\]
Taking the $L^2(\Omega\times\mathbb R_v^3)$ inner product of the first
equation in \eqref{eq-Gd-eps-sigma} with $\psi_b$, we obtain
\begin{align}
\|\mathbf b\|_{L^2_{\mathbf x}}^2
={}&
\underbrace{
\int_{\gamma_+}G(x,\pm1,z,v)\psi_b(x,\pm1,z,v)\,d\gamma
-
\int_{\gamma_-}G(x,\pm1,z,v)\psi_b(x,\pm1,z,v)\,d\gamma
}_{\mathcal B_b}
\notag\\
&-
\langle\mathbf P_1G,v\cdot\nabla_{\mathbf x}\psi_b\rangle_{\mathbf x,v}
+
(\epsilon+1-\sigma)\langle\nu G,\psi_b\rangle_{\mathbf x,v}
+
\sigma\langle LG,\psi_b\rangle_{\mathbf x,v}
-
\langle\mathfrak F,\psi_b\rangle_{\mathbf x,v}.
\label{eq-Gd-bj-test}
\end{align}
Since $\phi_j=0$ on $y=\pm1$, we have
\[
\psi_b(x,\pm1,z,v)
=
\Big[
-v_xv_y\partial_y\phi_x(x,\pm1,z)
+
(2-v_y^2)\partial_y\phi_y(x,\pm1,z)
-
v_zv_y\partial_y\phi_z(x,\pm1,z)
\Big]\sqrt\mu.
\]
Substituting this expression into the boundary integrals and computing directly, we obtain
\[
\begin{aligned}
\mathcal B_b
={}&
\int_{\gamma_+}(I-P_\gamma)G(x,\pm1,z,v)
\psi_b(x,\pm1,z,v)\,d\gamma
\\
&-
\int_{\gamma_-}
\Big[
\big(\Phi(x,\pm1,z,v)-\sqrt{2\pi}\sqrt\mu\big)
\mathcal J[G](x,\pm1,z)
+
r(x,\pm1,z,v)
\Big]
\psi_b(x,\pm1,z,v)\,d\gamma.
\end{aligned}
\]
Hence, by \eqref{eq-Gd-Phi-small},
\[
\begin{aligned}
|\mathcal B_b|
&\lesssim
\|\mathbf b\|_{L^2_{\mathbf x}}
\Big[
|(I-P_\gamma)G|_{L^2(\gamma_+)}
+
(\alpha+\delta)\|\mathcal J[G]\|_{L^2(\partial\Omega)}
+
|r|_{L^2(\gamma_-)}
\Big]
\\
&\le
\varepsilon\|\mathbf b\|_{L^2_{\mathbf x}}^2
+
C_\varepsilon\Big[
|(I-P_\gamma)G|_{L^2(\gamma_+)}^2
+
(\alpha+\delta)^2\|\mathcal J[G]\|_{L^2(\partial\Omega)}^2
+
|r|_{L^2(\gamma_-)}^2
\Big].
\end{aligned}
\]
Moreover, estimating the remaining terms in \eqref{eq-Gd-bj-test}
in the same way as \eqref{est-L2-Gd-a-zs}--\eqref{est-L2-Gd-a-rem}, we obtain
\begin{align}
&
(\epsilon+1-\sigma)
\big|\langle\nu G,\psi_b\rangle_{\mathbf x,v}\big|
+
\big|\langle\mathbf P_1G,v\cdot\nabla_{\mathbf x}\psi_b\rangle_{\mathbf x,v}\big|
+
\sigma\big|\langle LG,\psi_b\rangle_{\mathbf x,v}\big|
+
\big|\langle\mathfrak F,\psi_b\rangle_{\mathbf x,v}\big|
\notag\\
&\lesssim
\Big[
(\epsilon+1-\sigma)\|\nu^{\frac12}G\|_{L^2_{\mathbf x,v}}
+
\|\nu^{\frac12}\mathbf P_1G\|_{L^2_{\mathbf x,v}}
+
\|\nu^{-\frac12}\mathfrak F\|_{L^2_{\mathbf x,v}}
\Big]
\|\mathbf b\|_{L^2_{\mathbf x}}
\notag\\
&\le
\varepsilon\|\mathbf b\|_{L^2_{\mathbf x}}^2
+
C_\varepsilon\Big[
(\epsilon+1-\sigma)\|\nu^{\frac12}G\|_{L^2_{\mathbf x,v}}^2
+
\|\nu^{\frac12}\mathbf P_1G\|_{L^2_{\mathbf x,v}}^2
+
\|\nu^{-\frac12}\mathfrak F\|_{L^2_{\mathbf x,v}}^2
\Big].
\label{est-L2-Gd-bj-rem}
\end{align}
Combining \eqref{eq-Gd-bj-test}--\eqref{est-L2-Gd-bj-rem},
and choosing $\varepsilon>0$ sufficiently small, we obtain
\begin{align}
\|\mathbf b\|_{L^2_{\mathbf x}}^2
\lesssim{}&
\|\nu^{\frac12}\mathbf P_1G\|_{L^2_{\mathbf x,v}}^2
+
(\epsilon+1-\sigma)\|\nu^{\frac12}G\|_{L^2_{\mathbf x,v}}^2
+
|(I-P_\gamma)G|_{L^2(\gamma_+)}^2
\notag\\
&+
\|\nu^{-\frac12}\mathfrak F\|_{L^2_{\mathbf x,v}}^2
+
(\alpha+\delta)^2\|\mathcal J[G]\|_{L^2(\partial\Omega)}^2
+
|r|_{L^2(\gamma_-)}^2.
\label{est-L2-Gd-bxz}
\end{align}

\textbf{Estimate of $c$.}
Let $\phi_c=\phi_c(\mathbf x)$ solve
\[
\left\{
\begin{aligned}
-\Delta_{\mathbf x}\phi_c(\mathbf x)
&=c,\quad \mathbf x\in\Omega,
\\
\phi_c(\mathbf x)\big|_{\partial\Omega}
&=0.
\end{aligned}
\right.
\]
Hence,
\[
\|\phi_c\|_{H^2_{\mathbf x}}
+
\|\partial_y\phi_c\|_{L^2(\partial\Omega)}
\lesssim
\|c\|_{L^2_{\mathbf x}}.
\]
Let
\begin{equation}\label{def-Gd-psic}
\psi_c(\mathbf x,v)
=
\big(|v|^2-5\big)\sqrt\mu
v\cdot\nabla_{\mathbf x}\phi_c(\mathbf x).
\end{equation}
Taking the $L^2(\Omega\times\mathbb R_v^3)$ inner product of the first
equation in \eqref{eq-Gd-eps-sigma} with $\psi_c$, we obtain
\begin{align}
\frac{10}{\sqrt6}\|c\|_{L^2_{\mathbf x}}^2
&=
-\underbrace{
\int_{\gamma_+}G(x,\pm1,z,v)\psi_c(x,\pm1,z,v)\,d\gamma
-
\int_{\gamma_-}G(x,\pm1,z,v)\psi_c(x,\pm1,z,v)\,d\gamma
}_{\mathcal B_c}
\notag\\
&\quad+
\langle\mathbf P_1G,v\cdot\nabla_{\mathbf x}\psi_c\rangle_{\mathbf x,v}
-
(\epsilon+1-\sigma)\langle\nu(v)G,\psi_c\rangle_{\mathbf x,v}
-
\sigma\langle LG,\psi_c\rangle_{\mathbf x,v}
+
\langle\mathfrak F,\psi_c\rangle_{\mathbf x,v}.
\label{eq-Gd-test-c}
\end{align}
A direct computation yields
\begin{align}
\mathcal B_c
={}&
\int_{\mathbb T^2}\partial_y\phi_c(x,1,z)
\bigg\{
\int_{v_y>0}
\big(|v|^2-5\big)v_y^2\sqrt\mu
(I-P_\gamma)G(x,1,z,v)\,dv
\notag\\
&\quad+
\int_{v_y<0}
\big(|v|^2-5\big)v_y^2\sqrt\mu
\Big[
\big(\Phi(x,1,z,v)-\sqrt{2\pi}\sqrt\mu\big)
\mathcal J[G](x,1,z)
+
r(x,1,z,v)
\Big]\,dv
\bigg\}\,dxdz
\notag\\
&-
\int_{\mathbb T^2}\partial_y\phi_c(x,-1,z)
\bigg\{
\int_{v_y<0}
\big(|v|^2-5\big)v_y^2\sqrt\mu
(I-P_\gamma)G(x,-1,z,v)\,dv
\notag\\
&\quad+
\int_{v_y>0}
\big(|v|^2-5\big)v_y^2\sqrt\mu
\Big[
\big(\Phi(x,-1,z,v)-\sqrt{2\pi}\sqrt\mu\big)
\mathcal J[G](x,-1,z)
+
r(x,-1,z,v)
\Big]\,dv
\bigg\}\,dxdz.
\label{eq-Gd-c-bd}
\end{align}
Hence,
\begin{align}
|\mathcal B_c|
&\lesssim
\|c\|_{L^2_{\mathbf x}}
\Big[
|(I-P_\gamma)G|_{L^2(\gamma_+)}
+
(\alpha+\delta)\|\mathcal J[G]\|_{L^2(\partial\Omega)}
+
|r|_{L^2(\gamma_-)}
\Big]
\notag\\
&\le
\varepsilon\|c\|_{L^2_{\mathbf x}}^2
+
C_\varepsilon\Big[
|(I-P_\gamma)G|_{L^2(\gamma_+)}^2
+
(\alpha+\delta)^2\|\mathcal J[G]\|_{L^2(\partial\Omega)}^2
+
|r|_{L^2(\gamma_-)}^2
\Big].
\label{est-L2-Gd-c-bd}
\end{align}
Estimating the remaining terms in \eqref{eq-Gd-test-c} as in
\eqref{est-L2-Gd-a-zs} and \eqref{est-L2-Gd-a-rem} and taking
$\varepsilon>0$ sufficiently small, we obtain
\begin{align}
\|c\|_{L^2_{\mathbf x}}^2
\lesssim{}&
\|\nu^{\frac12}\mathbf P_1G\|_{L^2_{\mathbf x,v}}^2
+
(\epsilon+1-\sigma)\|\nu^{\frac12}G\|_{L^2_{\mathbf x,v}}^2
+
|(I-P_\gamma)G|_{L^2(\gamma_+)}^2
\notag\\
&+
\|\nu^{-\frac12}\mathfrak F\|_{L^2_{\mathbf x,v}}^2
+
(\alpha+\delta)^2\|\mathcal J[G]\|_{L^2(\partial\Omega)}^2
+
|r|_{L^2(\gamma_-)}^2.
\label{est-L2-Gd-c}
\end{align}

Therefore, \eqref{est-L2-Gd-a}, \eqref{est-L2-Gd-bxz},
and \eqref{est-L2-Gd-c} imply
\begin{align}
\|[a,\mathbf b,c]\|_{L^2_{\mathbf x}}^2
\lesssim{}&
\|\nu^{\frac12}\mathbf P_1G\|_{L^2_{\mathbf x,v}}^2
+
(\epsilon+1-\sigma)\|\nu^{\frac12}G\|_{L^2_{\mathbf x,v}}^2
+
|(I-P_\gamma)G|_{L^2(\gamma_+)}^2
\notag\\
&+
\|\nu^{-\frac12}\mathfrak F\|_{L^2_{\mathbf x,v}}^2
+
(\alpha+\delta)^2\|\mathcal J[G]\|_{L^2(\partial\Omega)}^2
+
|r|_{L^2(\gamma_-)}^2.
\label{est-L2-Gd-abc}
\end{align}

By \eqref{est-Gd-wall-B}, the contraction
condition holds uniformly in $\epsilon$. Applying the boundary flux
estimate \eqref{est-wall-flux-L2} in Lemma~\ref{lem-wall-flux}  with $\mathcal S=\sigma KG+\mathfrak F$, we obtain
\begin{align}
\|\mathcal J[G]\|_{L^2(\partial\Omega)}
&\lesssim
\|\nu^{-\frac12}KG\|_{L^2_{\mathbf x,v}}
+
\|\nu^{-\frac12}\mathfrak F\|_{L^2_{\mathbf x,v}}
+
|r|_{L^2(\gamma_-)}
\notag\\
&\lesssim
\|\nu^{\frac12}G\|_{L^2_{\mathbf x,v}}
+
\|\nu^{-\frac12}\mathfrak F\|_{L^2_{\mathbf x,v}}
+
|r|_{L^2(\gamma_-)}.
\label{est-Gd-J-L2}
\end{align}

By \eqref{def-P01}, we have
\[
\|\nu^{\frac12}\mathbf P_0G\|_{L^2_{\mathbf x,v}}^2
\sim
\|[a,\mathbf b,c]\|_{L^2_{\mathbf x}}^2,
\qquad
\|\nu^{\frac12}\mathbf P_1G\|_{L^2_{\mathbf x,v}}
\lesssim
\|\nu^{\frac12}G\|_{L^2_{\mathbf x,v}}.
\]
Since $\epsilon+1-\sigma+\sigma\lambda_0\ge\min\{1,\lambda_0\}$,
combining \eqref{est-L2-Gd-energy} and \eqref{est-L2-Gd-abc}
and then using \eqref{est-Gd-J-L2}, we obtain
\[
\begin{aligned}
&
\|\nu^{\frac12}G\|_{L^2_{\mathbf x,v}}^2
+
|(I-P_\gamma)G|_{L^2(\gamma_+)}^2
\\
&\lesssim
\varepsilon\|\nu^{\frac12}G\|_{L^2_{\mathbf x,v}}^2
+
C_\varepsilon\|\nu^{-\frac12}\mathfrak F\|_{L^2_{\mathbf x,v}}^2
+
(\alpha+\delta)^2\|\mathcal J[G]\|_{L^2(\partial\Omega)}^2
+
|r|_{L^2(\gamma_-)}^2
\\
&\lesssim
\big[\varepsilon+(\alpha+\delta)^2\big]
\|\nu^{\frac12}G\|_{L^2_{\mathbf x,v}}^2
+
C_\varepsilon\Big[
\|\nu^{-\frac12}\mathfrak F\|_{L^2_{\mathbf x,v}}^2
+
|r|_{L^2(\gamma_-)}^2
\Big].
\end{aligned}
\]
Taking $\varepsilon>0$ and then $\alpha,\delta>0$ sufficiently small, using
$\nu^{\frac12}w_l^{-1}\in L^2_v$ and
$\int_{\gamma_-}w_l^{-2}(v)\,d\gamma<\infty$ for $l>1$, we obtain
\[
\|w_lG\|_{L^\infty_{\mathbf x,v}}
\lesssim
\|\nu^{-1}w_l\mathfrak F\|_{L^\infty_{\mathbf x,v}}
+
|w_lr|_{L^\infty(\gamma_-)}.
\]
This completes the proof of Lemma~\ref{lem-Gd-apriori}.
\end{proof}

\subsection{Existence and uniqueness}

Having established the uniform a priori estimate in
Lemma~\ref{lem-Gd-apriori}, we now prove the existence and uniqueness of
solutions to the approximate problem \eqref{eq-Gd-eps-sigma} with $\sigma=1$. We then pass
to the limit as $\epsilon\to0$ and construct a solution to
\eqref{eq-Gd-linear} satisfying the zero-mass condition \eqref{eq-Gd-linear-mass}.

For $l>1$, we introduce the auxiliary space
\[
\mathcal X_{\nu,l}
:=
\Bigg\{
\mathsf g=\mathsf g(\mathbf x,v)\ \Bigg|\
\|w_l\mathsf g\|_{L^\infty_{\mathbf x,v}}<\infty,\quad
\int_\Omega\int_{\mathbb R^3}
\nu(v)\sqrt\mu\mathsf g(\mathbf x,v)\,dv d\mathbf x=0
\Bigg\}.
\]

\begin{lemma}\label{lem-Gd-sigma-cont}
Let $0<\epsilon\le1$ and $l>1$, and let $\alpha,\delta>0$
be sufficiently small. Assume that $\mathfrak F$ and $r$
satisfy the mass compatibility condition \eqref{eq-Gd-source-mass}.
Then the boundary value problem \eqref{eq-Gd-eps-sigma} with $\sigma=1$
admits a unique solution $G=G(\mathbf x,v)\in\mathcal X_{\nu,l}$. Moreover,
\begin{equation}\label{est-Gd-sigma-unif}
\|w_lG\|_{L^\infty_{\mathbf x,v}}
\le
C_l
\Big[
\|\nu^{-1}w_l\mathfrak F\|_{L^\infty_{\mathbf x,v}}
+
|w_lr|_{L^\infty(\gamma_-)}
\Big],
\end{equation}
where $C_l>0$ is independent of $\epsilon$.
\end{lemma}

\begin{proof}
\noindent{\textbf{Step 1. The case $\sigma=0$.}}
Fix $0<\epsilon\le1$. We first consider the boundary value problem
obtained by setting $\sigma=0$ in \eqref{eq-Gd-eps-sigma}:
\begin{equation}\label{eq-Gd-zero}
\left\{
\begin{aligned}
&\epsilon\nu(v)G^0
+v\cdot\nabla_{\mathbf x}G^0
+\nu(v)G^0
=
\mathfrak F,
\\
&G^0(x,\pm1,z,v)\big|_{\pm v_y<0}
=
\Phi(x,\pm1,z,v)\mathcal J[G^0](x,\pm1,z)
+r(x,\pm1,z,v).
\end{aligned}
\right.
\end{equation}
The solution $G^0$ constructed below will serve as the
initial iterate in the Picard scheme of Step~2. The characteristic
representation of \eqref{eq-Gd-zero} is
\[
\begin{aligned}
G^0(\mathbf x,v)
={}&
e^{-\frac{(1+\epsilon)\nu(v)(1\pm y)}{|v_y|}}
\Phi\bigg(
x-\frac{(1\pm y)v_x}{|v_y|},
\mp1,
z-\frac{(1\pm y)v_z}{|v_y|},
v
\bigg)
\\
&\qquad\qquad\times
\mathcal J[G^0]
\bigg(
x-\frac{(1\pm y)v_x}{|v_y|},
\mp1,
z-\frac{(1\pm y)v_z}{|v_y|}
\bigg)
\\
&+
e^{-\frac{(1+\epsilon)\nu(v)(1\pm y)}{|v_y|}}
r\bigg(
x-\frac{(1\pm y)v_x}{|v_y|},
\mp1,
z-\frac{(1\pm y)v_z}{|v_y|},
v
\bigg)
\\
&+
\int_0^{\frac{1\pm y}{|v_y|}}
e^{-(1+\epsilon)\nu(v)s}
\mathfrak F(x-v_xs,y-v_ys,z-v_zs,v)\,ds,
\quad \pm v_y>0.
\end{aligned}
\]

Using \eqref{eq-wall-flux-form} in Lemma~\ref{lem-wall-flux}, we obtain
\begin{equation}\label{eq-Gd-zero-flux}
\left\{
\begin{aligned}
\mathcal J[G^0](x,1,z)
={}&
\Big[I-\mathscr B_\lambda(-1)\mathscr B_\lambda(1)\Big]^{-1}
\Big[\mathcal R_{\epsilon,+}^{0}(x,z)
+\mathscr B_\lambda(-1)\mathcal R_{\epsilon,-}^{0}(x,z)\Big],
\\
\mathcal J[G^0](x,-1,z)
={}&
\Big[I-\mathscr B_\lambda(1)\mathscr B_\lambda(-1)\Big]^{-1}
\Big[\mathcal R_{\epsilon,-}^{0}(x,z)
+\mathscr B_\lambda(1)\mathcal R_{\epsilon,+}^{0}(x,z)\Big],
\end{aligned}
\right.
\end{equation}
where
\begin{align}
\mathcal R_{\epsilon,\pm}^{0}(x,z)
:={}&
\int_{\pm u_y>0}
e^{-\frac{2(1+\epsilon)\nu(u)}{|u_y|}}
r\bigg(x-\frac{2u_x}{|u_y|},\mp1,z-\frac{2u_z}{|u_y|},u\bigg)
\sqrt\mu|u_y|\,du
\notag\\
&+
\int_{\pm u_y>0}\int_0^{\frac{2}{|u_y|}}
e^{-(1+\epsilon)\nu(u)s}
\mathfrak F(x-u_xs,\pm1-u_ys,z-u_zs,u)
\sqrt\mu|u_y|\,dsdu.
\label{def-Gd-R-zero}
\end{align}
By \eqref{est-Gd-wall-B}, the inverse operators in
\eqref{eq-Gd-zero-flux} are uniformly bounded in $\epsilon$ and $\sigma$.
Hence, for given $\mathfrak F$ and $r$, the system
\eqref{eq-Gd-zero-flux} admits a unique solution. Substituting these
into the characteristic representation above yields the unique
solution to \eqref{eq-Gd-zero}.

Moreover, \eqref{eq-Gd-mass} gives
$\displaystyle\int_\Omega\int_{\mathbb R^3}\nu(v)\sqrt\mu G^0\,dv d\mathbf x=0$.
Applying Lemma~\ref{lem-Gd-apriori} with $\sigma=0$, we obtain
\begin{equation}\label{est-Gd-zero}
\|w_lG^0\|_{L^\infty_{\mathbf x,v}}
\le
C_l
\Big[
\|\nu^{-1}w_l\mathfrak F\|_{L^\infty_{\mathbf x,v}}
+
|w_lr|_{L^\infty(\gamma_-)}
\Big],
\end{equation}
where $C_l>0$ is independent of $\epsilon$. This implies
that $G^0\in\mathcal X_{\nu,l}$.

\noindent{\textbf{Step 2. The case $0\le\sigma\le\sigma_*$.}}
Let $0<\sigma_*\le1$ be chosen below, and let
$G^0\in\mathcal X_{\nu,l}$ be the unique solution constructed in
Step~1. For a fixed $0\le\sigma\le\sigma_*$, we define the Picard
sequence associated with \eqref{eq-Gd-eps-sigma} recursively by
\begin{equation}\label{eq-Gd-Picard}
\left\{
\begin{aligned}
&\epsilon\nu(v)G^{k+1}
+v\cdot\nabla_{\mathbf x}G^{k+1}
+\nu(v)G^{k+1}
=
\underbrace{\sigma KG^k+\mathfrak F}_{\mathfrak F_1^{(k)}},
\\
&G^{k+1}(x,\pm1,z,v)\big|_{\pm v_y<0}
=
\Phi(x,\pm1,z,v)\mathcal J[G^{k+1}](x,\pm1,z)
+r(x,\pm1,z,v),
\end{aligned}
\right.
\quad k\ge0.
\end{equation}
Suppose that $G^k\in\mathcal X_{\nu,l}$. By
\eqref{eq-Gd-source-mass} and the mass condition of
$G^k$, the modified source $\mathfrak F_1^{(k)}$ satisfies
\[
\begin{aligned}
\int_\Omega\int_{\mathbb R^3}
\sqrt\mu\mathfrak F_1^{(k)}\,dv d\mathbf x
&=
\sigma\int_\Omega\int_{\mathbb R^3}
\sqrt\mu KG^k\,dv d\mathbf x
+
\int_\Omega\int_{\mathbb R^3}
\sqrt\mu\mathfrak F\,dv d\mathbf x
\\
&=
\sigma\int_\Omega\int_{\mathbb R^3}
\nu(v)\sqrt\mu G^k\,dv d\mathbf x
\\
&=0.
\end{aligned}
\]
Therefore, the construction in Step~1 can be applied successively to
uniquely determine $G^{k+1}\in\mathcal X_{\nu,l}$. Hence,
$\{G^k\}_{k\ge0}$ is well defined in $\mathcal X_{\nu,l}$.
For $k\ge1$, subtracting two consecutive equations in
\eqref{eq-Gd-Picard} gives
\[
\left\{
\begin{aligned}
&\epsilon\nu(v)
\big(G^{k+1}-G^k\big)
+
v\cdot\nabla_{\mathbf x}
\big(G^{k+1}-G^k\big)
+
\nu(v)
\big(G^{k+1}-G^k\big)
\\
&\quad=
\sigma K
\big(G^k-G^{k-1}\big),
\\
&\big(G^{k+1}-G^k\big)
(x,\pm1,z,v)\big|_{\pm v_y<0}
=
\Phi(x,\pm1,z,v)
\mathcal J[G^{k+1}-G^k](x,\pm1,z).
\end{aligned}
\right.
\]
Applying Lemma~\ref{lem-Gd-apriori} to this difference problem with
$\sigma=0$ and homogeneous boundary data, and using
\eqref{est-kernel-weight}, we obtain
\begin{align}
\big\|w_l(G^{k+1}-G^k)\big\|_{L^\infty_{\mathbf x,v}}
&\le
C_l\sigma
\big\|\nu^{-1}w_lK(G^k-G^{k-1})\big\|_{L^\infty_{\mathbf x,v}}
\notag\\
&\le
C_l\sigma
\big\|w_l(G^k-G^{k-1})\big\|_{L^\infty_{\mathbf x,v}},
\label{est-Gd-Picard-contract}
\end{align}
where $C_l>0$ is independent of $\epsilon$ and $\sigma$.
Choose $\sigma_*=\sigma_*(l)\in(0,1]$ such that
$C_l\sigma_*\le1/2$. Then, for every
$0\le\sigma\le\sigma_*$, \eqref{est-Gd-Picard-contract} gives
\[
\big\|w_l(G^{k+1}-G^k)\big\|_{L^\infty_{\mathbf x,v}}
\le
\frac12
\big\|w_l(G^k-G^{k-1})\big\|_{L^\infty_{\mathbf x,v}}.
\]
Hence, for every $0\le\sigma\le\sigma_*$, the Picard sequence
$\{G^k\}_{k\ge0}$ converges in $\mathcal X_{\nu,l}$ to the unique
solution $G$ of \eqref{eq-Gd-eps-sigma}. Finally, Lemma~\ref{lem-Gd-apriori}
gives
\[
\|w_lG\|_{L^\infty_{\mathbf x,v}}
\le
C_l
\Big[
\|\nu^{-1}w_l\mathfrak F\|_{L^\infty_{\mathbf x,v}}
+
|w_lr|_{L^\infty(\gamma_-)}
\Big],
\]
where $C_l>0$ is independent of $\epsilon$ and $\sigma$.

\noindent{\textbf{Step 3. Continuation to $\sigma=1$.}}
Reducing $\sigma_*$ if necessary, we may assume that
$2\sigma_*\le1$. Fix $\sigma_*\le\sigma\le2\sigma_*$.
Let $\widetilde G^0\in\mathcal X_{\nu,l}$ be the solution
constructed in Step~2 at $\sigma=\sigma_*$, and define
$\widetilde G^{k+1}$ for $k\ge0$ by
\begin{equation}\label{eq-Gd-Picard-cont}
\left\{
\begin{aligned}
&\epsilon\nu(v)\widetilde G^{k+1}
+v\cdot\nabla_{\mathbf x}\widetilde G^{k+1}
+\nu(v)\widetilde G^{k+1}
-\sigma_*K\widetilde G^{k+1}
=
\underbrace{(\sigma-\sigma_*)K\widetilde G^k
+\mathfrak F}_{\mathfrak F_2^{(k)}},
\\
&\widetilde G^{k+1}(x,\pm1,z,v)\big|_{\pm v_y<0}
=
\Phi(x,\pm1,z,v)\mathcal J[\widetilde G^{k+1}](x,\pm1,z)
+r(x,\pm1,z,v).
\end{aligned}
\right.
\end{equation}
Since $\mathfrak F_2^{(k)}$ satisfies the mass condition
\eqref{eq-Gd-source-mass}, Step~2 yields a unique solution
$\widetilde G^{k+1}\in\mathcal X_{\nu,l}$ for every $k\ge0$.
Thus, the Picard sequence
$\{\widetilde G^k\}_{k\ge0}$ is well defined.
Applying Lemma~\ref{lem-Gd-apriori} to the difference of two
consecutive equations in \eqref{eq-Gd-Picard-cont}, we obtain
\[
\begin{aligned}
\big\|w_l(\widetilde G^{k+1}-\widetilde G^k)\big\|_{L^\infty_{\mathbf x,v}}
&\le
C_l(\sigma-\sigma_*)
\big\|w_l(\widetilde G^k-\widetilde G^{k-1})\big\|_{L^\infty_{\mathbf x,v}}
\\
&\le
\frac12
\big\|w_l(\widetilde G^k-\widetilde G^{k-1})\big\|_{L^\infty_{\mathbf x,v}},
\end{aligned}
\]
where
\[
0\le\sigma-\sigma_*\le\sigma_*,
\quad
C_l\sigma_*\le\frac12.
\]
Hence, arguing as in Step~2, the Picard sequence converges
in $\mathcal X_{\nu,l}$ to the unique solution of
\eqref{eq-Gd-eps-sigma} for every $\sigma_*\le\sigma\le2\sigma_*$.

For $j\ge0$, set
\[
\sigma_j
:=
\min\{j\sigma_*,1\}.
\]
Since the constants above are independent of $\epsilon$ and
$\sigma$, the same argument with $\sigma_j$ as the base parameter
extends the existence and uniqueness of solutions from
$0\le\sigma\le\sigma_j$ to
$0\le\sigma\le\sigma_{j+1}$. After finitely many steps,
$\sigma_j=1$. Therefore, \eqref{eq-Gd-eps-sigma} admits a unique
solution in $\mathcal X_{\nu,l}$ for every $0\le\sigma\le1$. In
particular, at $\sigma=1$, Lemma~\ref{lem-Gd-apriori} yields
\eqref{est-Gd-sigma-unif}. This completes the proof of
Lemma~\ref{lem-Gd-sigma-cont}.
\end{proof}

Let $G_\epsilon$ be the unique solution to \eqref{eq-Gd-eps-sigma}
with $\sigma=1$ given by Lemma~\ref{lem-Gd-sigma-cont}. We first pass to the limit as $\epsilon\to0$ to obtain a unique solution to the following problem in $\mathcal X_{\nu,l}$:
\begin{equation}\label{eq-Gd-linear-nu}
\left\{
\begin{aligned}
&v\cdot\nabla_{\mathbf x}G+\nu(v)G=KG+\mathfrak F,
\\
&G(x,\pm1,z,v)\big|_{\pm v_y<0}
=
\Phi(x,\pm1,z,v)\mathcal J[G](x,\pm1,z)
+r(x,\pm1,z,v),
\\
&\int_\Omega\int_{\mathbb R^3}
\nu(v)\sqrt\mu G\,dv d\mathbf x=0.
\end{aligned}
\right.
\end{equation}
We then recover the zero-mass condition \eqref{eq-Gd-linear-mass}
to obtain the solution to \eqref{eq-Gd-linear} in $\mathcal X_l$.

\begin{proof}[\itshape\bfseries Proof of Lemma~\ref{lem-Gd-linear}]
Let $\epsilon_k:=k^{-1}$ for $k\ge1$.
For each $k\ge1$, let $G_{\epsilon_k}\in\mathcal X_{\nu,l}$
be the unique solution provided by Lemma~\ref{lem-Gd-sigma-cont}
with $\epsilon=\epsilon_k$ and $\sigma=1$. Thus,
\[
\left\{
\begin{aligned}
&\epsilon_k\nu(v)G_{\epsilon_k}
+v\cdot\nabla_{\mathbf x}G_{\epsilon_k}
+\nu(v)G_{\epsilon_k}
=
KG_{\epsilon_k}+\mathfrak F,
\\
&G_{\epsilon_k}(x,\pm1,z,v)\big|_{\pm v_y<0}
=
\Phi(x,\pm1,z,v)\mathcal J[G_{\epsilon_k}](x,\pm1,z)
+r(x,\pm1,z,v).
\end{aligned}
\right.
\]
Moreover, Lemma~\ref{lem-Gd-sigma-cont} gives
\[
\|w_lG_{\epsilon_k}\|_{L^\infty_{\mathbf x,v}}
\le
C_l\Big[
\|\nu^{-1}w_l\mathfrak F\|_{L^\infty_{\mathbf x,v}}
+
|w_lr|_{L^\infty(\gamma_-)}
\Big],
\]
where $C_l>0$ is independent of $\epsilon_k$.

For $k,k'\ge1$, subtracting the corresponding equations yields
\[
\left\{
\begin{aligned}
&\epsilon_k\nu(v)(G_{\epsilon_k}-G_{\epsilon_{k'}})
+v\cdot\nabla_{\mathbf x}(G_{\epsilon_k}-G_{\epsilon_{k'}})
+\nu(v)(G_{\epsilon_k}-G_{\epsilon_{k'}})
\\
&\quad=
K(G_{\epsilon_k}-G_{\epsilon_{k'}})
+(\epsilon_{k'}-\epsilon_k)\nu(v)G_{\epsilon_{k'}},
\\
&(G_{\epsilon_k}-G_{\epsilon_{k'}})(x,\pm1,z,v)\big|_{\pm v_y<0}
=
\Phi(x,\pm1,z,v)
\mathcal J[G_{\epsilon_k}-G_{\epsilon_{k'}}](x,\pm1,z).
\end{aligned}
\right.
\]
Applying Lemma~\ref{lem-Gd-sigma-cont}, we obtain
\begin{align}
\|w_l(G_{\epsilon_k}-G_{\epsilon_{k'}})\|_{L^\infty_{\mathbf x,v}}
&\le
C_l|\epsilon_k-\epsilon_{k'}|
\|w_lG_{\epsilon_{k'}}\|_{L^\infty_{\mathbf x,v}}
\le
\frac{C_l}{\min\{k,k'\}}
\Big[
\|\nu^{-1}w_l\mathfrak F\|_{L^\infty_{\mathbf x,v}}
+
|w_lr|_{L^\infty(\gamma_-)}
\Big].
\label{est-Linf-Gd-eps-Cauchy}
\end{align}
Hence, for some $G_\nu\in\mathcal X_{\nu,l}$, we have
$G_{\epsilon_k}\to G_\nu$ in $\mathcal X_{\nu,l}$.
Passing to the limit in the characteristic representation of the
approximate problem, we conclude that $G_\nu$ solves \eqref{eq-Gd-linear-nu}.
Moreover, the uniform estimate \eqref{est-Gd-sigma-unif}
passes to the limit and yields
\[
\|w_lG_\nu\|_{L^\infty_{\mathbf x,v}}
\le
C_l\Big[
\|\nu^{-1}w_l\mathfrak F\|_{L^\infty_{\mathbf x,v}}
+
|w_lr|_{L^\infty(\gamma_-)}
\Big].
\]

For uniqueness in $\mathcal X_{\nu,l}$, let $G_\nu$ and $\widetilde G_\nu$
be two solutions to \eqref{eq-Gd-linear-nu}, and set
$\mathfrak G:=G_\nu-\widetilde G_\nu$.
Choose $0<\epsilon_0\le1$ such that $C_l\epsilon_0\le1/2$.
Then $\mathfrak G$ satisfies
\begin{equation}\label{eq-Gd-uniq}
\left\{
\begin{aligned}
&\epsilon_0\nu(v)\mathfrak G
+v\cdot\nabla_{\mathbf x}\mathfrak G
+\nu(v)\mathfrak G
=
K\mathfrak G+\epsilon_0\nu(v)\mathfrak G,
\\
&\mathfrak G(x,\pm1,z,v)\big|_{\pm v_y<0}
=
\Phi(x,\pm1,z,v)\mathcal J[\mathfrak G](x,\pm1,z).
\end{aligned}
\right.
\end{equation}
Applying Lemma~\ref{lem-Gd-sigma-cont} to \eqref{eq-Gd-uniq}, we obtain
\[
\|w_l\mathfrak G\|_{L^\infty_{\mathbf x,v}}
\le
C_l\epsilon_0\|w_l\mathfrak G\|_{L^\infty_{\mathbf x,v}}
\le
\frac12\|w_l\mathfrak G\|_{L^\infty_{\mathbf x,v}}.
\]
Thus $G_\nu=\widetilde G_\nu$.

Finally, we recover the mass condition \eqref{eq-Gd-linear-mass}. We claim that there exist a constant $c_0$ and a function $G_h$ such that $G:=G_\nu-c_0G_h$ is the unique solution to \eqref{eq-Gd-linear}--\eqref{eq-Gd-linear-mass}. Indeed, by linearity, $G$ still satisfies the equation and boundary condition in \eqref{eq-Gd-linear-nu} if $G_h$ solves
\begin{equation}\label{eq-Gd-hom}
\left\{
\begin{aligned}
&v\cdot\nabla_{\mathbf x}G_h+\nu(v)G_h=KG_h,
\\
&G_h(x,\pm1,z,v)\big|_{\pm v_y<0}
=
\Phi(x,\pm1,z,v)\mathcal J[G_h](x,\pm1,z).
\end{aligned}
\right.
\end{equation}
If $G_h$ has nonzero mass, let
\begin{equation}\label{def-Gd-c}
c_0=
\frac{\displaystyle\int_\Omega\int_{\mathbb R^3}
\sqrt\mu G_\nu\,dv d\mathbf x}
{\displaystyle\int_\Omega\int_{\mathbb R^3}
\sqrt\mu G_h\,dv d\mathbf x}.
\end{equation}
It therefore remains to find a solution to \eqref{eq-Gd-hom}
with nonzero mass.

Note that $\sqrt\mu$ solves \eqref{eq-Gd-hom} when
$\Phi=\sqrt{2\pi}\sqrt\mu$.
For the perturbed boundary profile, we write $G_h=\sqrt\mu+G_*$
and solve \eqref{eq-Gd-linear-nu} for $G_*$ with
\[
\mathfrak F=0,\qquad
r(x,\pm1,z,v)
=
\frac1{\sqrt{2\pi}}
\big[\Phi(x,\pm1,z,v)-\sqrt{2\pi}\sqrt\mu\big].
\]
By \eqref{eq-Gd-Phi-norm}, these data satisfy
\eqref{eq-Gd-source-mass}.
The preceding existence result and \eqref{eq-Gd-Phi-small}
give a solution $G_*$ satisfying
\[
\|w_lG_*\|_{L^\infty_{\mathbf x,v}}
\le C_l(\alpha+\delta).
\]
Consequently, $G_h=\sqrt\mu+G_*$ solves \eqref{eq-Gd-hom},
satisfies $\|w_lG_h\|_{L^\infty_{\mathbf x,v}}\le C_l$, and,
for $\alpha,\delta$ sufficiently small,
\begin{equation}\label{est-Gd-hom-mass}
\int_\Omega\int_{\mathbb R^3}\sqrt\mu G_h\,dv d\mathbf x
=
|\Omega|+O(\alpha+\delta)
\ge\frac{|\Omega|}{2}.
\end{equation}
For uniqueness, let $\bar G\in\mathcal X_l$ be the difference
of two solutions to \eqref{eq-Gd-linear}. Since
$G_*\in\mathcal X_{\nu,l}$,
\[
\int_\Omega\int_{\mathbb R^3}
\nu(v)\sqrt\mu G_h\,dv d\mathbf x
=
|\Omega|\int_{\mathbb R^3}\nu(v)\mu\,dv>0.
\]
Taking
\[
d_0=
\frac{\displaystyle\int_\Omega\int_{\mathbb R^3}
\nu(v)\sqrt\mu\bar G\,dv d\mathbf x}
{\displaystyle\int_\Omega\int_{\mathbb R^3}
\nu(v)\sqrt\mu G_h\,dv d\mathbf x}.
\]
Then $\bar G-d_0G_h$ belongs to $\mathcal X_{\nu,l}$ and solves
\eqref{eq-Gd-linear-nu} with $\mathfrak F=r=0$, so uniqueness
for that problem gives $\bar G=d_0G_h$. The mass condition \eqref{eq-Gd-linear-mass} and \eqref{est-Gd-hom-mass} then imply $d_0=0$. This completes the proof of Lemma~\ref{lem-Gd-linear}.
\end{proof}

\section{Construction of steady solutions}\label{sec-5}

In this section, we first construct the first-order shear and thermal
corrections $G_1^\alpha$ and $G_1^\delta$ and establish their uniqueness and uniform weighted $L^\infty$ estimates. We then construct the higher-order remainder $G_R^{\alpha,\delta}$.

\begin{proposition}\label{prop-Ga}
Let $l>1$, and let $\alpha,\delta>0$ be sufficiently small.
The boundary value problem \eqref{eq-Ga} admits a unique
solution $G_1^\alpha=G_1^\alpha(y,v)\in\mathcal X_l$, which is
independent of $x$ and $z$ and satisfies
\begin{equation}\label{eq-Ga-sym}
\left\{
\begin{aligned}
G_1^\alpha(y,-v_x,v_y,v_z)
&=-G_1^\alpha(y,v_x,v_y,v_z),
\\
G_1^\alpha(y,v_x,v_y,-v_z)
&=G_1^\alpha(y,v_x,v_y,v_z).
\end{aligned}
\right.
\end{equation}
Moreover,
\begin{equation}\label{est-Linf-Ga}
\big\|w_lG_1^\alpha\big\|_{L^\infty_{y,v}}
\le C_l,
\end{equation}
where $C_l>0$ is independent of $\alpha$ and $\delta$.
\end{proposition}

\begin{proof}
Taking $\mathfrak F=0$, $\Phi=\sqrt{2\pi}\sqrt\mu$, and
$r=\eta_1^\alpha$ in \eqref{eq-Gd-linear}, we obtain \eqref{eq-Ga}.
By \eqref{def-eta1},
\[
\int_{\pm v_y<0}\eta_1^\alpha(\pm1,v)\sqrt\mu|v_y|\,dv=0,
\qquad
|w_l\eta_1^\alpha|_{L^\infty(\gamma_-)}\le C_l.
\]
Lemma~\ref{lem-Gd-linear} therefore gives a unique solution in
$\mathcal X_l$ satisfying \eqref{est-Linf-Ga}.

The problem is invariant under translations in $x$ and $z$, so
uniqueness implies that $G_1^\alpha$ is independent of these variables.
Since $-G_1^\alpha(y,-v_x,v_y,v_z)$ and
$G_1^\alpha(y,v_x,v_y,-v_z)$ also solve \eqref{eq-Ga},
uniqueness yields \eqref{eq-Ga-sym}.
This completes the proof of Proposition~\ref{prop-Ga}.
\end{proof}

Arguing as in the proof of Proposition~\ref{prop-Ga}, we obtain the
following result for the first-order thermal correction $G_1^\delta$.

\begin{proposition}\label{prop-Gd}
Let $l>1$, and let $\alpha,\delta>0$ be sufficiently small. The boundary value problem \eqref{eq-Gd} admits a unique solution $G_1^\delta=G_1^\delta(\mathbf x,v)\in\mathcal X_l$ satisfying
\begin{equation}\label{est-Linf-Gd}
\big\|w_lG_1^\delta\big\|_{L^\infty_{\mathbf x,v}}
\le C_l\|\vartheta\|_{L^\infty(\partial\Omega)},
\end{equation}
where $C_l>0$ is independent of $\alpha$ and $\delta$.
\end{proposition}

With the first-order corrections $G_1^\delta$ and $G_1^\alpha$
constructed in Propositions~\ref{prop-Gd} and \ref{prop-Ga}, we now
construct the higher-order remainder $G_R^{\alpha,\delta}$ by solving
\eqref{eq-GR}.

\begin{proposition}\label{prop-GR}
Let $l>1$, and let $\alpha,\delta>0$ be sufficiently small.
The boundary value problem \eqref{eq-GR} admits a unique solution
$G_R^{\alpha,\delta}=G_R^{\alpha,\delta}(\mathbf x,v)\in\mathcal X_l$
satisfying
\begin{equation}\label{est-Linf-GR}
\big\|w_lG_R^{\alpha,\delta}\big\|_{L^\infty_{\mathbf x,v}}
\le C_l,
\end{equation}
where $C_l>0$ is independent of $\alpha$ and $\delta$.
\end{proposition}

\begin{proof}
For brevity, set
\begin{align}
\mathfrak S(G_1^\alpha,G_1^\delta)
:={}&
\frac{\alpha^2}{\alpha^2+\delta^2}\Gamma(G_1^\alpha,G_1^\alpha)
+
\frac{\alpha\delta}{\alpha^2+\delta^2}
\Big[\Gamma(G_1^\alpha,G_1^\delta)+\Gamma(G_1^\delta,G_1^\alpha)\Big]
+
\frac{\delta^2}{\alpha^2+\delta^2}\Gamma(G_1^\delta,G_1^\delta),
\label{def-GR-S}
\end{align}
and
\begin{align}
\mathcal F_b(x,\pm1,z,v)
:={}&
\eta_R^{\alpha,\delta}(x,\pm1,z,v)
+
\frac{\sqrt{2\pi}\delta}{\alpha^2+\delta^2}
\Big[\alpha\eta_1^\alpha+\delta\eta_1^\delta
+(\alpha^2+\delta^2)\eta_R^{\alpha,\delta}\Big]
\mathcal J[G_1^\delta](x,\pm1,z).
\label{def-GR-r}
\end{align}
Set $G_R^0=0$. For $k\ge0$, having constructed $G_R^k$, let
\begin{align}
\mathcal N_{G_R}^k
:={}&
\alpha\Big[\Gamma(G_1^\alpha,G_R^k)+\Gamma(G_R^k,G_1^\alpha)\Big]
+
\delta\Big[\Gamma(G_1^\delta,G_R^k)+\Gamma(G_R^k,G_1^\delta)\Big]
+
(\alpha^2+\delta^2)\Gamma(G_R^k,G_R^k).
\label{def-GR-N}
\end{align}
Define $G_R^{k+1}\in\mathcal X_l$ by
\begin{equation}\label{eq-GR-iter}
\left\{
\begin{aligned}
v\cdot\nabla_{\mathbf x}G_R^{k+1}+\nu(v)G_R^{k+1}
&=KG_R^{k+1}+\mathfrak S(G_1^\alpha,G_1^\delta)+\mathcal N_{G_R}^k,
\\
G_R^{k+1}(x,\pm1,z,v)\big|_{\pm v_y<0}
&=
\mu^{-\frac12}\mathcal M_{T(x,\pm1,z)}(v_x-U(\pm1),v_y,v_z)
\mathcal J[G_R^{k+1}](x,\pm1,z)
\\
&\quad+\mathcal F_b(x,\pm1,z,v).
\end{aligned}
\right.
\end{equation}
Here we have used $\mathcal J[G_1^\alpha]=0$, which follows from
\eqref{eq-Ga-sym}.

Mass conservation of the collision operator gives
\begin{equation}\label{eq-GR-source-mass}
\int_\Omega\int_{\mathbb R^3}
\sqrt\mu\Big[\mathfrak S(G_1^\alpha,G_1^\delta)
+\mathcal N_{G_R}^k\Big]\,dv d\mathbf x=0.
\end{equation}
Moreover, by \eqref{eq-MT-exp}--\eqref{def-eta1}, we have
\begin{equation}\label{eq-GR-boundary-mass}
\int_{\pm v_y<0}
\eta_R^{\alpha,\delta}\sqrt\mu|v_y|\,dv=0,
\qquad
\int_{\pm v_y<0}\mathcal F_b\sqrt\mu|v_y|\,dv=0.
\end{equation}
Hence \eqref{eq-GR-source-mass}--\eqref{eq-GR-boundary-mass}
give the compatibility conditions \eqref{eq-Gd-source-mass}
for \eqref{eq-GR-iter}.
Lemma~\ref{lem-Gd-linear} therefore determines
$G_R^{k+1}\in\mathcal X_l$ uniquely, so the sequence
$\{G_R^k\}_{k\ge0}$ is well defined.

Using \eqref{est-Linf-Gamma} in Lemma~\ref{lem-Gamma},
together with \eqref{est-Linf-Ga} and \eqref{est-Linf-Gd}, we obtain
\begin{align}
\big\|\nu^{-1}w_l\mathfrak S(G_1^\alpha,G_1^\delta)\big\|_{L^\infty_{\mathbf x,v}}
\le{}&
\frac{C_l\alpha^2}{\alpha^2+\delta^2}
\|w_lG_1^\alpha\|_{L^\infty_{\mathbf x,v}}^2
+
\frac{C_l\alpha\delta}{\alpha^2+\delta^2}
\|w_lG_1^\alpha\|_{L^\infty_{\mathbf x,v}}
\|w_lG_1^\delta\|_{L^\infty_{\mathbf x,v}}
\notag\\
&+
\frac{C_l\alpha\delta}{\alpha^2+\delta^2}
\|w_lG_1^\delta\|_{L^\infty_{\mathbf x,v}}
\|w_lG_1^\alpha\|_{L^\infty_{\mathbf x,v}}
+
\frac{C_l\delta^2}{\alpha^2+\delta^2}
\|w_lG_1^\delta\|_{L^\infty_{\mathbf x,v}}^2
\notag\\
\le{}&
C_l^3\frac{\alpha^2+2\alpha\delta+\delta^2}{\alpha^2+\delta^2}
\le 2C_l^3
\le C_l.
\label{est-GR-S}
\end{align}
Similarly,
\begin{align}
\big\|\nu^{-1}w_l\mathcal N_{G_R}^k\big\|_{L^\infty_{\mathbf x,v}}
\le{}&
C_l\alpha
\|w_lG_1^\alpha\|_{L^\infty_{\mathbf x,v}}
\|w_lG_R^k\|_{L^\infty_{\mathbf x,v}}
+
C_l\alpha
\|w_lG_R^k\|_{L^\infty_{\mathbf x,v}}
\|w_lG_1^\alpha\|_{L^\infty_{\mathbf x,v}}
\notag\\
&+
C_l\delta
\|w_lG_1^\delta\|_{L^\infty_{\mathbf x,v}}
\|w_lG_R^k\|_{L^\infty_{\mathbf x,v}}
+
C_l\delta
\|w_lG_R^k\|_{L^\infty_{\mathbf x,v}}
\|w_lG_1^\delta\|_{L^\infty_{\mathbf x,v}}
+
C_l(\alpha^2+\delta^2)
\|w_lG_R^k\|_{L^\infty_{\mathbf x,v}}^2
\notag\\
\le{}&
C_l(\alpha+\delta)\|w_lG_R^k\|_{L^\infty_{\mathbf x,v}}
+
C_l(\alpha^2+\delta^2)\|w_lG_R^k\|_{L^\infty_{\mathbf x,v}}^2.
\label{est-Linf-GR-N}
\end{align}
By \eqref{eq-MT-exp}, \eqref{def-eta1}, and \eqref{est-Linf-Gd},
\begin{align}
|w_l\mathcal F_b|_{L^\infty(\gamma_-)}
\lesssim{}&
|w_l\eta_R^{\alpha,\delta}|_{L^\infty(\gamma_-)}
+
\Big[
|w_l\eta_1^\alpha|_{L^\infty(\gamma_-)}
+
|w_l\eta_1^\delta|_{L^\infty(\gamma_-)}
+
|w_l\eta_R^{\alpha,\delta}|_{L^\infty(\gamma_-)}
\Big]
\|\mathcal J[G_1^\delta]\|_{L^\infty(\partial\Omega)}
\notag\\
\le{}&
C_l+C_l^2
\le C_l.
\label{est-GR-r}
\end{align}
Combining \eqref{est-GR-S}--\eqref{est-GR-r}, and \eqref{est-Linf-Gd-linear}, we obtain
\begin{align}
\|w_lG_R^{k+1}\|_{L^\infty_{\mathbf x,v}}
\le{}&
C_l+C_l(\alpha+\delta)\|w_lG_R^k\|_{L^\infty_{\mathbf x,v}}
\notag\\
&+
C_l(\alpha^2+\delta^2)\|w_lG_R^k\|_{L^\infty_{\mathbf x,v}}^2.
\label{est-GR-eps-iter}
\end{align}
Choose $R_0\ge C_l$ independent of $\alpha$, $\delta$, and $k$.
Since $G_R^0=0$, taking $\alpha,\delta>0$ sufficiently small
and applying induction gives
\begin{equation}\label{est-GR-eps-uniform}
\|w_lG_R^k\|_{L^\infty_{\mathbf x,v}}
\le2R_0,\qquad k\ge0.
\end{equation}

For $k\ge0$, by \eqref{eq-GR-iter}, we obtain
\begin{equation}\label{eq-GR-iter-diff}
\left\{
\begin{aligned}
&v\cdot\nabla_{\mathbf x}(G_R^{k+2}-G_R^{k+1})
+\nu(v)(G_R^{k+2}-G_R^{k+1})
=K(G_R^{k+2}-G_R^{k+1})
+(\mathcal N_{G_R}^{k+1}-\mathcal N_{G_R}^k),
\\
&(G_R^{k+2}-G_R^{k+1})(x,\pm1,z,v)\big|_{\pm v_y<0}
=
\mu^{-\frac12}\mathcal M_{T(x,\pm1,z)}(v_x-U(\pm1),v_y,v_z)
\mathcal J[G_R^{k+2}-G_R^{k+1}](x,\pm1,z),
\\
&\int_\Omega\int_{\mathbb R^3}
\sqrt\mu(G_R^{k+2}-G_R^{k+1})\,dv d\mathbf x
=0.
\end{aligned}
\right.
\end{equation}
By \eqref{def-GR-N}, we have
\begin{align}
\mathcal N_{G_R}^{k+1}-\mathcal N_{G_R}^k
={}&
\alpha\Big[
\Gamma(G_1^\alpha,G_R^{k+1}-G_R^k)
+
\Gamma(G_R^{k+1}-G_R^k,G_1^\alpha)
\Big]
\notag\\
&+
\delta\Big[
\Gamma(G_1^\delta,G_R^{k+1}-G_R^k)
+
\Gamma(G_R^{k+1}-G_R^k,G_1^\delta)
\Big]
\notag\\
&+
(\alpha^2+\delta^2)\Big[
\Gamma(G_R^{k+1},G_R^{k+1}-G_R^k)
+
\Gamma(G_R^{k+1}-G_R^k,G_R^k)
\Big].
\label{eq-GR-N-diff}
\end{align}
Arguing as in \eqref{est-GR-S}--\eqref{est-Linf-GR-N}
and using \eqref{est-GR-eps-uniform}, we obtain
\begin{equation}\label{est-GR-N-diff}
\big\|\nu^{-1}w_l
(\mathcal N_{G_R}^{k+1}-\mathcal N_{G_R}^k)
\big\|_{L^\infty_{\mathbf x,v}}
\le
C_l\Big[\alpha+\delta+R_0(\alpha^2+\delta^2)\Big]
\|w_l(G_R^{k+1}-G_R^k)\|_{L^\infty_{\mathbf x,v}}.
\end{equation}
Applying Lemma~\ref{lem-Gd-linear} to \eqref{eq-GR-iter-diff}
and taking $\alpha,\delta>0$ sufficiently small, we have
\begin{align}
\|w_l(G_R^{k+2}-G_R^{k+1})\|_{L^\infty_{\mathbf x,v}}
&\le
C_l\Big[\alpha+\delta+R_0(\alpha^2+\delta^2)\Big]
\|w_l(G_R^{k+1}-G_R^k)\|_{L^\infty_{\mathbf x,v}}
\notag\\
&\le
\frac12\|w_l(G_R^{k+1}-G_R^k)\|_{L^\infty_{\mathbf x,v}}.
\notag
\end{align}
Hence $\{G_R^k\}_{k\ge0}$ is Cauchy in $\mathcal X_l$
and converges to the unique solution $G_R^{\alpha,\delta}$ of \eqref{eq-GR}
satisfying \eqref{est-Linf-GR}. This completes the proof of Proposition~\ref{prop-GR}.
\end{proof}

\section{Local existence for the unsteady problem}\label{sec-6}

We first establish the local well-posedness of the unsteady problem \eqref{eq-prob-3d} in a weighted $L^\infty$ framework. For $T>0$, we introduce the Banach space
\[
\mathcal Y_T
:=
\Big\{
[\mathsf g_1,\mathsf g_2]:
\big\|[\mathsf g_1,\mathsf g_2]\big\|_{\mathcal Y_T}<\infty
\Big\},
\]
equipped with the norm
\begin{equation}\label{def-YT}
\big\|[\mathsf g_1,\mathsf g_2]\big\|_{\mathcal Y_T}
:=
\sup_{0\le t\le T}
\Big[
\|w_l\mathsf g_1(t)\|_{L^\infty_{\mathbf x,v}}
+
\|w_l\mathsf g_2(t)\|_{L^\infty_{\mathbf x,v}}
\Big].
\end{equation}

To study \eqref{eq-prob-3d} around the steady state $F_{\alpha,\delta}^s$ constructed in Theorem~\ref{thm-1}, we write
\[
F(t,\mathbf x,v)=F_{\alpha,\delta}^s(\mathbf x,v)+\sqrt\mu g(t,\mathbf x,v),
\]
where
\[
F_{\alpha,\delta}^s=\mu+\sqrt\mu G^s,
\quad
G^s=\alpha G_1^\alpha+\delta G_1^\delta+(\alpha^2+\delta^2)G_R^{\alpha,\delta}.
\]
By Theorem~\ref{thm-1},
\begin{equation}\label{est-Linf-Gs}
\big\|w_l\sqrt\mu G^s\big\|_{L^\infty_{\mathbf x,v}}
\le C_l(\alpha+\delta).
\end{equation}
Subtracting the steady equation for $F_{\alpha,\delta}^s$ from \eqref{eq-prob-3d} and using \eqref{eq-MT-exp}, we obtain
\begin{equation}\label{eq-g-unsteady}
\left\{
\begin{aligned}
&\partial_tg+v\cdot\nabla_{\mathbf x}g+Lg
=
\Gamma(G^s,g)+\Gamma(g,G^s)+\Gamma(g,g),
\quad
(t,\mathbf x,v)\in\mathbb R_+\times\Omega\times\mathbb R^3,
\\
&g(t,x,\pm1,z,v)\big|_{\pm v_y<0}
=
P_\gamma g(t,x,\pm1,z,v)
\\
&\quad+
\sqrt{2\pi}\Big[
\alpha\eta_1^\alpha(\pm1,v)
+\delta\eta_1^\delta(x,\pm1,z,v)
+(\alpha^2+\delta^2)\eta_R^{\alpha,\delta}(x,\pm1,z,v)
\Big]\mathcal J[g](t,x,\pm1,z),
\\
&g(0,\mathbf x,v)
=
\mu^{-\frac12}\Big[F_0(\mathbf x,v)-F_{\alpha,\delta}^s(\mathbf x,v)\Big],
\quad
\int_\Omega\int_{\mathbb R^3}\sqrt\mu g(t,\mathbf x,v)\,dv d\mathbf x=0.
\end{aligned}
\right.
\end{equation}

Let $M>0$ be chosen sufficiently large, and let $\chi_M=\chi_M(v)$ be a smooth cutoff function with $0\le\chi_M\le1$ satisfying
\begin{equation}\label{def-ChiM}
\chi_M(v)
=
\begin{cases}
0,& |v|\le M,\\
1,& |v|\ge M+1.
\end{cases}
\end{equation}
We further use Caflisch's decomposition $\sqrt\mu g=g_1+\sqrt\mu g_2$, and take $g_1$ and $g_2$ to satisfy the following coupled systems:
\begin{equation}\label{eq-g1-unsteady}
\left\{
\begin{aligned}
&\partial_tg_1+v\cdot\nabla_{\mathbf x}g_1+\nu(v)g_1
=
\chi_M\mathcal Kg_1+H(g_1,g_2),
\\
&g_1(t,x,\pm1,z,v)\big|_{\pm v_y<0}=0,
\quad
g_1(0,\mathbf x,v)=F_0(\mathbf x,v)-F_{\alpha,\delta}^s(\mathbf x,v),
\end{aligned}
\right.
\end{equation}
and
\begin{equation}\label{eq-g2-unsteady}
\left\{
\begin{aligned}
&\partial_tg_2+v\cdot\nabla_{\mathbf x}g_2+\nu(v)g_2
=
Kg_2+(1-\chi_M)\mu^{-\frac12}\mathcal Kg_1,
\\
&g_2(t,x,\pm1,z,v)\big|_{\pm v_y<0}
=
P_\gamma g(t,x,\pm1,z,v)
\\
&\quad+
\sqrt{2\pi}\Big[
\alpha\eta_1^\alpha(\pm1,v)
+\delta\eta_1^\delta(x,\pm1,z,v)
+(\alpha^2+\delta^2)\eta_R^{\alpha,\delta}(x,\pm1,z,v)
\Big]\mathcal J[g](t,x,\pm1,z),
\\
&g_2(0,\mathbf x,v)=0.
\end{aligned}
\right.
\end{equation}
Here
\begin{align}
H(g_1,g_2)
={}&
Q(\sqrt\mu G^s,g_1)+Q(g_1,\sqrt\mu G^s)
+Q(\sqrt\mu G^s,\sqrt\mu g_2)+Q(\sqrt\mu g_2,\sqrt\mu G^s)
\notag\\
&+
Q(g_1,g_1)+Q(g_1,\sqrt\mu g_2)
+Q(\sqrt\mu g_2,g_1)+Q(\sqrt\mu g_2,\sqrt\mu g_2).
\label{def-H}
\end{align}

\begin{lemma}[Local existence]\label{lem-local-existence}
For $-3<\gamma\le1$, suppose that $l>0$ is sufficiently large. There exist constants $A=A(l)>1$ and $\varepsilon_0=\varepsilon_0(l)>0$ such that, if $\alpha>0$ and $\delta>0$ are sufficiently small, then for any initial datum $F_0=F_0(\mathbf x,v)$ satisfying
\[
\big\|w_l(F_0-F_{\alpha,\delta}^s)\big\|_{L^\infty_{\mathbf x,v}}
\le \varepsilon_0,
\]
there exists $T>0$ such that the coupled system \eqref{eq-g1-unsteady}--\eqref{eq-g2-unsteady} admits a unique solution $[g_1,g_2]\in\mathcal Y_T$ satisfying
\[
\big\|[g_1,g_2]\big\|_{\mathcal Y_T}
\le
2A\big\|w_l(F_0-F_{\alpha,\delta}^s)\big\|_{L^\infty_{\mathbf x,v}}.
\]
\end{lemma}

\begin{proof}
Set
\[
h_1(t,\mathbf x,v):=w_l(v)g_1(t,\mathbf x,v),
\quad
h_2(t,\mathbf x,v):=w_l(v)g_2(t,\mathbf x,v).
\]
Then $h_1$ satisfies
\begin{equation}\label{eq-local-h1}
\left\{
\begin{aligned}
&\partial_th_1+v\cdot\nabla_{\mathbf x}h_1+\nu(v)h_1
=
\chi_Mw_l\mathcal K(w_l^{-1}h_1)+w_lH(g_1,g_2),
\\
&h_1(t,x,\pm1,z,v)\big|_{\pm v_y<0}=0,
\\
&h_1(0,\mathbf x,v)
=
w_l\big[F_0(\mathbf x,v)-F_{\alpha,\delta}^s(\mathbf x,v)\big],
\end{aligned}
\right.
\end{equation}
while $h_2$ satisfies
\begin{equation}\label{eq-local-h2}
\left\{
\begin{aligned}
&\partial_th_2+v\cdot\nabla_{\mathbf x}h_2+\nu(v)h_2
=
w_lKg_2+(1-\chi_M)w_l\mu^{-\frac12}\mathcal K(w_l^{-1}h_1),
\\
&h_2(t,x,\pm1,z,v)\big|_{\pm v_y<0}
=
\sqrt{2\pi}w_l\Big[
\sqrt\mu+\alpha\eta_1^\alpha(\pm1,v)+\delta\eta_1^\delta(x,\pm1,z,v)
\\
&\qquad\qquad\qquad\qquad\qquad\quad
+(\alpha^2+\delta^2)\eta_R^{\alpha,\delta}(x,\pm1,z,v)
\Big]\mathcal J[g](t,x,\pm1,z),
\\
&h_2(0,\mathbf x,v)=0.
\end{aligned}
\right.
\end{equation}

Along characteristics, we write
\begin{equation}\label{eq-local-Duhamel}
h_1(t,\mathbf x,v)
=
\sum_{i=0}^2\mathcal I_i^{(1)}(t,\mathbf x,v),
\quad
h_2(t,\mathbf x,v)
=
\sum_{i=1}^3\mathcal I_i^{(2)}(t,\mathbf x,v),
\end{equation}
where
\begin{equation}\label{def-local-I1}
\left\{
\begin{aligned}
\mathcal I_0^{(1)}
={}&
\mathbf 1_{\{t_b\ge t\}}e^{-\nu(v)t}h_1\big(0,\mathbf X(0),v\big),
\\
\mathcal I_1^{(1)}
={}&
\int_{\max\{0,t-t_b\}}^t
e^{-\nu(v)(t-s)}
\Big[\chi_Mw_l\mathcal K(w_l^{-1}h_1)\Big]\big(s,\mathbf X(s),v\big)\,ds,
\\
\mathcal I_2^{(1)}
={}&
\int_{\max\{0,t-t_b\}}^t
e^{-\nu(v)(t-s)}
\Big[w_lH(g_1,g_2)\Big]\big(s,\mathbf X(s),v\big)\,ds,
\end{aligned}
\right.
\end{equation}
and
\begin{equation}\label{def-local-I2}
\left\{
\begin{aligned}
\mathcal I_1^{(2)}
={}&
\mathbf 1_{\{t_b<t\}}e^{-\nu(v)t_b}\sqrt{2\pi}w_l(v)
\Big[
\sqrt\mu+\alpha\eta_1^\alpha\big(Y(t-t_b),v\big)
+\delta\eta_1^\delta\big(\mathbf X(t-t_b),v\big)
\\
&\quad+
(\alpha^2+\delta^2)\eta_R^{\alpha,\delta}\big(\mathbf X(t-t_b),v\big)
\Big]\mathcal J[g]\big(t-t_b,\mathbf X(t-t_b)\big),
\\
\mathcal I_2^{(2)}
={}&
\int_{\max\{0,t-t_b\}}^t
e^{-\nu(v)(t-s)}
\big[w_lKg_2\big]\big(s,\mathbf X(s),v\big)\,ds,
\\
\mathcal I_3^{(2)}
={}&
\int_{\max\{0,t-t_b\}}^t
e^{-\nu(v)(t-s)}
\Big[(1-\chi_M)w_l\mu^{-\frac12}\mathcal K(w_l^{-1}h_1)\Big]
\big(s,\mathbf X(s),v\big)\,ds.
\end{aligned}
\right.
\end{equation}
Here $\mathbf X(s):=\mathbf x-(t-s)v$, and $t_b=t_b(\mathbf x,v)$ is given by \eqref{def-tby}.

To derive a more explicit representation of $\mathcal I_1^{(2)}$, we recall from \eqref{def-wall-B-time} that
\[
\begin{aligned}
&\big(\mathscr B_\nu^t(\pm1)f\big)(s,x,\mp1,z)
\\&=
\sqrt{2\pi}\int_{\pm u_y<0}
\mathbf 1_{\{t_b<s\}}e^{-\nu(u)t_b}
\Big[\sqrt\mu+\alpha\eta_1^\alpha(\pm1,u)
+
\delta\eta_1^\delta(x-u_xt_b,\pm1,z-u_zt_b,u)
\\
&\quad+(\alpha^2+\delta^2)\eta_R^{\alpha,\delta}(x-u_xt_b,\pm1,z-u_zt_b,u)
\Big]\sqrt\mu|u_y|
f\big(s-t_b,x-u_xt_b,\pm1,z-u_zt_b\big)\,du.
\end{aligned}
\]
Arguing as in \eqref{est-Gd-wall-B}, we obtain, for $\alpha$ and $\delta$ sufficiently small,
\begin{align}
&\max_\pm
\|\mathscr B_\nu^t(\pm1)\|_{L^\infty(0,t;L^\infty(\partial\Omega))\to L^\infty(0,t;L^\infty(\partial\Omega))}
\notag\\
&\le
\sqrt{2\pi}\int_{u_y>0}
e^{-\frac{2\nu(u)}{u_y}}\mu|u_y|\,du
+C_l(\alpha+\delta)
<1,
\qquad \forall 0<t\le T.
\label{est-local-wall-B}
\end{align}

Consequently, applying the argument leading to \eqref{eq-wall-flux-form-time}
in Lemma~\ref{lem-wall-flux-time} to \eqref{eq-g-unsteady}, we obtain
\begin{equation}\label{eq-g-wall-flux-time}
\left\{
\begin{aligned}
\mathcal J[g](t,x,1,z)
={}&
\bigg\{
\Big[I-\mathscr B_\nu^t(-1)\mathscr B_\nu^t(1)\Big]^{-1}
\Big[\mathcal R_{\nu,+}+\mathscr B_\nu^t(-1)\mathcal R_{\nu,-}\Big]
\bigg\}(t,x,1,z),
\\
\mathcal J[g](t,x,-1,z)
={}&
\bigg\{
\Big[I-\mathscr B_\nu^t(1)\mathscr B_\nu^t(-1)\Big]^{-1}
\Big[\mathcal R_{\nu,-}+\mathscr B_\nu^t(1)\mathcal R_{\nu,+}\Big]
\bigg\}(t,x,-1,z),
\end{aligned}
\right.
\end{equation}
where
\begin{align}
\mathcal R_{\nu,\pm}(t,x,z)
:={}&
\int_{\pm u_y>0}\int_{\max\{0,t-t_b\}}^t
e^{-\nu(u)(t-\tau)}
\Big[\mathcal Kg_1+\sqrt\mu Kg_2+H(g_1,g_2)\Big]
\big(\tau,\mathbf X(\tau),u\big)
|u_y|\,d\tau du
\notag\\
&+
\int_{\pm u_y>0}
\mathbf 1_{\{t_b\ge t\}}e^{-\nu(u)t}
\Big[F_0-F_{\alpha,\delta}^s\Big]
\big(\mathbf X(0),u\big)|u_y|\,du.
\label{def-local-wall-R}
\end{align}

Set
\[
\mathbb B_T
:=
\Big\{
[g_1,g_2]\in\mathcal Y_T:
\big\|[g_1,g_2]\big\|_{\mathcal Y_T}
\le
2A\big\|w_l(F_0-F_{\alpha,\delta}^s)\big\|_{L^\infty_{\mathbf x,v}}
\Big\}.
\]
For any $[g_1,g_2]\in\mathbb B_T$, define the nonlinear mapping
$\mathscr T$ by the Duhamel formulas
\eqref{def-local-I1}--\eqref{def-local-I2}, with the boundary fluxes
given by \eqref{eq-g-wall-flux-time}, where $\mathcal R_{\nu,\pm}$
are evaluated at the given input $[g_1,g_2]$. Then the coupled system
is equivalent to the fixed-point equation
\[
[g_1,g_2]=\mathscr T[g_1,g_2].
\]

\noindent\textbf{Step 1. Self-mapping.}
Let $l>\dfrac12\max\{3,3+\gamma\}$. Since
\[
\int_{\max\{0,t-t_b\}}^t e^{-\nu(v)(t-s)}\,ds
\le
\int_0^t e^{-\nu(v)(t-s)}\,ds
\le
\min\{T,\nu^{-1}(v)\},
\]
it follows from \eqref{est-Linf-Q} in Lemma~\ref{lem-Gamma},
Lemma~\ref{lem-calK-high}, and \eqref{est-Linf-Gs} that
\begin{align}
\sum_{i=0}^2\big|\mathcal I_i^{(1)}\big|
\le{}&
\big\|w_l(F_0-F_{\alpha,\delta}^s)\big\|_{L^\infty_{\mathbf x,v}}
+
\sup_{0\le s\le T}
\big\|\nu^{-1}\chi_Mw_l\mathcal K g_1(s)\big\|_{L^\infty_{\mathbf x,v}}
+
\sup_{0\le s\le T}
\big\|\nu^{-1}w_lH(g_1,g_2)(s)\big\|_{L^\infty_{\mathbf x,v}}
\notag\\
\le{}&
\big\|w_l(F_0-F_{\alpha,\delta}^s)\big\|_{L^\infty_{\mathbf x,v}}
+
\bigg(\frac{C}{l^{\frac{\gamma+3}{4}}}+\frac{C_l}{M^2}\bigg)
\sup_{0\le s\le T}\big\|w_lg_1(s)\big\|_{L^\infty_{\mathbf x,v}}
\notag\\
&+
C_l\sup_{0\le s\le T}
\bigg\{
\Big[
\big\|w_l\sqrt\mu G^s\big\|_{L^\infty_{\mathbf x,v}}
+
\big\|w_l(g_1+\sqrt\mu g_2)(s)\big\|_{L^\infty_{\mathbf x,v}}
\Big]
\big\|w_l(g_1+\sqrt\mu g_2)(s)\big\|_{L^\infty_{\mathbf x,v}}
\bigg\}
\notag\\
\le{}&
\big\|w_l(F_0-F_{\alpha,\delta}^s)\big\|_{L^\infty_{\mathbf x,v}}
+
\bigg(\frac{C}{l^{\frac{\gamma+3}{4}}}+\frac{C_l}{M^2}\bigg)
\big\|[g_1,g_2]\big\|_{\mathcal Y_T}
+
C_l\Big[
(\alpha+\delta)\big\|[g_1,g_2]\big\|_{\mathcal Y_T}
+\big\|[g_1,g_2]\big\|_{\mathcal Y_T}^2
\Big].
\label{est-local-I1}
\end{align}

For $\mathcal I_2^{(2)}$ and $\mathcal I_3^{(2)}$, by
\eqref{est-kernel-weight} and \eqref{est-Linf-Q} we obtain
\begin{align}
\big|\mathcal I_2^{(2)}\big|
+
\big|\mathcal I_3^{(2)}\big|
\le{}&
C_lT\frac{\nu(v)}{(1+|v|)^2}
\sup_{0\le s\le T}\big\|w_lg_2(s)\big\|_{L^\infty_{\mathbf x,v}}
+
C_{l,M}T
\sup_{0\le s\le T}\big\|w_lg_1(s)\big\|_{L^\infty_{\mathbf x,v}}
\notag\\
\le{}&
C_{l,M}T\big\|[g_1,g_2]\big\|_{\mathcal Y_T}.
\label{est-local-I23-2}
\end{align}

For $\mathcal I_1^{(2)}$, since $(1+\nu(u))w_l^{-1}(u)|u_y|\in L^1(\mathbb R_v^3)$ for $l$ sufficiently large, it follows from \eqref{est-local-wall-B}-- \eqref{def-local-wall-R}, and the estimates leading to \eqref{est-local-I1} and \eqref{est-local-I23-2} that
\begin{align}
\sum_\pm\big|\mathcal J[g](t,x,\pm1,z)\big|
\le{}&
C\big\|w_l(F_0-F_{\alpha,\delta}^s)\big\|_{L^\infty_{\mathbf x,v}}
+
C_lT\Big[
\big\|[g_1,g_2]\big\|_{\mathcal Y_T}
+\big\|[g_1,g_2]\big\|_{\mathcal Y_T}^2
\Big].
\label{est-local-flux}
\end{align}
Substituting \eqref{est-local-flux} into the expression for
$\mathcal I_1^{(2)}$ in \eqref{def-local-I2}, we obtain
\begin{equation}\label{est-local-I1-2}
\big|\mathcal I_1^{(2)}(t,\mathbf x,v)\big|
\le
C_l\big\|w_l(F_0-F_{\alpha,\delta}^s)\big\|_{L^\infty_{\mathbf x,v}}
+
C_lT\Big[
\big\|[g_1,g_2]\big\|_{\mathcal Y_T}
+\big\|[g_1,g_2]\big\|_{\mathcal Y_T}^2
\Big].
\end{equation}

We first take $l$ sufficiently large and then fix $M=M(l)$ so that
\[
\frac{C}{l^{\frac{\gamma+3}{4}}}+\frac{C_l}{M^2}\le\frac18.
\]
Choose $A=A(l)>1$ sufficiently large such that $C_l\le A/2$, and then take $\alpha,\delta>0$ sufficiently small so that $C_l(\alpha+\delta)\le1/8$.
Further, take $\varepsilon_0>0$ and $T>0$ sufficiently small.
Then, for any $[g_1,g_2]\in\mathbb B_T$, combining
\eqref{est-local-I1}, \eqref{est-local-I23-2}, and
\eqref{est-local-I1-2}, we obtain
\begin{align}
\big\|\mathscr T[g_1,g_2]\big\|_{\mathcal Y_T}
&\le
C_l\big\|w_l(F_0-F_{\alpha,\delta}^s)\big\|_{L^\infty_{\mathbf x,v}}
+\frac12\big\|[g_1,g_2]\big\|_{\mathcal Y_T}
+C_l\big\|[g_1,g_2]\big\|_{\mathcal Y_T}^2
\notag\\
&\le
2A\big\|w_l(F_0-F_{\alpha,\delta}^s)\big\|_{L^\infty_{\mathbf x,v}}.
\label{est-local-map}
\end{align}
Hence, $\mathscr T(\mathbb B_T)\subset\mathbb B_T$.

\noindent\textbf{Step 2. Contraction.}
For $[g_1,g_2],[g_1',g_2']\in\mathbb B_T$, set
\[
\bar g_1=g_1-g_1',
\quad
\bar g_2=g_2-g_2'.
\]
By the same estimates as above, and Lemma~\ref{lem-wall-flux-time}
applied to the difference equation, we obtain
\begin{align}
&\big\|\mathscr T[g_1,g_2]-\mathscr T[g_1',g_2']\big\|_{\mathcal Y_T}
\notag\\
&\le
\bigg[
\frac{C}{l^{\frac{\gamma+3}{4}}}+\frac{C_l}{M^2}
+C_l(\alpha+\delta)+C_{l,M}T
+C_l\Big(
\big\|[g_1,g_2]\big\|_{\mathcal Y_T}
+\big\|[g_1',g_2']\big\|_{\mathcal Y_T}
\Big)
\bigg]
\big\|[\bar g_1,\bar g_2]\big\|_{\mathcal Y_T}
\notag\\
&\le
\bigg[
\frac{C}{l^{\frac{\gamma+3}{4}}}+\frac{C_l}{M^2}
+C_l(\alpha+\delta)+C_{l,M}T+4C_lA\varepsilon_0
\bigg]
\big\|[\bar g_1,\bar g_2]\big\|_{\mathcal Y_T}.
\label{est-local-contraction}
\end{align}
After fixing $l$ and $M$ and taking $\alpha,\delta>0$ sufficiently small
as above, we take $\varepsilon_0>0$ sufficiently small and further
decrease $T$ if necessary, so that
\[
\big\|\mathscr T[g_1,g_2]-\mathscr T[g_1',g_2']\big\|_{\mathcal Y_T}
\le
\frac12\big\|[\bar g_1,\bar g_2]\big\|_{\mathcal Y_T}.
\]
Thus $\mathscr T$ is a contraction on $\mathbb B_T$. By the Banach
fixed-point theorem, $\mathscr T$ admits a unique fixed point
$[g_1,g_2]\in\mathbb B_T$. This completes the proof of Lemma~\ref{lem-local-existence}.
\end{proof}

\section{Global stability for the unsteady problem}\label{sec-7}

By Lemma~\ref{lem-local-existence}, the unsteady problem admits a unique local solution. Let $T^*\in(0,\infty]$ denote its maximal existence time. We shall prove that $T^*=\infty$ and that the solution converges to the steady state $F_{\alpha,\delta}^s$. For any $0<T<T^*$, assume that, for some sufficiently small constant $\kappa>0$,
\begin{equation}\label{apriori}
\sup_{0\le t\le T}e^{\lambda_\gamma t}\Big[\|w_lg_1(t)\|_{L^\infty_{\mathbf x,v}}+\|w_lg_2(t)\|_{L^\infty_{\mathbf x,v}}\Big]\le\kappa,
\end{equation}
where
\[
\lambda_\gamma=
\left\{
\begin{aligned}
\lambda, &\quad 0\le\gamma\le1,\qquad 0<\lambda<\frac12\min\big\{\nu_0,\widetilde\lambda_0\big\},\\
0, &\quad -3<\gamma<0,
\end{aligned}
\right.
\]
with $\nu_0>0$ a lower bound for $\nu(v)$ when $0\le\gamma\le1$ and $\widetilde\lambda_0>0$ to be specified in \eqref{est-global-L2}.
We derive estimates uniform in $T$ that improve the a priori assumption \eqref{apriori}. These estimates yield global existence by a standard continuation argument and give exponential decay for $0\le\gamma\le1$. For $-3<\gamma<0$, we then establish algebraic decay with a loss of velocity weight.

Recalling that $h_1=w_lg_1$ and $h_2=w_lg_2$, we multiply the Duhamel representations \eqref{eq-local-Duhamel} by $e^{\lambda_\gamma t}$ and write
\begin{equation}\label{eq-global-Duhamel}
e^{\lambda_\gamma t}h_1(t,\mathbf x,v)=\sum_{i=0}^2\mathcal I_{i,\gamma}^{(1)}(t,\mathbf x,v),
\quad
e^{\lambda_\gamma t}h_2(t,\mathbf x,v)=\sum_{i=1}^3\mathcal I_{i,\gamma}^{(2)}(t,\mathbf x,v),
\end{equation}
where
\[
\mathcal I_{i,\gamma}^{(1)}:=e^{\lambda_\gamma t}\mathcal I_i^{(1)},
\quad
\mathcal I_{i,\gamma}^{(2)}:=e^{\lambda_\gamma t}\mathcal I_i^{(2)}.
\]
The estimates for the terms in the above decompositions are collected in the following lemma.

\begin{lemma}\label{lem-global-Linf}
For $-3<\gamma\le1$, let $[g_1,g_2]$ be the solution to \eqref{eq-g1-unsteady}--\eqref{eq-g2-unsteady} satisfying the a priori assumption \eqref{apriori}. There exists $l_0>0$ such that, for any $l\ge l_0$ and sufficiently small $\alpha,\delta,\kappa>0$, we have
\begin{equation}\label{est-global-h1}
\sup_{0\le s\le t}e^{\lambda_\gamma s}\|h_1(s)\|_{L^\infty_{\mathbf x,v}}
\le
C_l\big\|w_l(F_0-F_{\alpha,\delta}^s)\big\|_{L^\infty_{\mathbf x,v}}
+
C_l(\alpha+\delta+\kappa)
\sup_{0\le s\le t}e^{\lambda_\gamma s}\|h_2(s)\|_{L^\infty_{\mathbf x,v}},
\end{equation}
and
\begin{equation}\label{est-global-h2}
\sup_{0\le s\le t}e^{\lambda_\gamma s}\|h_2(s)\|_{L^\infty_{\mathbf x,v}}
\le
C_l\big\|w_l(F_0-F_{\alpha,\delta}^s)\big\|_{L^\infty_{\mathbf x,v}}
+
C_l\sup_{0\le s\le t}e^{\lambda_\gamma s}\|h_1(s)\|_{L^\infty_{\mathbf x,v}}
+
C_l\sup_{0\le s\le t}e^{\lambda_\gamma s}\|g_2(s)\|_{L^2_{\mathbf x,v}}.
\end{equation}
\end{lemma}

\begin{proof}
We first estimate the first component in \eqref{eq-global-Duhamel}. Since $\nu(v)-\lambda_\gamma\ge\nu(v)/2>0$, Lemma~\ref{lem-calK-high} yields
\begin{align}
\big|\mathcal I_{0,\gamma}^{(1)}\big|+\big|\mathcal I_{1,\gamma}^{(1)}\big|
\le{}&
e^{-(\nu(v)-\lambda_\gamma)t}
\big\|w_l(F_0-F_{\alpha,\delta}^s)\big\|_{L^\infty_{\mathbf x,v}}
\notag\\
&+
\bigg(\int_{\max\{0,t-t_b\}}^t
e^{-(\nu(v)-\lambda_\gamma)(t-s)}\nu(v)\,ds\bigg)
\sup_{0\le s\le t}e^{\lambda_\gamma s}
\big\|\nu^{-1}\chi_Mw_l\mathcal K g_1(s)\big\|_{L^\infty_{\mathbf x,v}}
\notag\\
\le{}&
\big\|w_l(F_0-F_{\alpha,\delta}^s)\big\|_{L^\infty_{\mathbf x,v}}
+
\bigg(\frac{C}{l^{\frac{\gamma+3}{4}}}+\frac{C_l}{M^2}\bigg)
\sup_{0\le s\le t}e^{\lambda_\gamma s}\|h_1(s)\|_{L^\infty_{\mathbf x,v}}.
\label{est-global-I012-1}
\end{align}
Moreover, by \eqref{est-Linf-Q} in Lemma~\ref{lem-Gamma}, \eqref{est-Linf-Gs}, and the a priori assumption \eqref{apriori}, we obtain
\begin{equation}\label{est-global-I3-1}
\big|\mathcal I_{2,\gamma}^{(1)}\big|
\le 
C_l(\alpha+\delta+\kappa)\sup_{0\le s\le t}e^{\lambda_\gamma s}\Big[\|h_1(s)\|_{L^\infty_{\mathbf x,v}}
+\|h_2(s)\|_{L^\infty_{\mathbf x,v}}\Big].
\end{equation}
Combining \eqref{est-global-I012-1} and \eqref{est-global-I3-1}, and then taking $l$ sufficiently large, $M$ sufficiently large depending on $l$, and $\alpha,\delta,\kappa>0$ sufficiently small, we obtain \eqref{est-global-h1}.

For $\mathcal I_{3,\gamma}^{(2)}$, we have
\begin{equation}\label{est-global-I3-2}
\big|\mathcal I_{3,\gamma}^{(2)}\big|
\le
C_l\sup_{0\le s\le t}e^{\lambda_\gamma s}
\|h_1(s)\|_{L^\infty_{\mathbf x,v}}.
\end{equation}

We next estimate the boundary term $\mathcal I_{1,\gamma}^{(2)}$. Since $\nu(v)-\lambda_\gamma\ge\nu(v)/2$, the argument in \eqref{est-local-wall-B} with damping $\nu(v)/2$ gives uniform bounds for the inverses with damping $\nu(v)-\lambda_\gamma$ for sufficiently small $\alpha+\delta$. Applying Lemma~\ref{lem-wall-flux-time} to $e^{\lambda_\gamma t}g$ in \eqref{eq-g-unsteady}, we obtain
\begin{equation}\label{est-global-flux-R}
\sum_\pm\big|\mathcal J[e^{\lambda_\gamma t}g](t,x,\pm1,z)\big|
\lesssim
\sup_{0\le s\le t}\sum_\pm
\|\mathcal R_{\nu-\lambda_\gamma,\pm}(s)\|_{L^\infty_{x,z}},
\end{equation}
where
\begin{align}
\mathcal R_{\nu-\lambda_\gamma,\pm}(t,x,z)
={}&
\int_{\pm u_y>0}\mathbf 1_{\{t_b\ge t\}}e^{-(\nu(u)-\lambda_\gamma)t}
\big[F_0-F_{\alpha,\delta}^s\big]\big(\mathbf X(0),u\big)|u_y|\,du
\notag\\
&+
\int_{\pm u_y>0}\int_{\max\{0,t-t_b\}}^t
e^{-(\nu(u)-\lambda_\gamma)(t-s)}e^{\lambda_\gamma s}
\big[\mathcal Kg_1\big]\big(s,\mathbf X(s),u\big)|u_y|\,dsdu
\notag\\
&+
\int_{\pm u_y>0}\int_{\max\{0,t-t_b\}}^t
e^{-(\nu(u)-\lambda_\gamma)(t-s)}e^{\lambda_\gamma s}
\big[\sqrt\mu Kg_2\big]\big(s,\mathbf X(s),u\big)|u_y|\,dsdu
\notag\\
&+
\int_{\pm u_y>0}\int_{\max\{0,t-t_b\}}^t
e^{-(\nu(u)-\lambda_\gamma)(t-s)}e^{\lambda_\gamma s}
H(g_1,g_2)\big(s,\mathbf X(s),u\big)|u_y|\,dsdu
\notag\\
=:{}&
\mathcal R_\pm^{(0)}
+\mathcal R_\pm^{(1)}
+\mathcal R_\pm^{(2)}
+\mathcal R_\pm^{(3)}.
\label{eq-global-R-decomp}
\end{align}
By the same estimates as in \eqref{est-global-I012-1}--\eqref{est-global-I3-2}, we have
\begin{align}
&\big|\mathcal R_\pm^{(0)}\big|
+\big|\mathcal R_\pm^{(1)}\big|
+\big|\mathcal R_\pm^{(3)}\big|
\notag\\
&\le
C_l\big\|w_l(F_0-F_{\alpha,\delta}^s)\big\|_{L^\infty_{\mathbf x,v}}
+
C_l\sup_{0\le s\le t}e^{\lambda_\gamma s}\|h_1(s)\|_{L^\infty_{\mathbf x,v}}
+
C_l(\alpha+\delta+\kappa)
\sup_{0\le s\le t}e^{\lambda_\gamma s}\|h_2(s)\|_{L^\infty_{\mathbf x,v}}.
\label{est-global-R-easy}
\end{align}

To estimate $\mathcal R_\pm^{(2)}$, we follow the decompositions in \eqref{eq-Gd-RK-decomp}. We first separate the region $0\le t-s<\varepsilon_4$, and on $t-s\ge\varepsilon_4$ further split according to $|u|\le M_4$ and $|u-v_*|<M_4^{-1}$, where $M_4>0$ is sufficiently large. For each fixed $s$ with $t-s\ge\varepsilon_4$, we use the change of variables
\[
x'=x-u_x(t-s),
\quad
y'=\pm1-u_y(t-s),
\quad
z'=z-u_z(t-s),
\]
for which
\[
\bigg|
\det\frac{\partial(x',y',z')}{\partial(u_x,u_y,u_z)}
\bigg|
=
(t-s)^3.
\]
Proceeding as in \eqref{est-Gd-RK-short}--\eqref{est-Gd-RK}, we obtain
\begin{equation}\label{est-global-R2}
\big|\mathcal R_\pm^{(2)}\big|
\le
C_l\varepsilon_4\sup_{0\le s\le t}e^{\lambda_\gamma s}\|h_2(s)\|_{L^\infty_{\mathbf x,v}}
+
\frac{C_l}{M_4^{\frac{\gamma+3}{2}}}\sup_{0\le s\le t}e^{\lambda_\gamma s}\|h_2(s)\|_{L^\infty_{\mathbf x,v}}
+
C_{l,\varepsilon_4,M_4}\sup_{0\le s\le t}e^{\lambda_\gamma s}\|g_2(s)\|_{L^2_{\mathbf x,v}}.
\end{equation}

Combining \eqref{est-global-flux-R}, \eqref{est-global-R-easy}, and \eqref{est-global-R2}, substituting the resulting estimate into $\mathcal I_{1,\gamma}^{(2)}$, and using \eqref{est-global-I3-2}, we obtain
\begin{align}
e^{\lambda_\gamma t}\big|h_2(t,\mathbf x,v)\big|
\le{}&
C_l\big\|w_l(F_0-F_{\alpha,\delta}^s)\big\|_{L^\infty_{\mathbf x,v}}
+
C_l\sup_{0\le s\le t}e^{\lambda_\gamma s}\|h_1(s)\|_{L^\infty_{\mathbf x,v}}
\notag\\
&+
C_l\Big(\varepsilon_4+\frac{1}{M_4^{\frac{\gamma+3}{2}}}+\alpha+\delta+\kappa\Big)
\sup_{0\le s\le t}e^{\lambda_\gamma s}\|h_2(s)\|_{L^\infty_{\mathbf x,v}}
\notag\\
&+
C_{l,\varepsilon_4,M_4}
\sup_{0\le s\le t}e^{\lambda_\gamma s}\|g_2(s)\|_{L^2_{\mathbf x,v}}
+
\big|\mathcal I_{2,\gamma}^{(2)}(t,\mathbf x,v)\big|.
\label{est-global-h2-pre-I2}
\end{align}

For $\mathcal I_{2,\gamma}^{(2)}$, we use the standard velocity decomposition in \eqref{eq-Gd-I2-decomp}, together with a short--long time decomposition in the Duhamel iteration. Substituting \eqref{est-global-h2-pre-I2} into the resulting second-flight term, choosing $l$ and then $M_4$ sufficiently large, and taking $\varepsilon_4,\alpha,\delta,\kappa>0$ sufficiently small, we obtain \eqref{est-global-h2}. This completes the proof of Lemma~\ref{lem-global-Linf}.
\end{proof}

The following lemma provides the $L^2$ estimate needed to close the weighted $L^\infty$ estimate.

\begin{lemma}\label{lem-global-L2}
For $-3<\gamma\le1$, let $[g_1,g_2]$ be the solution to \eqref{eq-g1-unsteady}--\eqref{eq-g2-unsteady} satisfying the a priori assumption \eqref{apriori}, and let $l\ge4$. Then, for any sufficiently small $\varepsilon_5>0$, if $\alpha,\delta,\kappa>0$ are sufficiently small, we have
\begin{align}
\sup_{0\le s\le t}e^{\lambda_\gamma s}\|g_2(s)\|_{L^2_{\mathbf x,v}}
\le{}&
C_l\big\|w_l(F_0-F_{\alpha,\delta}^s)\big\|_{L^\infty_{\mathbf x,v}}
+
C_{\varepsilon_5,l}\sup_{0\le s\le t}e^{\lambda_\gamma s}\|h_1(s)\|_{L^\infty_{\mathbf x,v}}
\notag\\
&+
C_l\varepsilon_5\sup_{0\le s\le t}e^{\lambda_\gamma s}\|h_2(s)\|_{L^\infty_{\mathbf x,v}}.
\label{est-global-g2-L2-weight}
\end{align}
\end{lemma}

\begin{proof}
\noindent\textbf{Step 1. Microscopic estimate.}
We first derive the microscopic energy estimate from the $g_2$ equation. For brevity, set
\begin{equation}\label{def-global-eta}
\eta^{\alpha,\delta}:=
\alpha\eta_1^\alpha+\delta\eta_1^\delta+(\alpha^2+\delta^2)\eta_R^{\alpha,\delta}.
\end{equation}
Recalling \eqref{eq-g2-unsteady}, we write
\begin{equation}\label{eq-global-g2-L2}
\left\{
\begin{aligned}
&\partial_tg_2+v\cdot\nabla_{\mathbf x}g_2+Lg_2
=(1-\chi_M)\mu^{-\frac12}\mathcal Kg_1,
\\
&g_2(t,x,\pm1,z,v)\big|_{\pm v_y<0}
=P_\gamma g(t,x,\pm1,z,v)
+\sqrt{2\pi}\eta^{\alpha,\delta}(x,\pm1,z,v)\mathcal J[g](t,x,\pm1,z),
\\
&g_2(0,\mathbf x,v)=0.
\end{aligned}
\right.
\end{equation}
Taking the $L^2(\Omega\times\mathbb R_v^3)$ inner product of the first equation in \eqref{eq-global-g2-L2} with $g_2$ and using Lemma~\ref{lem-L-coer}, we obtain
\begin{align}
&\frac{d}{dt}\|g_2\|_{L^2_{\mathbf x,v}}^2
+2\lambda_0\|\nu^{\frac12}\mathbf P_1g_2\|_{L^2_{\mathbf x,v}}^2
+\sum_\pm\int_{\mathbb T^2}\int_{\mathbb R^3}
|g_2(t,x,\pm1,z,v)|^2(\pm v_y)\,dv dxdz
\notag\\
&\le 2\big|\langle(1-\chi_M)\mu^{-\frac12}\mathcal Kg_1,g_2\rangle_{\mathbf x,v}\big|.
\label{est-global-g2-energy-pre}
\end{align}
We readily obtain
\begin{align}
2\big|\langle(1-\chi_M)\mu^{-\frac12}\mathcal Kg_1,g_2\rangle_{\mathbf x,v}\big|
&\le
C_l\|h_1\|_{L^\infty_{\mathbf x,v}}
\Big[\|\nu^{\frac12}\mathbf P_0g_2\|_{L^2_{\mathbf x,v}}
+\|\nu^{\frac12}\mathbf P_1g_2\|_{L^2_{\mathbf x,v}}\Big]
\notag\\
&\le
C_l\|h_1\|_{L^\infty_{\mathbf x,v}}
\Big[\|\mathbf P_0g_2\|_{L^2_{\mathbf x,v}}
+\|\nu^{\frac12}\mathbf P_1g_2\|_{L^2_{\mathbf x,v}}\Big]
\notag\\
&\le
\varepsilon\|\mathbf P_0g_2\|_{L^2_{\mathbf x,v}}^2
+\varepsilon\|\nu^{\frac12}\mathbf P_1g_2\|_{L^2_{\mathbf x,v}}^2
+C_{\varepsilon,l}\|h_1\|_{L^\infty_{\mathbf x,v}}^2.
\label{est-global-g2-source}
\end{align}

By \eqref{eq-MT-exp} and \eqref{eq-Gd-Phi-norm}--\eqref{eq-Gd-Phi-small},
\begin{equation}\label{est-global-eta-boundary}
\int_{\pm v_y<0}\eta^{\alpha,\delta}(x,\pm1,z,v)\sqrt\mu|v_y|\,dv=0,
\qquad
|\eta^{\alpha,\delta}\mathcal J[g]|_{L^2(\gamma_-)}
\le C_l(\alpha+\delta)\|\mathcal J[g]\|_{L^2(\partial\Omega)}.
\end{equation}
Moreover, Caflisch's decomposition gives
$\mathcal J[g_2]=\mathcal J[g]-\mathcal J[\mu^{-\frac12}g_1]$,
\[
\begin{aligned}
\big|\mathcal J[\mu^{-\frac12}g_1](t,x,\pm1,z)\big|
&=\bigg|\int_{\pm u_y>0}g_1(t,x,\pm1,z,u)|u_y|\,du\bigg|
\\
&\le \|h_1(t)\|_{L^\infty_{\mathbf x,v}}
\int_{\mathbb R^3}\frac{|u_y|}{w_l(u)}\,du
\le C_l\|h_1(t)\|_{L^\infty_{\mathbf x,v}}.
\end{aligned}
\]
Substituting the incoming boundary condition in \eqref{eq-global-g2-L2} into the incoming integral and using the orthogonality of $P_\gamma g_2$ and $(I-P_\gamma)g_2$ in $L^2(\gamma_+)$ and the zero incoming flux in \eqref{est-global-eta-boundary}, for any $\varepsilon_5>0$, we have
\begin{align}
&\sum_\pm\int_{\mathbb T^2}\int_{\mathbb R^3}
|g_2(t,x,\pm1,z,v)|^2(\pm v_y)\,dv dxdz
\notag\\
&=
\sum_\pm\int_{\mathbb T^2}\bigg[
\int_{\pm v_y>0}|g_2(t,x,\pm1,z,v)|^2|v_y|\,dv
-\int_{\pm v_y<0}|g_2(t,x,\pm1,z,v)|^2|v_y|\,dv
\bigg]\,dxdz
\notag\\
&=
|(I-P_\gamma)g_2|_{L^2(\gamma_+)}^2
+\sqrt{2\pi}\Big[\|\mathcal J[g_2]\|_{L^2(\partial\Omega)}^2
-\|\mathcal J[g]\|_{L^2(\partial\Omega)}^2\Big]
-2\pi|\eta^{\alpha,\delta}\mathcal J[g]|_{L^2(\gamma_-)}^2
\notag\\
&\ge
|(I-P_\gamma)g_2|_{L^2(\gamma_+)}^2
-2\sqrt{2\pi}\|\mathcal J[g]\|_{L^2(\partial\Omega)}
\|\mathcal J[\mu^{-\frac12}g_1]\|_{L^2(\partial\Omega)}
-C_l(\alpha+\delta)^2\|\mathcal J[g]\|_{L^2(\partial\Omega)}^2
\notag\\
&\ge
|(I-P_\gamma)g_2|_{L^2(\gamma_+)}^2
-C_l\big[\varepsilon_5+(\alpha+\delta)^2\big]\|\mathcal J[g]\|_{L^2(\partial\Omega)}^2
-C_{\varepsilon_5,l}\|h_1\|_{L^\infty_{\mathbf x,v}}^2.
\label{est-global-g2-boundary}
\end{align}
Combining \eqref{est-global-g2-energy-pre}, \eqref{est-global-g2-source}, and \eqref{est-global-g2-boundary}, and taking $\varepsilon>0$ sufficiently small, we obtain
\begin{align}
&\frac{d}{dt}\|g_2\|_{L^2_{\mathbf x,v}}^2
+|(I-P_\gamma)g_2|_{L^2(\gamma_+)}^2
+\lambda_0\|\nu^{\frac12}\mathbf P_1g_2\|_{L^2_{\mathbf x,v}}^2
\notag\\
&\le
\varepsilon\|\mathbf P_0g_2\|_{L^2_{\mathbf x,v}}^2
+C_{\varepsilon,\varepsilon_5,l}\|h_1\|_{L^\infty_{\mathbf x,v}}^2
+C_l\big[\varepsilon_5+(\alpha+\delta)^2\big]\|\mathcal J[g]\|_{L^2(\partial\Omega)}^2.
\label{eq-global-g2-energy}
\end{align}

\noindent\textbf{Step 2. Macroscopic estimates.}
To recover the macroscopic dissipation, we return to the reconstructed perturbation $g=\mu^{-\frac12}g_1+g_2$. Recalling \eqref{def-chi} and \eqref{def-P01}, we write
\begin{equation}\label{def-global-macro-g}
\left\{
\begin{aligned}
\mathbf P_0g&=\bigg[a+\mathbf b\cdot v+c\frac{|v|^2-3}{\sqrt6}\bigg]\sqrt\mu,
\qquad \mathbf b=(b_x,b_y,b_z),\qquad \mathbf P_1g=g-\mathbf P_0g
\\
a&=\langle g,\chi_0\rangle_v,\quad
b_x=\langle g,\chi_1\rangle_v,\quad b_y=\langle g,\chi_2\rangle_v,
\quad b_z=\langle g,\chi_3\rangle_v,\quad c=\langle g,\chi_4\rangle_v.
\end{aligned}
\right.
\end{equation}
Since $\mathbf P_0g_2=\mathbf P_0g-\mathbf P_0(\mu^{-\frac12}g_1)$, we have
\begin{align}
\|\mathbf P_0g_2\|_{L^2_{\mathbf x,v}}
+\|\nu^{\frac12}\mathbf P_0g_2\|_{L^2_{\mathbf x,v}}
&\lesssim
\|[a,\mathbf b,c]\|_{L^2_{\mathbf x}}
+
\bigg\|
\int_{\mathbb R^3}(1+|v|^2)|g_1|\,dv
\bigg\|_{L^2_{\mathbf x}}
\notag\\
&\lesssim
\|[a,\mathbf b,c]\|_{L^2_{\mathbf x}}
+
\|h_1\|_{L^\infty_{\mathbf x,v}},
\qquad l\ge3.
\label{est-global-P0-g2}
\end{align}

Taking the $L^2(\mathbb R_v^3)$ inner products of \eqref{eq-g-unsteady} with $\chi_j$, $0\le j\le4$, we obtain
\begin{equation}\label{eq-global-macro-system}
\left\{
\begin{aligned}
&\partial_t a+\nabla_{\mathbf x}\cdot\mathbf b=0,
\\
&\partial_t b_i+\partial_i\Big(a+\frac{2}{\sqrt6}c\Big)
+\nabla_{\mathbf x}\cdot\int_{\mathbb R^3}v_i v\sqrt\mu\mathbf P_1g\,dv=0,
\qquad i\in\{x,y,z\},
\\
&\partial_t c+\frac{2}{\sqrt6}\nabla_{\mathbf x}\cdot\mathbf b
+\frac{1}{\sqrt6}\nabla_{\mathbf x}\cdot\int_{\mathbb R^3}(|v|^2-5)v\sqrt\mu\mathbf P_1g\,dv=0.
\end{aligned}
\right.
\end{equation}
By the mass condition in \eqref{eq-g-unsteady} and the initial condition \eqref{eq-t0-datum}, we have
\begin{equation}\label{eq-global-a-mass}
\int_\Omega a(t,\mathbf x)\,d\mathbf x=0,\qquad \forall t\ge0.
\end{equation}

\textbf{Estimate of $a$.}
By \eqref{eq-global-a-mass}, we introduce $\phi_a=\phi_a(t,\mathbf x)$ and $\psi_a=\psi_a(t,\mathbf x,v)$ by
\begin{equation}\label{def-global-psia}
\left\{
\begin{aligned}
-\Delta_{\mathbf x}\phi_a
&=
a,\quad \mathbf x\in\Omega,
\\
\partial_y\phi_a\big|_{\partial\Omega}
&=
0,\quad
\int_\Omega\phi_a(t,\mathbf x)\,d\mathbf x=0,
\\
\psi_a(t,\mathbf x,v)
&=
(10-|v|^2)\sqrt\mu v\cdot\nabla_{\mathbf x}\phi_a(t,\mathbf x).
\end{aligned}
\right.
\end{equation}
Here and below, the potentials are periodic in $x,z$. Then,
\begin{equation}\label{est-global-phi}
\|\phi_a\|_{H^2_{\mathbf x}}
+\|\partial_x\phi_a\|_{L^2(\partial\Omega)}
+\|\partial_z\phi_a\|_{L^2(\partial\Omega)}
\lesssim
\|a\|_{L^2_{\mathbf x}}.
\end{equation}

Taking the $L^2(\Omega\times\mathbb R_v^3)$ inner product of \eqref{eq-g-unsteady} with $\psi_a$, we obtain
\begin{align}
\frac{d}{dt}\langle g,\psi_a\rangle_{\mathbf x,v}
+
5\|a\|_{L^2_{\mathbf x}}^2
={}&
\langle g,\partial_t\psi_a\rangle_{\mathbf x,v}
+
\langle\mathbf P_1g,v\cdot\nabla_{\mathbf x}\psi_a\rangle_{\mathbf x,v}
\notag\\
&-
\underbrace{
\int_{\gamma_+}g(t,x,\pm1,z,v)\psi_a\,d\gamma
-
\int_{\gamma_-}g(t,x,\pm1,z,v)\psi_a\,d\gamma
}_{\mathcal B_a}
\notag\\
&-
\langle Lg,\psi_a\rangle_{\mathbf x,v}
+
\langle\mu^{-\frac12}H(g_1,g_2),\psi_a\rangle_{\mathbf x,v}.
\label{eq-global-a-test}
\end{align}

We next estimate the terms one by one. By the Neumann boundary condition in \eqref{def-global-psia},
\[
\psi_a(t,x,\pm1,z,v)
=
(10-|v|^2)\sqrt\mu
\Big[
v_x\partial_x\phi_a(t,x,\pm1,z)
+
v_z\partial_z\phi_a(t,x,\pm1,z)
\Big].
\]
Arguing as in \eqref{est-L2-G1d-incom} and using the incoming boundary condition in \eqref{eq-g-unsteady}, we obtain
\begin{align}
\mathcal B_a
={}&
\int_{\gamma_+}\mu^{-\frac12}g_1(t,x,\pm1,z,v)\psi_a\,d\gamma
+
\int_{\gamma_+}(I-P_\gamma)g_2(t,x,\pm1,z,v)\psi_a\,d\gamma
\notag\\
&-
\sqrt{2\pi}\int_{\gamma_-}\eta^{\alpha,\delta}\mathcal J[g]\psi_a\,d\gamma.
\label{eq-global-a-bd-decomp}
\end{align}
By \eqref{est-global-phi},
\[
\begin{aligned}
\bigg|\int_{\gamma_+}\mu^{-\frac12}g_1\psi_a\,d\gamma\bigg|
={}&
\bigg|\int_{\gamma_+}g_1(10-|v|^2)v\cdot\nabla_{\mathbf x}\phi_a\,d\gamma\bigg|
\\
\lesssim{}&
\|h_1\|_{L^\infty_{\mathbf x,v}}
\int_{\mathbb R^3}\frac{1+|v|^4}{w_l(v)}\,dv
\Big[
\|\partial_x\phi_a\|_{L^2(\partial\Omega)}
+
\|\partial_z\phi_a\|_{L^2(\partial\Omega)}
\Big]
\\
\lesssim{}&
\|h_1\|_{L^\infty_{\mathbf x,v}}\|a\|_{L^2_{\mathbf x}}.
\end{aligned}
\]
Then, by \eqref{est-global-phi} agian and \eqref{est-global-eta-boundary}, we obtain
\begin{align}
|\mathcal B_a|
\lesssim{}&
\Big[
\|h_1\|_{L^\infty_{\mathbf x,v}}
+
|(I-P_\gamma)g_2|_{L^2(\gamma_+)}
+
|\eta^{\alpha,\delta}\mathcal J[g]|_{L^2(\gamma_-)}
\Big]\|a\|_{L^2_{\mathbf x}}
\notag\\
\le{}&
\varepsilon\|a\|_{L^2_{\mathbf x}}^2
+
C_{\varepsilon,l}\Big[
\|h_1\|_{L^\infty_{\mathbf x,v}}^2
+
|(I-P_\gamma)g_2|_{L^2(\gamma_+)}^2
+
(\alpha+\delta)^2\|\mathcal J[g]\|_{L^2(\partial\Omega)}^2
\Big].
\label{est-global-a-boundary}
\end{align}

For the term $\langle g,\partial_t\psi_a\rangle_{\mathbf x,v}$ in \eqref{eq-global-a-test}, we apply $\partial_t$ to the first equation in \eqref{def-global-psia} and use the first equation in \eqref{eq-global-macro-system} to obtain
\[
-\Delta_{\mathbf x}\partial_t\phi_a=-\nabla_{\mathbf x}\cdot\mathbf b.
\]
As in \eqref{eq-Gd-by-bc}, the boundary condition in \eqref{eq-g-unsteady} and \eqref{est-global-eta-boundary} give $b_y(t,x,\pm1,z)=0$. Testing with $\partial_t\phi_a$ and using $\partial_y\partial_t\phi_a\big|_{\partial\Omega}=0$, we have
\[
\|\nabla_{\mathbf x}\partial_t\phi_a\|_{L^2_{\mathbf x}}^2
=
\langle\mathbf b,\nabla_{\mathbf x}\partial_t\phi_a\rangle_{\mathbf x}.
\]
It follows that
\begin{equation}\label{est-global-dt-phi}
\|\nabla_{\mathbf x}\partial_t\phi_a\|_{L^2_{\mathbf x}}
\le
\|\mathbf b\|_{L^2_{\mathbf x}}.
\end{equation}
By \eqref{est-global-dt-phi}, we have
\begin{align}
\big|\langle g,\partial_t\psi_a\rangle_{\mathbf x,v}\big|
={}&
\bigg|
\bigg\langle
5\mathbf b-\int_{\mathbb R^3}(|v|^2-5)v\sqrt\mu g\,dv,
\nabla_{\mathbf x}\partial_t\phi_a
\bigg\rangle_{\mathbf x}
\bigg|
\notag\\
={}&
\bigg|
\bigg\langle
5\mathbf b-\int_{\mathbb R^3}(|v|^2-5)v\big[g_1+\sqrt\mu \mathbf P_1g_2\big]\,dv,
\nabla_{\mathbf x}\partial_t\phi_a
\bigg\rangle_{\mathbf x}
\bigg|
\notag\\
\lesssim{}&
\bigg[
\|\mathbf b\|_{L^2_{\mathbf x}}
+
\|\nu^{-\frac12}(|v|^2-5)v\sqrt\mu\|_{L^2_v}
\|\nu^{\frac12}\mathbf P_1g_2\|_{L^2_{\mathbf x,v}}
+
\|h_1\|_{L^\infty_{\mathbf x,v}}
\int_{\mathbb R^3}\frac{1+|v|^3}{w_l(v)}\,dv
\bigg]\|\mathbf b\|_{L^2_{\mathbf x}}
\notag\\
\lesssim{}&
\|\mathbf b\|_{L^2_{\mathbf x}}^2
+
\|\nu^{\frac12}\mathbf P_1g_2\|_{L^2_{\mathbf x,v}}^2
+
\|h_1\|_{L^\infty_{\mathbf x,v}}^2,
\qquad l\ge4.
\label{est-global-a-dt}
\end{align}

The self-adjointness of $\mathbf P_1$ give
\begin{align}
\big|\langle\mathbf P_1g,v\cdot\nabla_{\mathbf x}\psi_a\rangle_{\mathbf x,v}\big|
={}&
\bigg|
\big\langle\mu^{-\frac12}g_1,\mathbf P_1(v\cdot\nabla_{\mathbf x}\psi_a)\big\rangle_{\mathbf x,v}
+
\big\langle\mathbf P_1g_2,v\cdot\nabla_{\mathbf x}\psi_a\big\rangle_{\mathbf x,v}
\bigg|
\notag\\
\lesssim{}&
\Big[
\|h_1\|_{L^\infty_{\mathbf x,v}}
\int_{\mathbb R^3}\frac{1+|v|^4}{w_l(v)}\,dv
+
\|\nu^{\frac12}\mathbf P_1g_2\|_{L^2_{\mathbf x,v}}
\Big]\|\phi_a\|_{H^2_{\mathbf x}}
\notag\\
\le{}&
\varepsilon\|a\|_{L^2_{\mathbf x}}^2
+
C_\varepsilon\Big[
\|h_1\|_{L^\infty_{\mathbf x,v}}^2
+
\|\nu^{\frac12}\mathbf P_1g_2\|_{L^2_{\mathbf x,v}}^2
\Big].
\label{est-global-a-transport}
\end{align}
Arguing as in \eqref{est-L2-Gd-a-rem} and using \eqref{def-calL-calK} and \eqref{est-Linf-Q}, we obtain
\begin{align}
\big|\langle Lg,\psi_a\rangle_{\mathbf x,v}\big|
&=
\big|\langle\mathcal Lg_1,\mu^{-\frac12}\psi_a\rangle_{\mathbf x,v}
+\langle L\mathbf P_1g_2,\psi_a\rangle_{\mathbf x,v}\big|
\notag\\
&\lesssim
\|\nu^{-1}w_l\mathcal Lg_1\|_{L^\infty_{\mathbf x,v}}
\bigg(\int_{\mathbb R^3}\frac{\nu(v)(1+|v|^3)}{w_l(v)}\,dv\bigg)
\|a\|_{L^2_{\mathbf x}}
+
\|\nu^{-\frac12}L\mathbf P_1g_2\|_{L^2_{\mathbf x,v}}
\|\nu^{\frac12}\psi_a\|_{L^2_{\mathbf x,v}}
\notag\\
&\le
\varepsilon\|a\|_{L^2_{\mathbf x}}^2
+
C_{\varepsilon,l}\Big[
\|h_1\|_{L^\infty_{\mathbf x,v}}^2
+
\|\nu^{\frac12}\mathbf P_1g_2\|_{L^2_{\mathbf x,v}}^2
\Big].
\label{est-global-a-linear}
\end{align}

By \eqref{def-H}, \eqref{est-Linf-Q}, \eqref{est-Linf-Gs}, and \eqref{apriori}, we obtain
\begin{align}
\big|\langle\mu^{-\frac12}H(g_1,g_2),\psi_a\rangle_{\mathbf x,v}\big|
&\lesssim
\|\nu^{-1}w_lH(g_1,g_2)\|_{L^\infty_{\mathbf x,v}}
\bigg(\int_{\mathbb R^3}\frac{\nu(v)(1+|v|^3)}{w_l(v)}\,dv\bigg)
\|a\|_{L^2_{\mathbf x}}
\notag\\
&\le
C_l(\alpha+\delta+\kappa)
\Big[\|h_1\|_{L^\infty_{\mathbf x,v}}+\|h_2\|_{L^\infty_{\mathbf x,v}}\Big]
\|a\|_{L^2_{\mathbf x}}
\notag\\
&\le
\varepsilon\|a\|_{L^2_{\mathbf x}}^2
+
C_{\varepsilon,l}(\alpha+\delta+\kappa)^2
\Big[\|h_1\|_{L^\infty_{\mathbf x,v}}^2+\|h_2\|_{L^\infty_{\mathbf x,v}}^2\Big].
\label{est-global-a-nonlinear}
\end{align}

Combining \eqref{eq-global-a-test}, \eqref{est-global-a-boundary},
\eqref{est-global-a-dt}, and \eqref{est-global-a-transport}--
\eqref{est-global-a-nonlinear}, taking $\varepsilon>0$
sufficiently small, we obtain
\begin{align}
\frac{d}{dt}\langle g,\psi_a\rangle_{\mathbf x,v}
+\|a\|_{L^2_{\mathbf x}}^2
\le{}&
C_l\Big[
\|\mathbf b\|_{L^2_{\mathbf x}}^2
+
\|\nu^{\frac12}\mathbf P_1g_2\|_{L^2_{\mathbf x,v}}^2
+
|(I-P_\gamma)g_2|_{L^2(\gamma_+)}^2
+
\|h_1\|_{L^\infty_{\mathbf x,v}}^2
\Big]
\notag\\
&+
C_l(\alpha+\delta)^2\|\mathcal J[g]\|_{L^2(\partial\Omega)}^2
+
C_l(\alpha+\delta+\kappa)^2\|h_2\|_{L^\infty_{\mathbf x,v}}^2.
\label{est-global-a}
\end{align}

\textbf{Estimate of $\mathbf b$.}
For $i=x,y,z$, we introduce $\phi_{b_i}=\phi_{b_i}(t,\mathbf x)$ and
$\psi_b=\psi_b(t,\mathbf x,v)$ by
\begin{equation}\label{def-global-psib}
\left\{
\begin{aligned}
-\Delta_{\mathbf x}\phi_{b_i}
&=
b_i,\quad \mathbf x\in\Omega,
\\
\phi_{b_i}\big|_{\partial\Omega}
&=
0,
\\
\psi_b(t,\mathbf x,v)
&=
\sum_{i=x,y,z}
\Big[
-2\partial_i\phi_{b_i}(t,\mathbf x)
+
v_iv\cdot\nabla_{\mathbf x}\phi_{b_i}(t,\mathbf x)
\Big]\sqrt\mu.
\end{aligned}
\right.
\end{equation}
Then,
\begin{equation}\label{est-global-phib}
\|\phi_{b_i}\|_{H^2_{\mathbf x}}
+
\|\partial_y\phi_{b_i}\|_{L^2(\partial\Omega)}
\lesssim
\|b_i\|_{L^2_{\mathbf x}},
\qquad i=x,y,z.
\end{equation}

Testing \eqref{eq-g-unsteady} with $\psi_b$, we obtain
\begin{align}
\frac{d}{dt}\langle g,\psi_b\rangle_{\mathbf x,v}
+
\|\mathbf b\|_{L^2_{\mathbf x}}^2
={}&
\langle g,\partial_t\psi_b\rangle_{\mathbf x,v}
+
\langle\mathbf P_1g,v\cdot\nabla_{\mathbf x}\psi_b\rangle_{\mathbf x,v}
\notag\\
&-
\underbrace{
\int_{\gamma_+}g(t,x,\pm1,z,v)\psi_b\,d\gamma
-
\int_{\gamma_-}g(t,x,\pm1,z,v)\psi_b\,d\gamma
}_{\mathcal B_b}
\notag\\
&-
\langle Lg,\psi_b\rangle_{\mathbf x,v}
+
\langle\mu^{-\frac12}H(g_1,g_2),\psi_b\rangle_{\mathbf x,v}.
\label{eq-global-b-test}
\end{align}

By \eqref{def-global-psib},
\[
\psi_b(t,x,\pm1,z,v)
=
\sqrt\mu\Big[
v_xv_y\partial_y\phi_{b_x}
+
(v_y^2-2)\partial_y\phi_{b_y}
+
v_zv_y\partial_y\phi_{b_z}
\Big](t,x,\pm1,z).
\]
A direct computation gives
\[
\int_{\gamma_-}P_\gamma g\,\psi_b\,d\gamma
=
\int_{\gamma_+}P_\gamma g_2\,\psi_b\,d\gamma
=
0.
\]
Together with the incoming boundary condition in \eqref{eq-g-unsteady}, these yield
\begin{align}
\mathcal B_b
={}&
\int_{\gamma_+}\mu^{-\frac12}g_1\psi_b\,d\gamma
+
\int_{\gamma_+}(I-P_\gamma)g_2\psi_b\,d\gamma
-
\sqrt{2\pi}\int_{\gamma_-}\eta^{\alpha,\delta}\mathcal J[g]\psi_b\,d\gamma.
\label{eq-global-b-bd-decomp}
\end{align}
Arguing as in \eqref{est-global-a-boundary} and using \eqref{est-global-eta-boundary} and \eqref{est-global-phib}, we obtain
\begin{align}
|\mathcal B_b|
\lesssim{}&
\Big[
\|h_1\|_{L^\infty_{\mathbf x,v}}
+
|(I-P_\gamma)g_2|_{L^2(\gamma_+)}
+
|\eta^{\alpha,\delta}\mathcal J[g]|_{L^2(\gamma_-)}
\Big]\|\mathbf b\|_{L^2_{\mathbf x}}
\notag\\
\le{}&
\varepsilon\|\mathbf b\|_{L^2_{\mathbf x}}^2
+
C_{\varepsilon,l}\Big[
\|h_1\|_{L^\infty_{\mathbf x,v}}^2
+
|(I-P_\gamma)g_2|_{L^2(\gamma_+)}^2
+
(\alpha+\delta)^2\|\mathcal J[g]\|_{L^2(\partial\Omega)}^2
\Big].
\label{est-global-b-boundary}
\end{align}

For $\langle g,\partial_t\psi_b\rangle_{\mathbf x,v}$ in
\eqref{eq-global-b-test}, we apply $\partial_t$ to the first equation
in \eqref{def-global-psib} and use the second equations in
\eqref{eq-global-macro-system} to obtain
\[
-\Delta_{\mathbf x}\partial_t\phi_{b_i}
=
-\partial_i\Big(a+\frac{2}{\sqrt6}c\Big)
-
\nabla_{\mathbf x}\cdot\int_{\mathbb R^3}v_i v\sqrt\mu \mathbf P_1g\,dv,
\qquad i=x,y,z.
\]
Testing with $\partial_t\phi_{b_i}$, we have
\begin{align}
\|\nabla_{\mathbf x}\partial_t\phi_{b_i}\|_{L^2_{\mathbf x}}^2
={}&
\langle a,\partial_i\partial_t\phi_{b_i}\rangle_{\mathbf x}
+
\frac{2}{\sqrt6}\langle c,\partial_i\partial_t\phi_{b_i}\rangle_{\mathbf x}
+
\bigg\langle
\int_{\mathbb R^3}v_i v\sqrt\mu \mathbf P_1g\,dv,
\nabla_{\mathbf x}\partial_t\phi_{b_i}
\bigg\rangle_{\mathbf x}.
\label{est-global-dt-phib}
\end{align}
By \eqref{def-global-macro-g} and \eqref{est-global-dt-phib},
\begin{align}
\langle g,\partial_t\psi_b\rangle_{\mathbf x,v}
={}&
\sum_{i=x,y,z}
\bigg\langle-a+\frac{2}{\sqrt6}c,\partial_i\partial_t\phi_{b_i}\bigg\rangle_{\mathbf x}
+
\sum_{i=x,y,z}\bigg\langle
\int_{\mathbb R^3}v_i v\sqrt\mu\mathbf P_1g\,dv,
\nabla_{\mathbf x}\partial_t\phi_{b_i}
\bigg\rangle_{\mathbf x}
\notag\\
={}&
-\sum_{i=x,y,z}\|\nabla_{\mathbf x}\partial_t\phi_{b_i}\|_{L^2_{\mathbf x}}^2
+
\frac{4}{\sqrt6}\sum_{i=x,y,z}\langle c,\partial_i\partial_t\phi_{b_i}\rangle_{\mathbf x}
+
2\sum_{i=x,y,z}\bigg\langle
\int_{\mathbb R^3}v_i v\sqrt\mu\mathbf P_1g\,dv,
\nabla_{\mathbf x}\partial_t\phi_{b_i}
\bigg\rangle_{\mathbf x}
\notag\\
\le{}&
-\sum_{i=x,y,z}\|\nabla_{\mathbf x}\partial_t\phi_{b_i}\|_{L^2_{\mathbf x}}^2
+
C\bigg[
\|c\|_{L^2_{\mathbf x}}
+
\|\nu^{\frac12}\mathbf P_1g_2\|_{L^2_{\mathbf x,v}}
+
\|h_1\|_{L^\infty_{\mathbf x,v}}
\int_{\mathbb R^3}\frac{1+|v|^2}{w_l(v)}\,dv
\bigg]
\notag\\
&\qquad\times
\sum_{i=x,y,z}\|\nabla_{\mathbf x}\partial_t\phi_{b_i}\|_{L^2_{\mathbf x}}
\notag\\
\le{}&
-\frac12\sum_{i=x,y,z}\|\nabla_{\mathbf x}\partial_t\phi_{b_i}\|_{L^2_{\mathbf x}}^2
+
C\Big[
\|c\|_{L^2_{\mathbf x}}^2
+
\|\nu^{\frac12}\mathbf P_1g_2\|_{L^2_{\mathbf x,v}}^2
+
\|h_1\|_{L^\infty_{\mathbf x,v}}^2
\Big].
\label{est-global-b-dt}
\end{align}

Arguing as in \eqref{est-global-a-transport}--\eqref{est-global-a-nonlinear}, we obtain
\begin{align}
&
\big|\langle\mathbf P_1g,v\cdot\nabla_{\mathbf x}\psi_b\rangle_{\mathbf x,v}\big|
+
\big|\langle Lg,\psi_b\rangle_{\mathbf x,v}\big|
+
\big|\langle\mu^{-\frac12}H(g_1,g_2),\psi_b\rangle_{\mathbf x,v}\big|
\notag\\
&\le
\varepsilon\|\mathbf b\|_{L^2_{\mathbf x}}^2
+
C_{\varepsilon,l}\Big[
\|\nu^{\frac12}\mathbf P_1g_2\|_{L^2_{\mathbf x,v}}^2
+
\|h_1\|_{L^\infty_{\mathbf x,v}}^2
+
(\alpha+\delta+\kappa)^2\|h_2\|_{L^\infty_{\mathbf x,v}}^2
\Big].
\label{est-global-b-bulk}
\end{align}

Combining \eqref{eq-global-b-test}, \eqref{est-global-b-boundary},
\eqref{est-global-b-dt} and \eqref{est-global-b-bulk} and taking
$\varepsilon>0$ sufficiently small, we obtain
\begin{align}
\frac{d}{dt}\langle g,\psi_b\rangle_{\mathbf x,v}
+\frac12\|\mathbf b\|_{L^2_{\mathbf x}}^2
\le{}&
C_l\Big[
\|c\|_{L^2_{\mathbf x}}^2
+
\|\nu^{\frac12}\mathbf P_1g_2\|_{L^2_{\mathbf x,v}}^2
+
|(I-P_\gamma)g_2|_{L^2(\gamma_+)}^2
+
\|h_1\|_{L^\infty_{\mathbf x,v}}^2
\Big]
\notag\\
&+
C_l(\alpha+\delta)^2\|\mathcal J[g]\|_{L^2(\partial\Omega)}^2
+
C_l(\alpha+\delta+\kappa)^2\|h_2\|_{L^\infty_{\mathbf x,v}}^2.
\label{est-global-b}
\end{align}

\textbf{Estimate of $c$.}
We next introduce $\phi_c=\phi_c(t,\mathbf x)$ and $\psi_c=\psi_c(t,\mathbf x,v)$ by
\begin{equation}\label{def-global-psic}
\left\{
\begin{aligned}
-\Delta_{\mathbf x}\phi_c
&=c,\quad \mathbf x\in\Omega,
\\
\phi_c\big|_{\partial\Omega}
&=0,
\\
\psi_c(t,\mathbf x,v)
&=(|v|^2-5)\sqrt\mu v\cdot\nabla_{\mathbf x}\phi_c(t,\mathbf x).
\end{aligned}
\right.
\end{equation}
Then,
\begin{equation}\label{est-global-phic}
\|\phi_c\|_{H^2_{\mathbf x}}
+
\|\partial_y\phi_c\|_{L^2(\partial\Omega)}
\lesssim
\|c\|_{L^2_{\mathbf x}}.
\end{equation}

Testing \eqref{eq-g-unsteady} with $\psi_c$, we obtain
\begin{align}
\frac{d}{dt}\langle g,\psi_c\rangle_{\mathbf x,v}
+
\frac{10}{\sqrt6}\|c\|_{L^2_{\mathbf x}}^2
={}&
\langle g,\partial_t\psi_c\rangle_{\mathbf x,v}
+
\langle\mathbf P_1g,v\cdot\nabla_{\mathbf x}\psi_c\rangle_{\mathbf x,v}
\notag\\
&-
\underbrace{
\int_{\gamma_+}g(t,x,\pm1,z,v)\psi_c\,d\gamma
-
\int_{\gamma_-}g(t,x,\pm1,z,v)\psi_c\,d\gamma
}_{\mathcal B_c}
\notag\\
&-
\langle Lg,\psi_c\rangle_{\mathbf x,v}
+
\langle\mu^{-\frac12}H(g_1,g_2),\psi_c\rangle_{\mathbf x,v}.
\label{eq-global-c-test}
\end{align}

By \eqref{def-global-psic},
\[
\psi_c(t,x,\pm1,z,v)
=
(|v|^2-5)\sqrt\mu v_y\partial_y\phi_c(t,x,\pm1,z).
\]
A direct computation gives
\[
\int_{\gamma_-}P_\gamma g\,\psi_c\,d\gamma
=
\int_{\gamma_+}P_\gamma g_2\,\psi_c\,d\gamma
=
0.
\]
Together with the boundary condition in \eqref{eq-g-unsteady}, we have
\[
\mathcal B_c
=
\int_{\gamma_+}\mu^{-\frac12}g_1\psi_c\,d\gamma
+
\int_{\gamma_+}(I-P_\gamma)g_2\psi_c\,d\gamma
-
\sqrt{2\pi}\int_{\gamma_-}\eta^{\alpha,\delta}\mathcal J[g]\psi_c\,d\gamma.
\]
Arguing as in \eqref{est-global-a-boundary} and using
\eqref{est-global-eta-boundary} and \eqref{est-global-phic}, we obtain
\begin{align}
|\mathcal B_c|
\le{}&
\varepsilon\|c\|_{L^2_{\mathbf x}}^2
+
C_{\varepsilon,l}\Big[
\|h_1\|_{L^\infty_{\mathbf x,v}}^2
+
|(I-P_\gamma)g_2|_{L^2(\gamma_+)}^2
+
(\alpha+\delta)^2\|\mathcal J[g]\|_{L^2(\partial\Omega)}^2
\Big].
\label{est-global-c-boundary}
\end{align}

For $\langle g,\partial_t\psi_c\rangle_{\mathbf x,v}$ in
\eqref{eq-global-c-test}, we apply $\partial_t$ to the first equation
in \eqref{def-global-psic} and use the last equation in
\eqref{eq-global-macro-system} to obtain
\[
-\Delta_{\mathbf x}\partial_t\phi_c
=
-\frac{2}{\sqrt6}\nabla_{\mathbf x}\cdot\mathbf b
-
\frac{1}{\sqrt6}\nabla_{\mathbf x}\cdot
\int_{\mathbb R^3}(|v|^2-5)v\sqrt\mu\mathbf P_1g\,dv.
\]
Testing with $\partial_t\phi_c$ and estimating the velocity moment as in
\eqref{est-global-a-dt}, we obtain
\begin{align}
\|\nabla_{\mathbf x}\partial_t\phi_c\|_{L^2_{\mathbf x}}
&\lesssim
\|\mathbf b\|_{L^2_{\mathbf x}}
+
\bigg\|\int_{\mathbb R^3}(|v|^2-5)v\sqrt\mu \mathbf P_1g\,dv\bigg\|_{L^2_{\mathbf x}}
\notag\\
&\lesssim
\|\mathbf b\|_{L^2_{\mathbf x}}
+
\|\nu^{\frac12}\mathbf P_1g_2\|_{L^2_{\mathbf x,v}}
+
\|h_1\|_{L^\infty_{\mathbf x,v}}.
\label{est-global-dt-phic}
\end{align}
By \eqref{est-global-dt-phic}, we have
\begin{align}
\big|\langle g,\partial_t\psi_c\rangle_{\mathbf x,v}\big|
={}&
\bigg|
\bigg\langle
\int_{\mathbb R^3}(|v|^2-5)v\sqrt\mu\mathbf P_1g\,dv,
\nabla_{\mathbf x}\partial_t\phi_c
\bigg\rangle_{\mathbf x}
\bigg|
\notag\\
\lesssim{}&
\Big[
\|\nu^{\frac12}\mathbf P_1g_2\|_{L^2_{\mathbf x,v}}
+
\|h_1\|_{L^\infty_{\mathbf x,v}}
\Big]\|\nabla_{\mathbf x}\partial_t\phi_c\|_{L^2_{\mathbf x}}
\notag\\
\le{}&
\varepsilon\|\mathbf b\|_{L^2_{\mathbf x}}^2
+
C_{\varepsilon,l}\Big[
\|\nu^{\frac12}\mathbf P_1g_2\|_{L^2_{\mathbf x,v}}^2
+
\|h_1\|_{L^\infty_{\mathbf x,v}}^2
\Big].
\label{est-global-c-dt}
\end{align}

Arguing as in \eqref{est-global-a-transport}--\eqref{est-global-a-nonlinear}, we obtain
\begin{align}
&
\big|\langle\mathbf P_1g,v\cdot\nabla_{\mathbf x}\psi_c\rangle_{\mathbf x,v}\big|
+
\big|\langle Lg,\psi_c\rangle_{\mathbf x,v}\big|
+
\big|\langle\mu^{-\frac12}H(g_1,g_2),\psi_c\rangle_{\mathbf x,v}\big|
\notag\\
&\le
\varepsilon\|c\|_{L^2_{\mathbf x}}^2
+
C_{\varepsilon,l}\Big[
\|\nu^{\frac12}\mathbf P_1g_2\|_{L^2_{\mathbf x,v}}^2
+
\|h_1\|_{L^\infty_{\mathbf x,v}}^2
+
(\alpha+\delta+\kappa)^2\|h_2\|_{L^\infty_{\mathbf x,v}}^2
\Big].
\label{est-global-c-bulk}
\end{align}

Combining \eqref{eq-global-c-test}, \eqref{est-global-c-boundary},
\eqref{est-global-c-dt} and \eqref{est-global-c-bulk} and taking
$\varepsilon>0$ sufficiently small, we obtain
\begin{align}
\frac{d}{dt}\langle g,\psi_c\rangle_{\mathbf x,v}
+
\|c\|_{L^2_{\mathbf x}}^2
\le{}&
\varepsilon\|\mathbf b\|_{L^2_{\mathbf x}}^2
+
C_{\varepsilon,l}\Big[
\|\nu^{\frac12}\mathbf P_1g_2\|_{L^2_{\mathbf x,v}}^2
+
|(I-P_\gamma)g_2|_{L^2(\gamma_+)}^2
+
\|h_1\|_{L^\infty_{\mathbf x,v}}^2
\Big]
\notag\\
&+
C_{\varepsilon,l}(\alpha+\delta)^2\|\mathcal J[g]\|_{L^2(\partial\Omega)}^2
+
C_{\varepsilon,l}(\alpha+\delta+\kappa)^2\|h_2\|_{L^\infty_{\mathbf x,v}}^2.
\label{est-global-c}
\end{align}

We now combine \eqref{est-global-a}, \eqref{est-global-b}, and \eqref{est-global-c}. Choose a fixed sufficiently small $\kappa_1>0$, define
\begin{equation}\label{def-global-E-int}
\mathcal E_{\rm int}(t)
:=
\kappa_1^2\langle g,\psi_a\rangle_{\mathbf x,v}
+
\kappa_1\langle g,\psi_b\rangle_{\mathbf x,v}
+
\langle g,\psi_c\rangle_{\mathbf x,v}.
\end{equation}
Taking $\kappa_1$ and then $\varepsilon$ sufficiently small, we obtain, for some $\lambda_1>0$,
\begin{align}
\frac{d}{dt}\mathcal E_{\rm int}(t)
+
\lambda_1\|[a,\mathbf b,c]\|_{L^2_{\mathbf x}}^2
\le{}&
C_l\Big[
\|\nu^{\frac12}\mathbf P_1g_2\|_{L^2_{\mathbf x,v}}^2
+
|(I-P_\gamma)g_2|_{L^2(\gamma_+)}^2
+
\|h_1\|_{L^\infty_{\mathbf x,v}}^2
\Big]
\notag\\
&+
C_l(\alpha+\delta)^2\|\mathcal J[g]\|_{L^2(\partial\Omega)}^2
+
C_l(\alpha+\delta+\kappa)^2\|h_2\|_{L^\infty_{\mathbf x,v}}^2.
\label{est-global-macro}
\end{align}
By \eqref{def-global-E-int}, \eqref{est-global-phi},
\eqref{est-global-phib}, and \eqref{est-global-phic},
\begin{equation}\label{est-global-L2-equi}
|\mathcal E_{\rm int}(t)|
\lesssim
\|g_2(t)\|_{L^2_{\mathbf x,v}}^2
+
\|h_1(t)\|_{L^\infty_{\mathbf x,v}}^2.
\end{equation}
Since $\sqrt\mu g=g_1+\sqrt\mu g_2$,
\begin{equation}\label{est-global-L2-flux}
\|\mathcal J[g](t)\|_{L^2(\partial\Omega)}^2
\lesssim
\|h_1(t)\|_{L^\infty_{\mathbf x,v}}^2
+
\|h_2(t)\|_{L^\infty_{\mathbf x,v}}^2.
\end{equation}

Adding $\kappa_2$ times \eqref{est-global-macro} to \eqref{eq-global-g2-energy}
and using \eqref{est-global-P0-g2}, \eqref{est-global-L2-equi}, \eqref{est-global-L2-flux}, we obtain, for $0<\kappa_2\ll 1$ and then
$\varepsilon$ in \eqref{eq-global-g2-energy} sufficiently small,
\begin{align}
&\frac{d}{dt}
\Big[\|g_2\|_{L^2_{\mathbf x,v}}^2+\kappa_2\mathcal E_{\rm int}(t)\Big]
+
\lambda_2\Big[
\|[a,\mathbf b,c]\|_{L^2_{\mathbf x}}^2
+
\|\nu^{\frac12}\mathbf P_1g_2\|_{L^2_{\mathbf x,v}}^2
\Big]
+
\frac12|(I-P_\gamma)g_2|_{L^2(\gamma_+)}^2
\notag\\
&\le
C_{\varepsilon_5,l}\|h_1\|_{L^\infty_{\mathbf x,v}}^2
+
C_l\big[\varepsilon_5+(\alpha+\delta+\kappa)^2\big]
\|h_2\|_{L^\infty_{\mathbf x,v}}^2
\label{est-global-L2-combined}
\end{align}
for some $\lambda_2>0$. 

\noindent\textbf{Case 1. $0\le\gamma\le1$.}
Since $\nu(v)\ge\nu_0>0$, we have
\begin{equation}\label{est-global-g2-micro-hard}
\|\mathbf P_1g_2\|_{L^2_{\mathbf x,v}}^2
\le
\nu_0^{-1}\|\nu^{\frac12}\mathbf P_1g_2\|_{L^2_{\mathbf x,v}}^2.
\end{equation}
It follows from \eqref{est-global-g2-micro-hard}, \eqref{est-global-P0-g2}, and \eqref{est-global-L2-equi} that, for $\widetilde\lambda_0>0$ sufficiently small,
\[
\widetilde\lambda_0
\Big[\|g_2\|_{L^2_{\mathbf x,v}}^2+\kappa_2\mathcal E_{\rm int}(t)\Big]
\le
\frac{\lambda_2}{2}\Big[
\|[a,\mathbf b,c]\|_{L^2_{\mathbf x}}^2
+
\|\nu^{\frac12}\mathbf P_1g_2\|_{L^2_{\mathbf x,v}}^2
\Big]
+
C_l\|h_1\|_{L^\infty_{\mathbf x,v}}^2.
\]
Combining this bound with \eqref{est-global-L2-combined}, we obtain
\begin{align}
&\frac{d}{dt}
\Big[\|g_2\|_{L^2_{\mathbf x,v}}^2+\kappa_2\mathcal E_{\rm int}(t)\Big]
+
\widetilde\lambda_0
\Big[\|g_2\|_{L^2_{\mathbf x,v}}^2+\kappa_2\mathcal E_{\rm int}(t)\Big]
\notag\\
&\le
C_{\varepsilon_5,l}\|h_1\|_{L^\infty_{\mathbf x,v}}^2
+
C_l\big[\varepsilon_5+(\alpha+\delta+\kappa)^2\big]
\|h_2\|_{L^\infty_{\mathbf x,v}}^2.
\label{est-global-L2}
\end{align}
Applying Gr\"onwall's inequality to \eqref{est-global-L2} and using \eqref{est-global-L2-equi} and $g_2(0)=0$, we obtain
\begin{align}
\|g_2(t)\|_{L^2_{\mathbf x,v}}^2
\le{}&
C_l e^{-\widetilde\lambda_0t}
\big\|w_l(F_0-F_{\alpha,\delta}^s)\big\|_{L^\infty_{\mathbf x,v}}^2
+
C_l\|h_1(t)\|_{L^\infty_{\mathbf x,v}}^2
\notag\\
&+
C_{\varepsilon_5,l}\int_0^t e^{-\widetilde\lambda_0(t-s)}
\|h_1(s)\|_{L^\infty_{\mathbf x,v}}^2\,ds
\notag\\
&+
C_l\big[\varepsilon_5+(\alpha+\delta+\kappa)^2\big]
\int_0^t e^{-\widetilde\lambda_0(t-s)}
\|h_2(s)\|_{L^\infty_{\mathbf x,v}}^2\,ds.
\label{est-global-g2-L2-decay}
\end{align}
Let $\lambda_\gamma=\lambda$ be as in the a priori assumption \eqref{apriori}, so that $2\lambda<\widetilde\lambda_0$, we multiply \eqref{est-global-g2-L2-decay} by $e^{2\lambda t}$ and take $\alpha,\delta,\kappa>0$ sufficiently small yield
\begin{align}
\sup_{0\le s\le t}e^{\lambda s}\|g_2(s)\|_{L^2_{\mathbf x,v}}
\le{}&
C_{\varepsilon_5,l}\bigg[
\big\|w_l(F_0-F_{\alpha,\delta}^s)\big\|_{L^\infty_{\mathbf x,v}}
+
\sup_{0\le s\le t}e^{\lambda s}\|h_1(s)\|_{L^\infty_{\mathbf x,v}}
\bigg]
+
C_l\varepsilon_5
\sup_{0\le s\le t}e^{\lambda s}\|h_2(s)\|_{L^\infty_{\mathbf x,v}}.
\label{est-global-g2-L2-hard}
\end{align}

\noindent\textbf{Case 2. $-3<\gamma<0$.}
The lower bound for $\nu(v)$ used in \eqref{est-global-g2-micro-hard} is now just local in velocity. For $M_5\gg1$, we have
\begin{align}
\|g_2\|_{L^2_{\mathbf x,v}}^2
={}&
\int_\Omega\int_{|v|\le M_5}|g_2|^2\,dv d\mathbf x
+
\int_\Omega\int_{|v|>M_5}w_l^{-2}(v)|h_2|^2\,dv d\mathbf x
\notag\\
\lesssim{}&
\frac{\displaystyle\int_\Omega\int_{|v|\le M_5}\nu(v)|g_2|^2\,dv d\mathbf x}
{\displaystyle\inf_{|v|\le M_5}\nu(v)}
+
\|h_2\|_{L^\infty_{\mathbf x,v}}^2
\int_{M_5}^\infty r^{2-4l}\,dr
\notag\\
\le{}&
\frac{\|\nu^{\frac12}g_2\|_{L^2_{\mathbf x,v}}^2}
{\inf_{|v|\le M_5}\nu(v)}
+
\frac{C_l}{M_5^2}\|h_2\|_{L^\infty_{\mathbf x,v}}^2.
\label{est-global-g2-velocity-split}
\end{align}

By \eqref{est-global-P0-g2}, \eqref{est-global-L2-equi} and \eqref{est-global-g2-velocity-split}, for each $M_5$ we can choose $\lambda_{M_5}>0$ sufficiently small so that
\begin{align}
&\lambda_{M_5}\Big[\|g_2\|_{L^2_{\mathbf x,v}}^2+\kappa_2\mathcal E_{\rm int}(t)\Big]
\notag\\
&\lesssim
\frac{\lambda_{M_5}}{\inf_{|v|\le M_5}\nu(v)}
\|\nu^{\frac12}g_2\|_{L^2_{\mathbf x,v}}^2
+
\lambda_{M_5}\|h_1\|_{L^\infty_{\mathbf x,v}}^2
+
\frac{\lambda_{M_5}}{M_5^2}\|h_2\|_{L^\infty_{\mathbf x,v}}^2
\notag\\
&\le
\frac{\lambda_2}{2}\Big[
\|[a,\mathbf b,c]\|_{L^2_{\mathbf x}}^2
+
\|\nu^{\frac12}\mathbf P_1g_2\|_{L^2_{\mathbf x,v}}^2
\Big]
+
C_l\|h_1\|_{L^\infty_{\mathbf x,v}}^2
+
\frac{C_l\lambda_{M_5}}{M_5^2}\|h_2\|_{L^\infty_{\mathbf x,v}}^2.
\label{est-global-L2-soft-control}
\end{align}
Combining \eqref{est-global-L2-soft-control} with \eqref{est-global-L2-combined}, we obtain
\begin{align}
&\frac{d}{dt}
\Big[\|g_2\|_{L^2_{\mathbf x,v}}^2+\kappa_2\mathcal E_{\rm int}(t)\Big]
+
\lambda_{M_5}
\Big[\|g_2\|_{L^2_{\mathbf x,v}}^2+\kappa_2\mathcal E_{\rm int}(t)\Big]
\notag\\
&\le
C_{\varepsilon_5,l}\|h_1\|_{L^\infty_{\mathbf x,v}}^2
+
C_l\bigg[\frac{\lambda_{M_5}}{M_5^2}+\varepsilon_5+(\alpha+\delta+\kappa)^2\bigg]
\|h_2\|_{L^\infty_{\mathbf x,v}}^2.
\label{est-global-L2-soft}
\end{align}
Choosing $M_5$ sufficiently large and then $\varepsilon_5,\alpha,\delta,\kappa$ sufficiently small, we obtain, with the resulting small parameter again denoted by $\varepsilon_5$,
\begin{align}
\sup_{0\le s\le t}\|g_2(s)\|_{L^2_{\mathbf x,v}}
\le{}&
C_l
\big\|w_l(F_0-F_{\alpha,\delta}^s)\big\|_{L^\infty_{\mathbf x,v}}
+
C_{\varepsilon_5,l}\sup_{0\le s\le t}\|h_1(s)\|_{L^\infty_{\mathbf x,v}}
+
C_l\varepsilon_5\sup_{0\le s\le t}\|h_2(s)\|_{L^\infty_{\mathbf x,v}}.
\label{est-global-g2-L2-soft}
\end{align}
Since $\lambda_\gamma=0$ for $-3<\gamma<0$, \eqref{est-global-g2-L2-hard} and \eqref{est-global-g2-L2-soft} yield \eqref{est-global-g2-L2-weight}. This completes the proof of Lemma~\ref{lem-global-L2}.
\end{proof}

\begin{proof}[\itshape\bfseries Proof of Theorem~\ref{thm-2}]
Combining \eqref{est-global-g2-L2-weight} in Lemma~\ref{lem-global-L2} with \eqref{est-global-h1} and \eqref{est-global-h2} in Lemma~\ref{lem-global-Linf}, and taking $\varepsilon_5,\alpha,\delta,\kappa$ sufficiently small, we obtain
\begin{equation}\label{est-global-close}
\sup_{0\le s\le t}e^{\lambda_\gamma s}
\Big[\|h_1(s)\|_{L^\infty_{\mathbf x,v}}+\|h_2(s)\|_{L^\infty_{\mathbf x,v}}\Big]
\le
C_l\big\|w_l(F_0-F_{\alpha,\delta}^s)\big\|_{L^\infty_{\mathbf x,v}}.
\end{equation}
Taking the initial perturbation sufficiently small, we improve the a priori assumption \eqref{apriori}. By Lemma~\ref{lem-local-existence} and the standard continuation argument, the solution extends uniquely to all times and satisfies \eqref{est-global-decay} for $0\le\gamma\le1$ and \eqref{est-global-bound-soft} for $-3<\gamma<0$. The nonnegativity of $F$ follows from the standard iteration along characteristics based on the decomposition $Q=Q_+-Q_-$, we omit the details and refer to \cite[Section~5]{DL21} and \cite[Section~5]{Guo10}.

It remains to prove \eqref{est-global-decay-soft}. For $-3<\gamma<0$, let $l\ge l_0+|\gamma|/2$, where $l_0$ is given in Theorem~\ref{thm-2}. Motivated by the soft-potential decay estimates in \cite{DLL26}, we seek a time weight $W(t)$ satisfying
\begin{equation}\label{def-soft-W}
W(0)=1,\qquad W'(t)\ge0.
\end{equation}
For $0<m\le l-l_0$ and $i=1,2$, interpolation gives
\begin{equation}\label{est-soft-interp}
\|w_{l-m}g_i(t)\|_{L^\infty_{\mathbf x,v}}
\le
\|w_{l_0}g_i(t)\|_{L^\infty_{\mathbf x,v}}^{\frac{m}{l-l_0}}
\|w_lg_i(t)\|_{L^\infty_{\mathbf x,v}}^{1-\frac{m}{l-l_0}}.
\end{equation}
Together with \eqref{est-global-close}, this shows that it suffices to establish decay for $w_{l_0}g_i$. Multiplying \eqref{eq-g1-unsteady} and \eqref{eq-g2-unsteady} by $W(t)$, we obtain
\begin{equation}\label{eq-soft-Ag}
\left\{
\begin{aligned}
&\partial_t[W(t)g_1]+v\cdot\nabla_{\mathbf x}[W(t)g_1]+\nu(v)W(t)g_1
=\chi_M\mathcal K[W(t)g_1]+W(t)H(g_1,g_2)+W'(t)g_1,
\\
&\partial_t[W(t)g_2]+v\cdot\nabla_{\mathbf x}[W(t)g_2]+L[W(t)g_2]
=(1-\chi_M)\mu^{-\frac12}\mathcal K[W(t)g_1]+W'(t)g_2,
\\
&[W(t)g_1]\big|_{\gamma_-}=0,\qquad
[W(t)g_2]\big|_{\gamma_-}
=P_\gamma[W(t)g]+\sqrt{2\pi}\eta^{\alpha,\delta}\mathcal J[W(t)g],
\\
&W(0)g_1(0,\mathbf x,v)=F_0(\mathbf x,v)-F_{\alpha,\delta}^s(\mathbf x,v),
\qquad W(0)g_2(0,\mathbf x,v)=0,
\end{aligned}
\right.
\end{equation}
where $\eta^{\alpha,\delta}$ is defined in \eqref{def-global-eta}.

Applying Duhamel's formula to \eqref{eq-soft-Ag}, we obtain
\begin{align}
W(t)g_1(t,\mathbf x,v)
={}&
\mathbf1_{\{t_b\ge t\}}e^{-\nu(v)t}
(F_0-F_{\alpha,\delta}^s)(\mathbf X(0),v)
\notag\\
&+
\int_{\max\{0,t-t_b\}}^t e^{-\nu(v)(t-s)}
\Big\{\chi_M\mathcal K[W(s)g_1]+W(s)H(g_1,g_2)\Big\}
(s,\mathbf X(s),v)\,ds
\notag\\
&+
\int_{\max\{0,t-t_b\}}^t e^{-\nu(v)(t-s)}
W'(s)g_1(s,\mathbf X(s),v)\,ds,
\notag\\
W(t)g_2(t,\mathbf x,v)
={}&
\mathbf1_{\{t_b<t\}}e^{-\nu(v)t_b}W(t-t_b)
\Big\{P_\gamma g+\sqrt{2\pi}\eta^{\alpha,\delta}\mathcal J[g]\Big\}
(t-t_b,\mathbf X(t-t_b),v)
\notag\\
&+
\int_{\max\{0,t-t_b\}}^t e^{-\nu(v)(t-s)}
\Big\{K[W(s)g_2]+(1-\chi_M)\mu^{-\frac12}\mathcal K[W(s)g_1]\Big\}
(s,\mathbf X(s),v)\,ds
\notag\\
&+
\int_{\max\{0,t-t_b\}}^t e^{-\nu(v)(t-s)}
W'(s)g_2(s,\mathbf X(s),v)\,ds.
\label{eq-soft-Duhamel}
\end{align}
Arguing as in the proofs of \eqref{est-global-h1} and \eqref{est-global-h2}, we obtain
\begin{align}
\sup_{0\le s\le t}W(s)\|w_{l_0}g_1(s)\|_{L^\infty_{\mathbf x,v}}
\le{}&
C_l\|w_l(F_0-F_{\alpha,\delta}^s)\|_{L^\infty_{\mathbf x,v}}
+
C_l(\alpha+\delta+\kappa)
\sup_{0\le s\le t}W(s)\|w_{l_0}g_2(s)\|_{L^\infty_{\mathbf x,v}}
\notag\\
&+
C_l\sup_{0\le s\le t}
W'(s)\|\nu^{-1}w_{l_0}g_1(s)\|_{L^\infty_{\mathbf x,v}},
\notag\\
\sup_{0\le s\le t}W(s)\|w_{l_0}g_2(s)\|_{L^\infty_{\mathbf x,v}}
\le{}&
C_l\|w_l(F_0-F_{\alpha,\delta}^s)\|_{L^\infty_{\mathbf x,v}}
+
C_l\sup_{0\le s\le t}W(s)\|w_{l_0}g_1(s)\|_{L^\infty_{\mathbf x,v}}
+
C_l\sup_{0\le s\le t}W(s)\|g_2(s)\|_{L^2_{\mathbf x,v}}
\notag\\
&+
C_l\sup_{0\le s\le t}
W'(s)\Big[
\|\nu^{-1}w_{l_0}g_1(s)\|_{L^\infty_{\mathbf x,v}}
+
\|\nu^{-1}w_{l_0}g_2(s)\|_{L^\infty_{\mathbf x,v}}
\Big].
\label{est-soft-decay-Linf}
\end{align}

We estimate the terms involving $W'(t)$ by interpolation. Using \eqref{est-global-close} with $\lambda_\gamma=0$, we obtain
\begin{align}
W'(s)\|\nu^{-1}w_{l_0}g_i(s)\|_{L^\infty_{\mathbf x,v}}
&\lesssim
W'(s)\big\|(1+|v|^2)^{l_0+\frac{|\gamma|}{2}}g_i(s)\big\|_{L^\infty_{\mathbf x,v}}
\notag\\
&\le
W'(s)W(s)^{\beta-1}
\big\|(1+|v|^2)^lg_i(s)\big\|_{L^\infty_{\mathbf x,v}}^\beta
\Big[W(s)\big\|(1+|v|^2)^{l_0}g_i(s)\big\|_{L^\infty_{\mathbf x,v}}\Big]^{1-\beta}
\notag\\
&=
\frac1\beta\frac{d}{ds}[W(s)^\beta]
\|w_lg_i(s)\|_{L^\infty_{\mathbf x,v}}^\beta
\Big[W(s)\|w_{l_0}g_i(s)\|_{L^\infty_{\mathbf x,v}}\Big]^{1-\beta}
\notag\\
&\le
C_l\frac{d}{ds}[W(s)^\beta]
\|w_l(F_0-F_{\alpha,\delta}^s)\|_{L^\infty_{\mathbf x,v}}^\beta
\Big[W(s)\|w_{l_0}g_i(s)\|_{L^\infty_{\mathbf x,v}}\Big]^{1-\beta}
\notag\\
&\le
C_l\frac{d}{ds}[W(s)^\beta]
\Big[
\|w_l(F_0-F_{\alpha,\delta}^s)\|_{L^\infty_{\mathbf x,v}}
+
W(s)\|w_{l_0}g_i(s)\|_{L^\infty_{\mathbf x,v}}
\Big],\qquad i=1,2
\label{est-soft-decay-source}
\end{align}
where $\beta$ is chosen so that
\begin{equation}\label{def-soft-decay-beta}
\beta l+(1-\beta)l_0=l_0+\frac{|\gamma|}{2},
\qquad
\beta=\frac{|\gamma|}{2(l-l_0)}\in(0,1].
\end{equation}
To absorb the resulting terms in \eqref{est-soft-decay-Linf}, we require
\begin{equation}\label{est-soft-W-source}
\frac{d}{dt}[W(t)^\beta]=\varepsilon.
\end{equation}

To estimate $W(t)\|g_2(t)\|_{L^2_{\mathbf x,v}}$ using \eqref{est-global-L2-soft} with $l=l_0$, we impose
\begin{equation}\label{est-soft-W-L2}
\frac{W'(t)}{W(t)}\le\frac{\lambda_{M_5}}4,
\end{equation}
which gives
\[
\int_0^t e^{-\lambda_{M_5}(t-s)}
\left(\frac{W(t)}{W(s)}\right)^2\,ds
\le
\int_0^t e^{-\frac{\lambda_{M_5}}2(t-s)}\,ds
\le\frac{2}{\lambda_{M_5}}.
\]
Arguing as in \eqref{est-global-g2-L2-soft} and using \eqref{est-global-L2-equi} and $g_2(0)=0$, we obtain
\begin{align}
\sup_{0\le s\le t}W(s)\|g_2(s)\|_{L^2_{\mathbf x,v}}
\le{}&
C_l\|w_l(F_0-F_{\alpha,\delta}^s)\|_{L^\infty_{\mathbf x,v}}
+
C_{\varepsilon_5,l}
\sup_{0\le s\le t}W(s)\|w_{l_0}g_1(s)\|_{L^\infty_{\mathbf x,v}}
\notag\\
&+
C_l\varepsilon_5
\sup_{0\le s\le t}W(s)\|w_{l_0}g_2(s)\|_{L^\infty_{\mathbf x,v}}.
\label{est-soft-decay-L2}
\end{align}

Substituting \eqref{est-soft-decay-L2} into \eqref{est-soft-decay-Linf}, using \eqref{est-soft-decay-source} and \eqref{est-soft-W-source}, and taking $\varepsilon_5,\varepsilon$ sufficiently small, we obtain
\begin{align}
\sup_{0\le s\le t}W(s)\|w_{l_0}g_2(s)\|_{L^\infty_{\mathbf x,v}}
\le{}&
C_{\varepsilon_5,l}
\Big[
\|w_l(F_0-F_{\alpha,\delta}^s)\|_{L^\infty_{\mathbf x,v}}
+
\sup_{0\le s\le t}W(s)\|w_{l_0}g_1(s)\|_{L^\infty_{\mathbf x,v}}
\Big]
\notag\\
\le{}&
C_{\varepsilon_5,l}\|w_l(F_0-F_{\alpha,\delta}^s)\|_{L^\infty_{\mathbf x,v}}
+
C_{\varepsilon_5,l}(\alpha+\delta+\kappa)
\sup_{0\le s\le t}W(s)\|w_{l_0}g_2(s)\|_{L^\infty_{\mathbf x,v}}.
\label{est-soft-decay-close}
\end{align}
Together with the first inequality in \eqref{est-soft-decay-Linf}, this yields
\begin{equation}\label{est-soft-decay-absorb}
W(t)\Big[
\|w_{l_0}g_1(t)\|_{L^\infty_{\mathbf x,v}}
+
\|w_{l_0}g_2(t)\|_{L^\infty_{\mathbf x,v}}
\Big]
\le
C_l\|w_l(F_0-F_{\alpha,\delta}^s)\|_{L^\infty_{\mathbf x,v}}.
\end{equation}

Solving for $W(t)$ under \eqref{def-soft-W}, \eqref{est-soft-W-source}, and \eqref{est-soft-W-L2}, we may take
\[
W(t)^\beta=1+\varepsilon t,
\qquad
W(t)=(1+\varepsilon t)^{\frac{2(l-l_0)}{|\gamma|}}.
\]
Finally, for $0<m\le l-l_0$, \eqref{est-soft-interp}, \eqref{est-global-close}, and \eqref{est-soft-decay-absorb} yield
\begin{align}
\|w_{l-m}g_1(t)\|_{L^\infty_{\mathbf x,v}}
+
\|w_{l-m}g_2(t)\|_{L^\infty_{\mathbf x,v}}
&\le
C_lW(t)^{-\frac{m}{l-l_0}}
\|w_l(F_0-F_{\alpha,\delta}^s)\|_{L^\infty_{\mathbf x,v}}
\notag\\
&\le
C_l(1+t)^{-\frac{2m}{|\gamma|}}
\|w_l(F_0-F_{\alpha,\delta}^s)\|_{L^\infty_{\mathbf x,v}}.
\label{est-soft-final}
\end{align}
This proves \eqref{est-global-decay-soft} and completes the proof of Theorem~\ref{thm-2}.
\end{proof}

\noindent\textbf{Acknowledgements.}
The authors are supported by the special foundation for Guangxi Ba Gui Scholars and the National Natural Science Foundation of China (No.~12171104).

\medskip
\noindent\textbf{Data availability:}
This article does not have any external supporting data.

\end{document}